\documentclass[12pt]{report}

\usepackage[german, american]{babel}  %Deutsch

\usepackage{times}
\usepackage{graphicx}
\usepackage{amsfonts}       %
\usepackage{amsmath} %[leqno]=Gleichungen werden links numeriert;
\usepackage{amssymb}        %
\usepackage{amstext}        %
\usepackage{amsthm}         %theorem-umgebungen

\usepackage{epsfig}
\psfigdriver{dvips}
\usepackage{latexsym}
\usepackage[matrix,curve]{xypic}
\usepackage{rotating}
\rotdriver{dvips}

\newcommand{\Proj}{\mathbb{P}}

\newcommand{\N}{\mathbb{N}}
\newcommand{\Z}{\mathbb{Z}}
\newcommand{\Q}{\mathbb{Q}}

\newcommand{\C}{\mathbb{C}}

\newtheorem{Lemma}{Lemma}
\newtheorem{Theorem}{Theorem}
\newtheorem{Corollary}{Corollary}

\DeclareMathSizes{12}{11}{8}{5} 
\usepackage{xypic,dsfont}

\begin{document}
\author{Rainer Weissauer}
\title{Brill-Noether Sheaves}
\maketitle

\tableofcontents

\section{Introduction}
Let $X$ be an abelian variety $X$ of dimension $g$ over an algebraically closed field $k$. For
irreducible subvarieties $Y\subseteq X$ the  perverse intersection cohomology sheaves
$\delta_Y$ define perverse sheaves on $X$.

\bigskip\noindent
Let $D_c^b(X)=D_c^b(X,\overline\Q_l)$ be the triangulated category  of etale
$\overline\Q_l$-sheaves on $X$ in the sense of \cite{KW}. It contains the abelian category
$Perv(X)$ of (middle) perverse sheaves on $X$ as a full abelian subcategory. The addition law
of the abelian variety $a:X\times X\to X$ defines the convolution product $K*L\in
D_c^b(X,\overline\Q_l)$ of two complexes $K$ and $L$ in $D_c^b(X,\overline\Q_l)$. Similarly one
can define the iterated products $\delta_{Y_1}* \cdots *\delta_{Y_r}$ for subvarieties
$Y_1,..,Y_r$ of $X$. These convolutions tend to become complicated, if the sum  of the
dimensions $dim(Y_1) + \cdots + dim(Y_r) $ exceeds $g$.

\bigskip\noindent
In the simplest case of $r$ smooth curves $Y_i$ in $X$, where $r-1$ of them generate the
abelian variety $X$, we show, that the convolution of $r$ perverse sheaves $\delta_{Y_i}$ is a
perverse sheaf on $X$ for $r\leq g$. We later extend this result.

\bigskip\noindent
If $X$ is the Jacobian of a curve $C$ of genus $g\geq 2$, and if all $Y_i$ coincide and are
curves, the above assertion remains true for $r>g$ in the following weaker sense: Each $r$-fold
iterated self convolution of $\delta_C$ is a direct sum of a perverse sheaf and a complex $T$,
such that $T$ is a direct sum of translates of constant sheaves on $X$. For this we decompose
the $r$-fold iterated self convolutions of $\delta_C$ into direct summands $\delta_\alpha$
parameterized by partitions $\alpha$ of degree $deg(\alpha)=r$. These $\delta_\alpha$ behave
nicely with respect to the convolution product on $X$. We prove, that $\delta_\alpha$ is a
direct sum of a perverse sheaf ${}^p\delta_\alpha$ and a complex $T_\alpha$, such that
$T_\alpha$ is a translate of constant sheaves on $X$, and we  compute $T_\alpha$. The sheaves
${}^p\delta_\alpha$ are related to the Brill-Noether varieties in $X$ in a natural way via the
supports of the cohomology sheaves. The semisimple category ${\cal BN}'$, defined by direct
sums of perverse constituents of the perverse sheaves ${}^p\delta_\alpha$, is a tensor
category, which is a Tannakian category ${\cal BN}$ equivalent for $g\geq 3$ to the category
$Rep(G)$ of representations of an algebraic group $G$ over $\overline\Q_l$. This group  either
is $Sp(2g-2,\overline\Q_l)$ or $Sl(2g-2,\overline\Q_l)$, depending on whether $C$ is
hyperelliptic or not, assuming that the Riemann constant has been normalized suitably. In
particular the ${}^p\delta_\alpha$ are irreducible perverse sheaves in the latter case. The
category ${\cal BN}$ is related to more general tensor categories ${\cal T}(X)$ defined for
arbitrary abelian varieties $X$ in case $k$ has characteristic zero or is the algebraic closure
of a finite field.

\bigskip\noindent
The Tannakian approach in chapter \ref{Tann} implies most results of chapter \ref{Sbr} by an
independent argument. However, since the arguments of chapter \ref{Sbr} are more elementary and
hold over arbitrary algebraically closed fields $k$, I thought it might be useful also to
present the more elementary proofs. See also \cite{W}.

\chapter{Preliminaries on perverse sheaves} For  a variety $X$
over an algebraically closed field $k$ and  a prime $l$ different
from the characteristic $p$ of $k$ let $D_c^b(X,\overline\Q_l)$
denote the triangulated category of complexes of etale $\overline
\Q_l$-sheaves on $X$ in the sense of \cite{KW}. This category
carries a standard $t$-structure. Truncation with respect to the
standard $t$-structure defines the etale $\overline\Q_l$-sheaves on
$X$, also called $l$-adic sheaves for simplicity. For a complex
$K\in D_c^b(X,\overline\Q_l)$ let ${\cal H}^\nu(K)$ denote its etale
cohomology $\overline \Q_l$-sheaves with respect to the standard
$t$-structure. Let
$$Perv(X)$$ denote the abelian subcategory of middle perverse
sheaves of the triangulated category $D_c^b(X,\overline\Q_l)$.  Recall $K\in Perv(X)$ if and
only if the complex $K$ and its Verdier dual $DK$ are contained in the full subcategory
${}^p\!D^{\leq 0}(X)$ of semi-perverse sheaves in $D_c^b(X,\overline\Q_l)$, where $K\in
D_c^b(X,\overline\Q_l)$ is semi-perverse if and only if  $dim(S_\nu) \leq \nu$ holds for all
integers $\nu \in \Z$, where $S_\nu$ denotes the support \label{supports} of the cohomology
sheaf ${\cal H}^{-\nu}(K)$ of $K$. Here, by convention, the support $S_\nu$ is defined to be
the Zariski closure of the locus of points $x$ for which the stalk cohomology ${\cal
H}^{-\nu}(K)_{\overline x}$ of a geometric point $\overline x$ over $x$ does not vanish. For a
complex $K\in D_c^b(X,\overline\Q_l)$ let $D(K)$ denote its Verdier dual. We say $K$ is
selfdual, if $K\cong D(K)$ holds.

\bigskip\noindent
If $k$ is the algebraic closure of a finite field $\kappa$, then a
complex $K$ of etale $\overline\Q_l$-Weil sheaves is called mixed of
weight $\leq w$, if all its cohomology sheaves ${\cal H}^\nu(K)$ are
mixed etale $\overline\Q_l$-sheaves with upper weights $w({\cal
H}^\nu(K))-\nu \leq w$ for all integers $\nu$. It is called pure of
weight $w$, if $K$ and its Verdier dual $DK$ are mixed of weight
$\leq w$. In this case let $Perv_m(X)$ denote the abelian category
of mixed perverse sheaves on $X$. Concerning base fields of
characteristic zero, the following should be remarked. Since all
sheaves relevant for this paper will be of geometric origin in the
sense of the last chapter of \cite{BBD}, we still dispose over the
notion of the weight filtration and purity and Gabber's
decomposition theorem even in case $char(k)=0$. Hoping that this
will not give too much confusion, we therefore will not make
distinctions between the case where the base field is of
characteristic zero or not. However it should be remarked, that all
results (e.g. section \ref{categor}) which use the curve lemma
(lemma \ref{curvel}) are only proven under the additional condition,
that either $k$ has characteristic zero or is the algebraic closure
of a finite field $\kappa$.

\bigskip\noindent
If $i:Y\hookrightarrow X$ is an irreducible subvariety of $X$, and if $E$ is a smooth $l$-adic
sheaf defined on an essentially smooth dense open subset $j:U\hookrightarrow Y$ of $Y$, then
the sheaf  defined by the intermediate extension $j_{!*}$
$$\delta_E= i_*j_{!*}(E[dim(X)])$$ is a perverse sheaf on $X$. $\delta_E$ is
an irreducible perverse sheaf on $X$ if and only if $E$ is an irreducible smooth
$\overline\Q_l$-sheaf on $U$. Every irreducible perverse sheaf is isomorphic to some $\delta_E$
for some triple $(E,U,Y)$ with the notations above. $U$ is not uniquely determined by
$\delta_E$, whereas the closure $Y=Y_E=\overline U$ only depends on $\delta_E$. If $E$ is the
constant sheaf $\overline\Q_{l,U}$ on $U$, we usually write $\delta_Y$ instead of $\delta_E$.
In this case $\delta_Y$ is a selfdual perverse sheaf on $X$. If $Y$ is irreducible, then
$\delta_Y$ is pure of weight $dim(Y)$. Finally, we write $\lambda_E$ respectively $\lambda_Y$
for the etale sheaf $i_*(j_*(E))$ on $X$ with support in $Y$. If $d$ is the dimension of $Y$ we
write $H^i(Y)$ for the ordinary etale cohomology group $H^{i}(Y,\overline\Q_{l,Y})$, and
$IH^i(Y)$ for the intersection cohomology group $H^{i-d}(Y,\delta_Y)$.

\bigskip\noindent
\begin{Lemma}\label{normal} {\it If $Y$ is a normal variety of dimension $d$, then $\lambda_Y=
\overline\Q_{l,Y}[d]$.} \end{Lemma}

\bigskip\noindent
\underbar{Proof}: Let $U$ be open dense smooth in $Y$ and let $j:U\hookrightarrow Y$ be the
inclusion map. Since $\lambda_Y={\cal H}^{-d}(\delta_Y)$ is $j_*\overline\Q_{l,Y}$, there
exists a natural morphism
$$ \overline\Q_{l,Y} \longrightarrow
j_*(\overline\Q_{l,U}) =\lambda_Y  \ .$$ By assumption $j$ is dominant and $Y$ is normal. Hence
\cite{SGA} IX, lemma 2.14.1 implies, that this sheaf homomorphism is an isomorphism.

\bigskip\noindent
\begin{Lemma} {\it For irreducible perverse sheaves $\delta_E$  there exists a distinguished
triangle $\psi_E \to \lambda_E \to \delta_E \to \psi_E[1]$, such that
$\psi_E,\lambda_E,\delta_E\in {}^p D^{\leq 0}(X)$ and \begin{enumerate}
\item $0\to {}^pH^0(\psi_E) \to {}^pH^0(\lambda_E) \to \delta_E \to
0$ is exact in $Perv(X)$. \item ${}^p H^{-\nu}(X,\psi_E)=0$ for $\nu\notin
\{0,1,..,dim(Y)-2\}.$
\item $dim \ supp({}^p H^{-\nu}(X,\psi_E)) \leq dim(Y)-\nu-2$ for all $\nu\geq 0$.
\item ${\cal H}^\mu ({}^p H^{-\nu}(X,\psi_E))=0$ for all $\mu <
-dim(Y)+\nu+2$.
\item ${}^p
H^\nu(X,\psi_E) = {}^p H^\nu(X,\lambda_E)$ for all $\nu\neq 0$.
\item If $\delta_E$ is pure of weight $w$, then $w(\lambda_E)\leq
w$ and $ w(\psi_E)\leq w-1$ $$dim\ supp(\psi_E) \leq dim(Y)-2\ ,$$ and the sequence in 1.
defines the first step of the weight filtration of ${}^p H^0(\lambda_E)$.
\end{enumerate}}
\end{Lemma}

\bigskip\noindent

\bigskip\noindent
\underbar{Proof}: This can be easily deduced from [KW], section III.5 by induction. In fact for
a complex $K\in {}^p D^{\leq 0}(X)$ the conditions $K\in {}^{st} D^{\geq -r}(X)$ and
${}^pH(K)\in {}^{st} D^{\geq -r+\nu}(X)$ for all $\nu$ are equivalent.

\bigskip\noindent
Notice $${\cal H}^{-\nu}(\delta_E) \cong {\cal H}^{-\nu+1}(\psi_E) \ $$ for $\nu\neq dim(Y)$.
\qed

\bigskip\noindent
Now assume $X$ to be a scheme of finite type over the algebraic closure $k$ of a finite field
with the following properties \begin{enumerate}
\item There exists an integer $d$ such that
$\overline\Q_{l,X}[d]$ is an irreducible perverse sheaf on $X$.
\item There exists a finite morphism $\pi: X\to Y$ and a finite
group $G$ acting on $X$, such that $\pi\circ g=\pi$ holds for all $g\in G$.
\item There exists an essentially smooth open dense subset $j:U\hookrightarrow Y$,
such that $j_V:V=\pi^{-1}(U)\hookrightarrow X$ is open an dense in $X$, and such that the
restriction of $\pi$ from $V\to U$ is an etale Galois covering with Galois group $G$.
\end{enumerate}

\bigskip\noindent
\underbar{Remark}: Suppose $X$ is an irreducible local complete intersection of dimension $d$,
then $\overline\Q_{l,X}[d]$ is a perverse sheaf. Then there exists a distinguished triangle
$(\psi_X,\lambda_X,\delta_X)$ with perverse $\psi_X$, in other words an exact sequence of
perverse sheaves
$$ 0 \to \psi_X \to \lambda_X \to \delta_X \to 0 \ ,$$
defining the highest part of the weight filtration of $\lambda_X$, and the following conditions
are equivalent: a) $\lambda_X=\overline\Q_{l,X}$ is irreducible b) $\lambda_X$ is isomorphic to
the intersection cohomology sheaf $\delta_X$ c) $\lambda_X$ is irreducible d) $\lambda_X$ is
pure of weight $d$.

\bigskip\noindent
{\bf Lemma 1}: {\it Under the assumptions 1-3) above the direct image complex $$ K=
\pi_*(\overline \Q_{l,X}[d]) $$ is a $G$-equivariant perverse sheaf on $Y$, which decomposes
into a direct sum of perverse sheaves
$$ K \ =\ \bigoplus_\rho \rho \boxtimes F_\rho \ ,$$
where the $F_\rho$ are irreducible perverse sheaves on $Y$, such that $$F_\rho[-d]$$ are etale
$\overline\Q_l$-sheaves on $Y$. Furthermore $F_{\rho_1}\cong F_{\rho_2}$ if and only if
$\rho_1\cong\rho_2$.}

\bigskip\noindent
{\bf Corollary}: {\it If $X$ is smooth irreducible of dimension $d$ and $\pi: X\to Y$ is a
ramified covering with Galois group $G$. Then $\pi_*(\overline\Q_{l,X}[d]) =\bigoplus_\rho
\rho\boxtimes F_\rho$ for irreducible perverse sheaves $F_\rho\in Perv(Y)$.}

\bigskip\noindent
\underbar{Proof}: $\overline\Q_{l,X}$ is irreducible perverse by assumption, hence also the
Verdier dual $D(\overline\Q_{l,X})$. Since $D(\overline\Q_{l,X})$ and $\overline\Q_{l,X}$ are
isomorphic over $V$, they are isomorphic over $X$ being irreducible perverse sheaves. Being
irreducible $\overline\Q_{l,X}$ is pure of weight $d$.

\bigskip\noindent
$\pi$ is finite and $\overline\Q_{l,X}$ is a perverse sheaf. Thus $K=\pi_*(\overline\Q_l)$ is
in ${}^p D^{\leq 0}(X)$. By duality therefore $K\in Perv(Y)$. Since $\pi$ is proper, $K$ is
pure of weight $d$. By the equivariance of $\pi$ we can decompose $K=\bigoplus_\rho
\rho\boxtimes F_\rho$ into isotypic components with respect to the irreducible
$\overline\Q_l$-representations of $G$. The $F_\rho$ are pure perverse sheaves on $Y$. The
restrictions $F_\rho\vert_U$ to $U$   are irreducible smooth etale sheaves ([KW] III.15.3 d))
up to a complex shift by $d$. Therefore by purity
$$ F_\rho = j_{!*}(F_\rho\vert_U) \oplus R_\rho $$
for some perverse sheaf $R_\rho$ supported in $Y\setminus U$. In particular $F_\rho\neq 0$ for
all $\rho$. Recall the sheaves $F_\rho\vert_U$ are smooth etale $\overline\Q_l$-sheaves of rank
$deg(\rho)$ (by [KW] III.15.3 d)). The adjunction formula ([KW], p.107) implies
$$ Hom_{D_c^b(X,\overline\Q_l)}\Bigl(\pi^*\pi_*\overline\Q_{l,X}[d],
\overline\Q_{l,X}[d]\Bigr) \cong Hom_{D_c^b(Y,\overline\Q_l)}\Bigl(\pi_*\overline\Q_{l,X}[d],
\pi_*\overline\Q_{l,X}[d]\Bigr) \ .$$ The dimension on the right is $\sum_{\rho_1,\rho_2}
dim(\rho_1)dim(\rho_2) \cdot dim(Hom_{D_c^b(Y,\overline\Q_l)}( F_{\rho_1}, F_{\rho_2}))$. Since
$dim(Hom_{D_c^b(Y,\overline\Q_l)}( F_\rho, F_\rho)) \geq 1$ as already seen, the claims would
follow from
$$ dim(Hom_{D_c^b(X,\overline\Q_l)}\Bigl(\pi^*\pi_*\overline\Q_{l,X}[d],
\overline\Q_{l,X}[d])\Bigr) \ \leq\  \vert G\vert $$ using the character formula $\vert G\vert
=\sum_\rho dim(\rho)^2$ and the decomposition theorem (purity of $K$). We remark, that $K=\pi_*
\overline\Q_{l,X}[d] \in Perv(Y)$ implies $E=\pi^*(K) = \pi^*\pi_* \overline\Q_{l,X}[d] \in
{}^p D^{\leq 0}(X)$ by [KW] lemma III.7.1. Therefore the distinguished triangle of perverse
truncation and the long exact Hom-sequence and the truncation axiom loc. cit. p.75) gives
$$ Hom_{D_c^b(X,\overline\Q_l)}\Bigl(\pi^*\pi_*\overline\Q_{l,X}[d],
\overline\Q_{l,X}[d]\Bigr) \hookrightarrow Hom_{D_c^b(X,\overline\Q_l)}\Bigl(B,
\overline\Q_{l,X}[d]\Bigr) \
$$ for $B= {}^p H^0(\pi^*\pi_*\overline\Q_{l,X}[d])$. Recall the process of truncation (using
notation from KW p. 140): It is clear, that for $E=\pi^*(K)=\pi^*\pi_*(\overline\Q_{l,X})$ the
exact triangle $({}^p\tau_{\leq -1}j_V^*(E), j_V^*(E), {}^pH^0(j_{V}^*(E)))$ defines an exact
triangle $(F,E,Rj_{V,*} {}^p\tau_{\geq 0}(j_V^*(E)))$ by adjunction. Notice ${}^p\tau_{\geq
0}(j_V^*(E))= \overline\Q_{l,V}[d]^{\vert G\vert}$. By the definition of truncation $B={}^p
H^0(E)$ sits in a diagram
$$ \xymatrix{ S_1 \ar[r] & B={}^p H^0(E) \ar[r] &
Rj_{V,*}(\overline \Q_{l,V}[d])^{\vert G \vert}  \cr & E \ar[u] & \cr  & A \ar[u] & \cr }
$$ with perverse sheaves $S_1$ and $A$ supported in $X\setminus V$. This implies
$$ B\vert_V = {}^p H^0(E)\vert_V \cong
\overline\Q_{l,V}[d]^{\vert G\vert} \ .$$

\bigskip\noindent
\underbar{Weight filtration of $B$}: Since $w(\overline\Q_l[d])\leq d$ also $w(K)\leq d$, using
that $\pi$ is proper. But then also $w(E)\leq d$ for $E=\pi^*(K)$. Now [KW] corollary 10.2
implies $w(B)\leq d$ for $B={}^p H^0(E)$. The weight filtration of the perverse sheaf $B$
defines an exact sequence of perverse sheaves
$$ 0 \to B_{w<d} \to B \to B_{w=d} \to 0 \ ,$$
where $B_{w=d}$ is pure of weight $d$ and $w(B_{w<d}) < d$. By assumption $\overline\Q_{l,X}$
has weight $\geq d$. Hence [KW] proposition II.12.6 and lemma III.4.3 imply
$$Hom_{D_c^b(X,\overline\Q_l)}(B_{w<d}, \overline\Q_{l,X}[d]))=0\ .$$ Hence the long exact
Hom-sequence gives an inclusion
$$ Hom_{D_c^b(X,\overline\Q_l)}\Bigl(B,
\overline\Q_{l,X}[d])\Bigr) \hookrightarrow Hom_{D_c^b(X,\overline\Q_l)}\Bigl(B_{w=0},
\overline\Q_{l,X}[d])\Bigr) \ .$$ To prove the lemma, it therefore is enough by the chain of
inclusion obtained to prove the estimate $$dim\Bigl(Hom_{D_c^b(X,\overline\Q_l)}(B_{w=0},
\overline\Q_{l,X}[d])\Bigr) \leq \vert G\vert \ .$$ By purity the perverse sheaf $B_{w=d}$ is a
direct sum of $j_{V,!*}(B\vert_V)$ and a perverse sheaf $C$ supported in $X\setminus V$. Since
$\overline\Q_{l,X}[d]$ is irreducible, there does not exist a nontrivial morphism from $C$ to
$\overline\Q_{l,X}$. We have also seen above $B\vert_V = \overline\Q_{l,V}^{\vert G\vert}$.
Since $\overline\Q_{l,X}[d]$ is irreducible, therefore $j_{V,!*}(B\vert_V) \cong
\overline\Q_{l,X}[d]^{\vert G\vert} $, since $B_{w<d}\vert_V =0$ and $B_{w=0}\vert_V =
\overline\Q_{l,V}^{\vert G\vert} $. In other words by Gabber's theorem
$$ dim\Bigl(Hom_{D_c^b(X,\overline\Q_l)}(B_{w=0},
\overline\Q_{l,X}[d])\Bigr) = \vert G\vert\ . $$ %By the inclusion (*) from above therefore
%$$ dim(Hom_{D_c^b(X,\overline\Q_l)}(B,
%\overline\Q_{l,X}[d])) \leq \vert G\vert. $$
This proves the lemma, i.e. forces all $F_\rho$ to irreducible.

%\bigskip\noindent
%\begin{Lemma}\label{normal} {\it If $Y$ is a normal variety of dimension $d$, and admits\footnote{Probably
%normality is enough.} a singularization of singularities $\pi:\tilde Y \to Y$, then $\lambda_Y=
%\overline\Q_{l,Y}[d]$.} \end{Lemma}
%
%\bigskip\noindent
%\underbar{Proof}: Let $U$ be open dense smooth in $Y$ and let
%$j:U\hookrightarrow Y$ be the inclusion map. Since
%$\lambda_Y={\cal H}^{-d}(\delta_Y)$ is $j_*\overline\Q_{l,Y}$, and
%since  there exists a natural inclusion
%$$ \overline\Q_{l,Y} \hookrightarrow
%j_*\overline\Q_{l,U} =\lambda_Y  \ ,$$ it therefore suffices, that
%the stalks of $j_*\overline\Q_{l,U}$ have all dimension $\leq 1$.
%Now, since $R\pi_*\overline \Q_{l,\tilde Y} = \delta_Y \oplus X $
%by the decomposition theorem. Hence
% $\lambda_Y={\cal H}^{-d}(\delta_Y)$ is a subsheaf of $R^0\pi_* \overline
%\Q_{l,\tilde Y}$. Since $Y$ is normal, Zariski's main theorem
%(Stein factorization) implies that the fibers of $\pi: \tilde Y\to
%Y$ are connected. Hence the stalks of $R^0\pi_* \overline
%\Q_{l,\tilde Y}$ have $\overline\Q_l$-dimension equal to one. This
%implies $\overline\Q_{l,Y}\cong j_*\overline\Q_{l,U}$, and proves
%the claim.

\chapter{Convolution for abelian varieties}

\bigskip\noindent
\section{Abelian varieties}\label{Begin}

\bigskip\noindent
Suppose $X$ is an abelian variety over $k$ with addition law
$$ X\times X \overset{a}{\to} X \ .$$
Let $g$ be its dimension. For complexes $K,L$ in $D_c^b(X,\overline\Q_l)$ the convolution $K*L$
in $ D_c^b(X,\overline\Q_l)$ is defined by $$ K*L = Ra_*(K\boxtimes L) \ .$$ Convolution
commutes with complex shifts $K[m]*L[n]=K*L[m+n]$.

\bigskip\noindent
\underbar{Commutativity}: The tensor product  of two sheaf complexes $F$ and $G$ (in our case
$F=pr_1^*(K)$ and $G=pr_2^*(L)$) is defined by $ (F\otimes G)^n \ = \ \bigoplus_{p+q=n} F^p
\otimes G^q $ with differentials $ d(F\otimes G)\big\vert (F^p \otimes G^q) = d(F)\vert F^p
\otimes Id + (-1)^p Id \otimes d(G)\vert G^q$. The commutativity constraints for complexes
$\psi_{F,G}=(\psi^n_{F,G})$
$$ \psi_{F,G}^n: (F\otimes G)^n \ = \ \bigoplus_{p+q=n} F^p \otimes G^q
\longrightarrow \ \bigoplus_{q+p=n} G^q \otimes F^p \ = \ (G\otimes F)^n \ $$
 are defined by $
\psi_{F,G}^n\big\vert F^p \otimes G^q = (-1)^{pq} \psi_{F^p,G^q}$,
where $\psi_{F^p,G^q}$ are the \lq{trivial\rq}\  commutativity
constraints for sheaves. In particular for a sheaf $F$ and a sheaf
complex $G$ one has $F\otimes G=G\otimes F$. This induces functorial
isomorphisms $\tilde\psi_{K,L}: \sigma_2^*(K\boxtimes L) \cong
L\boxtimes K$ via the  canonical isomorphisms via
$\tilde\psi_{K,L}=\sigma_{12}^*(\psi_{pr_1^*(K),pr_2^*(L)})$
$$ K\boxtimes L  = pr_1^*(K)\otimes^L pr_2^*(L)  \ \cong\
\sigma_{12}^*(pr_1^*(L)\otimes^L pr_2^*(K)) = \sigma_{12}^*(L\boxtimes K)  \ .$$ Notice $
pr_1^*(K)\otimes^L pr_2^*(L) = \sigma_{12}^*( pr_2^*(K)\otimes^L pr_1^*(L))$ for the
permutation $\sigma_{12}:X^2 \to X^2$, which are defined by $\sigma_{12}(x,y)=(y,x)$. The
isomorphism $\tilde\psi_{K,L}$ above is the pullback via $\sigma_{12}^*$ of the isomorphism $
pr_1^*(K)\otimes^L pr_2^*(L) \cong pr_2^*(L) \otimes^L pr_1^*(K)$ for the ordinary tensor
product $\otimes^L$ of the triangulated category $D_c^b(X^2,\overline\Q_l)$ described by the
ordinary tenor product of sheaf complexes $F=pr_1^*(K)$ and $G=pr_2^*(L)$. This implies the
existence of functorial commutativity constraints $\Psi_{K,L}: K*L\cong L*K$ for the
convolution product via
$$ \xymatrix@C+1cm{ K*L \ar@{=}[r]\ar[d]^\sim_{\Psi_{K,L}} & Ra_*(K\boxtimes L)\ar[r]^{Ra_*\tilde\psi_{K,L}^{-1}} & Ra_*(\sigma_{12}^*(L\boxtimes
K))    \ar@{=}[d] \cr L*K \ar@{=}[r] & Ra_*(L\boxtimes K) \ar@{=}[r] & Ra_*(L\boxtimes K) \cr}
$$ The commutativity constraints $\Psi_{K,L}=Ra_*\tilde\psi_{K,L}^{-1}: K*L \cong L*K$ are derived as
direct images under $a:X\times X\to X$ from the commutativity constraints of the outer tensor
products, and from  $Ra_*\circ \sigma_{12}^*=Ra_*$

%via the  canonical isomorphisms $\tilde\psi_{K,L}=\sigma_{12}^*(\psi_{pr_1^*(K),pr_2^*(L)})$
%$$ K\boxtimes L  = pr_1^*(K)\otimes^L pr_2^*(L)  \ \cong\
%\sigma_{12}^*(pr_1^*(L)\otimes^L pr_2^*(K)) = \sigma_{12}^*(L\boxtimes K)  \ .$$ Notice $
%pr_1^*(K)\otimes^L pr_2^*(L) = \sigma_{12}^*( pr_2^*(K)\otimes^L pr_1^*(L))$ for the
%permutation $\sigma_{12}:X^2 \to X^2$, which are defined by $\sigma_{12}(x,y)=(y,x)$. The
%isomorphism $\tilde\psi_{K,L}$ above is the pullback via $\sigma_{12}^*$ of the isomorphism $
%pr_1^*(K)\otimes^L pr_2^*(L) \cong pr_2^*(L) \otimes^L pr_1^*(K)$ for the ordinary tensor
%product $\otimes^L$ of the triangulated category $D_c^b(X^2,\overline\Q_l)$ described by the
%ordinary tenor product of sheaf complexes as follows: The tensor product of two sheaf complexes
%$F$ and $G$ (in our case $F=pr_1^*(K)$ and $G=pr_2^*(L)$) is defined by
%$$ (F\otimes G)^n \ = \ \bigoplus_{p+q=n} F^p \otimes G^q $$
%$$ d(F\otimes G)\big\vert (F^p \otimes G^q) =
%d(F)\vert F^p \otimes Id + (-1)^p Id \otimes d(G)\vert G^q \ .$$ The commutativity constraints
%for complexes $\psi_{F,G}=(\psi^n_{F,G})$
%$$ \psi_{F,G}^n: (F\otimes G)^n \ = \ \bigoplus_{p+q=n} F^p \otimes G^q
%\longrightarrow \ \bigoplus_{q+p=n} G^q \otimes F^p \ = \ (G\otimes F)^n \ $$ are defined by $$
%\psi_{F,G}^n\big\vert F^p \otimes G^q = (-1)^{pq} \psi_{F^p,G^q} \ ,$$ where $\psi_{F^p,G^q}$
%are the \lq{trivial\rq}\  commutativity constraints for sheaves. In particular for a sheaf $F$
%and a sheaf complex $G$ one has $F\otimes G=G\otimes F$.

\bigskip\noindent
\underbar{Associativity}: The existence of functorial associativity isomorphisms
$\varphi_{K,L,M}: (K\boxtimes L)\boxtimes M \cong K\boxtimes (L\boxtimes M)$ for the exterior
tensor product of complexes implies, that the convolution is associative with certain
functorial associativity constraints $\Phi_{K,L,M}:(K*L)*M \cong K*(L*M)$ via
$$ \xymatrix{ (K*L)*M \ar@{=}[r]\ar@{=}[dd]_{\Phi_{K,L,M}} & Ra_* R(a\times id)_* ((K\boxtimes L)\boxtimes M)\ar@{=}[d]^{Ra_*R(a\times id)_* id_{K,L,M}} \cr
           & Ra_* R(a\times id)_* (K\boxtimes (L\boxtimes M))\ar@{=}[d] \cr
K*(L*M) \ar@{=}[r] & Ra_* R(id\times a)_* (K\boxtimes (L\boxtimes M)) \cr} \ $$ since $a\circ
(a\times id)= a\circ (id\times a)$. Associativity $(F\otimes G)\otimes H = F\otimes (G\otimes
H)$ is satisfied for sheaf complexes, hence also in $D_c^b(X,\overline\Q_l)$ in the strict
sense
$$ d((F\otimes G)\otimes H)\vert (F^p\otimes G^q \otimes H^r) $$
$$ = d(F)\vert F^p \otimes id\otimes id + (-1)^p\cdot id\otimes d(G)\vert
G^q \otimes id + (-1)^{p+q}\cdot id\otimes id\otimes d(H)\vert H^r
$$
$$ = d(F\otimes (G\otimes H))\vert (F^p\otimes G^q \otimes H^r) \ .
$$
The commutativity and associativity constraints for the convolution product satisfy the usual
hexagon and pentagon axiom.

\bigskip\noindent
\goodbreak
\bigskip\noindent
\underbar{Supports}: By the proper base change theorem the stalk of
the cohomology sheaf $ {\cal H}^\nu(K*L)_z$ of $K*L$ at a point $z$
of $X$ is $H^\nu(X,K\otimes^L \varphi_z^*(L))$, where
$$\varphi_z:X\to X\quad , \quad \varphi_z(x)=z-x \ $$ and $\varphi_z^2=id_X$. For a complex $K\in
D^b_c(X,\overline\Q_l)$ define $supp(K)= \cup_\nu S^\nu$ be the
union of the supports $S^\nu$ of the cohomology sheaves ${\cal
H}^{-\nu}(K)$. By definition $S^\nu$ is  Zariski closed. Then $
K\otimes^L \varphi_z^*(L) $ is supported in $supp(K) \cap (z -
supp(L))$. Therefore by induction
$$ supp(K_1*\cdots * K_r) \subseteq supp(K_1) + \cdots + supp(K_r)
\ .$$ If $Y,Y'$ are irreducible normal closed subvarieties of the abelian variety $X$, we
obtain
$$ supp(\delta_{Y}*\delta_{Y'}) = supp(\delta_Y) +
supp(\delta_{Y'})= Y+Y' \ .$$ In fact ${\cal H}^{-\nu}(\delta_{Y}*\delta_{Y'})_z=0$ holds for
all $\nu > dim(Y) + dim(Y')$ and ${\cal H}^{-\nu}(\delta_{Y}*\delta_{Y'})_z$ is equal to
$H^0(Y\cap (z-Y'),\overline\Q_l)\neq 0$ for $\nu = dim(Y) + dim(Y')$ and $z\in Y+Y'$, which is
Zariski closed.

\bigskip\noindent
Concerning Verdier duality $D:D_c^b(X,\overline\Q_l)\to
D_c^b(X,\overline\Q_l)$ we have $D(K*L)=DRa_*(K\boxtimes L) =
Ra_*(D(K\boxtimes L))= Ra_*(DK\boxtimes DL) = DK*DL$ by Poincare
duality, since $a$ is proper. The convolution $K*L$ of selfdual
complexes $K\cong DK$ and $L\cong DL$ is selfdual. Similarly the
convolution $K*L$ of mixed complexes $K$ and $L$ of weight $\leq w$
respectively weight $\leq w'$ is mixed of weight $\leq w+w'$. If $K$
and $L$ are pure of weights $w$ respectively $w'$, then $K*L$ is
pure of weight $w+w'$.

\bigskip\noindent
\underbar{Relative K\"unneth theorem}: We now review \cite{SGA}, XVII.5.4 (in a slightly
modified context). Given morphisms $f:X\to Y$ and $g:X'\to Y'$ complexes $K\in
D_c^b(X,\overline\Q_l)$ and $L\in D_c^b(Y,\overline\Q_l)$ consider $p^*K\otimes q^*L \cong
K\boxtimes L\in D_c^b(X\times Y,\overline\Q_l)$ and the morphism $h=f\times g:X\times Y\to
X'\times Y'$. There exists a natural morphism
$$ \nu:\ Rf_*(K)\boxtimes Rg_*(L) \to Rh_*(K\boxtimes L)\ ,$$
which is induced by the basechange morphisms $(p')^*R^rf_* K \to R^r h_*(p^*(K))$ and
$(q')^*R^sf_* L \to R^s h_*(q^*(L))$ and cup-product $R^rh_*(p^*(K))\times R^sh_*(q^*(L)) \to
R^{r+s}h_*(K\boxtimes L)$ (see \cite{Mi2}, p.172). If both morphisms $f$ and $g$ are proper,
then $\nu$ is an isomorphism. This follows from the proper basechange theorem, which allows to
reduce this to the case $X'=Y'=Spec(k)$ (see \cite{SGA}, XVII Theorem 5.4.3, or \cite{Mi2}, thm
8.5 and p.262 loc. cit.).

\bigskip\noindent
\goodbreak
\bigskip\noindent
\underbar{Functoriality}: Suppose $f:X\to Y$ is a homomorphism of abelian varieties. For
complexes $K,L\in D_c^b(X,\overline\Q_l)$ the direct image $$Rf_*: D_c^b(X,\overline\Q_l)\to
D_c^b(Y,\overline\Q_l)$$ commutes with convolution
$$ Rf_*(K*L) \ = Rf_*(K)\ *\ Rf_*(L) \ .$$
This is evident from the commutativity $a\circ (f\times f)=f\circ a$ of the diagram
$$ \xymatrix@+0.5cm{ X\times X \ar[r]^a\ar[d]_{f\times f}\ & X \ar[d]^f \cr
Y\times Y \ar[r]^a & Y \cr } \  $$ and the relative K\"unneth formula for $f\times f$: Indeed
$Rf_*(K*L)= Rf_*Ra_*(K\boxtimes L) = Ra_* R(f\times f)_*(K\boxtimes L)= Ra_*(Rf_*(K)\boxtimes
Rf_*(L)) = Rf_*(K)*Rf_*(L)$.

\bigskip\noindent
\underbar{The dual}: For $K\in D_c^b(X,\overline \Q_l)$ we define the adjoint dual $K^\vee \in
D_c^b(X,\overline\Q_l)$ by
$$K^\vee = (-id_X)^*D(K)\ .$$ $(-id_X)^* D= D(-id_X)^*$,
$(-id_X)^*(-id_X)^*K=K$ and existence of a natural isomorphism $j_K: K\cong D(D(K))$ implies,
that there exists a natural isomorphism $\lambda_K:K\cong (K^\vee){}^\vee $. This defines a
natural transformation $\lambda:id \to {}^\vee{}^\vee$, which allows to identify  the identity
functor with the bidual functor $\lambda_K:K\cong K^\vee{}^\vee$.

\begin{Lemma}\label{eas}
$(\lambda_K)^\vee = (\lambda_{K^\vee})^{-1}$. \end{Lemma}

\bigskip\noindent
\underbar{Proof}: This can be reduced to the corresponding statement $D(j_K) = (j_{D(K)})^{-1}$
for Verdier duality, which follows from \cite{Sav} (3.2.3.9) and \cite{Sav} p. 51- 57, since
for the $\otimes$-category $(D_c^b(X,\overline\Q_l),\otimes^L)$ the functor $Z\mapsto
Hom_{D_c^b(X,\overline\Q_l)}(Z\otimes^LX,Y)$ is representable by
$Hom_{D_c^b(X,\overline\Q_l)}(Z,R{\cal H}om(X,Y))$ for $ R{\cal H}om(X,Y) = D(X\otimes^L
D(Y))$. We notice in passing, that $K\mapsto D(K)[-2g]$ is related to the functor $f\mapsto
f^t$ in \cite{Sav} 3.2.3.1 using the unit object $\overline\Q_{l,X}$ of
$(D_c^b(X,\overline\Q_l),\otimes^L)$.

\bigskip\noindent
$K\mapsto K^\vee$ is a contravariant functor on $D_c^b(X,\overline\Q_l)$
$$ f: K\to L  \quad , \quad f^\vee: L^\vee \to K^\vee \ .$$
 There exist functorial isomorphisms
$$ i_K: H^\bullet(X,K^\vee) \cong H^\bullet(X,K)^\vee \ $$
defined by the composite morphism
$$ \xymatrix{ H^\bullet(X,K^\vee)\ar@{=}[r] & H^\bullet(X,(-id_X)^*(DK)) & H^\bullet(X,DK)\ar[l]_-{(-id_X)^*}^-{\sim}
\ar[r]^{PD}_\sim & H^\bullet(X,K)^\vee\cr} \ $$ with the last isomorphism given by Poincare
duality.

\begin{Lemma}\label{compdia} For complexes $A,B\in D_c^b(X,\overline\Q_l)$ and complex maps $\varphi: A\to B$
and dual $\varphi^\vee: B^\vee \to A^\vee$ the following diagram commutes ($i$ is a functorial
isomorphism)
$$ \xymatrix{ H^\bullet(X,B^\vee)\ \ar[d]_{H^\bullet(\varphi^\vee)} \ar[r]^-{i_B}_-\sim & \ H^\bullet(X,B)^\vee \ar[d]^{H^\bullet(\varphi)^\vee}\cr
H^\bullet(X,A^\vee)\ \ar[r]^-{i_A}_-\sim & \ H^\bullet(X,A)^\vee \cr}$$ and the natural
identification $\lambda_A: A\to A^\vee{}^\vee$ induces a commutative diagram
$$ \xymatrix{ H^\bullet(X,A)\ \ar@{=}[d]_{can}\ar[r]^-{H^\bullet(\lambda_A)}_-\sim & \ H^\bullet(X,A^\vee{}^\vee)\ar[d]^{i_{A^\vee}}_-\sim \cr
H^\bullet(X,A)^\vee{}^\vee\ \ar[r]^-{(i_A)^\vee}_-\sim & \ H^\bullet(X,A^\vee)^\vee. \cr} $$
\end{Lemma}

\bigskip\noindent
\underbar{Proof}: The first statement is easily reduced to the functoriality of Poincare
duality. Also for the second statement one reduces to the commutativity of the diagram
$$ \xymatrix{ H^\bullet(X,A)\ \ar@{=}[d]_{can}\ar[r]^-{H^\bullet(j_A)}_-\sim & \ H^\bullet(X,D(D(A)))\ar[d]^{PD_{A^\vee}}_-\sim \cr
H^\bullet(X,A)^\vee{}^\vee\ \ar[r]^-{(PD_A)^\vee}_-\sim & \ H^\bullet(X,D(A))^\vee. \cr} $$ By
\cite{SGA45}, 1.2, p.155 and \cite{SGA}, lemma 3.2.3, p.583 Poincare duality $$PD_A:
H^\bullet(X,D(A)) \cong H^\bullet(X,A)^\vee$$ is induced from the pairing $
H^\bullet(X,DA)\times H^\bullet(X,A) \to H^\bullet(X,D(A)\otimes^L A)$ followed by
$H^\bullet(X,ev_A)$ using the evaluation $ev_A: D(A)\otimes^L A\to K_X=\Q_{l,X}[2g]$. Here we
in addition used, that $X$ is smooth and therefore $K_X=\Q_{l,X}[2g]$. Now compare with $
H^\bullet(X,D(D(A)))\times H^\bullet(X,D(A)) \to H^\bullet(X,D(A)\otimes^L A)$ (up to switching
sides) as required in the last diagram, by the commutative diagram of \cite{Sav}, 3.2.3.8
$$ \xymatrix{ A\otimes^L D(A)[-2g]\ \ar[d]^{\psi_{A,D(A)}}\ar[r]^-{j_A\otimes^L id} & \ D(D(A))\otimes^L
D(A)[-2g] \ar[d]_{ev_{D(A)}} \cr D(A)[-2g]\otimes^L A\ \ar[r]^{ev_A} & \ \Q_l.\cr} $$ The
condition of loc. cit are satisfied for the category $(D_c^b(X,\overline\Q_l),\otimes^L)$. For
this see proof of lemma \ref{eas}. This completes the proof.

\bigskip\noindent
For a homomorphism $f:X\to Y$ between abelian varieties the relative Poincare duality theorem
implies
$$ Rf_*(K)^\vee = Rf_*(K^\vee) \ .$$
Since homomorphisms are proper, $Rf_*$ preserves purity of complexes and maps direct sums of
translates of pure perverse sheaves to direct sums of translates of pure perverse sheaves by
Gabber's decomposition theorem. Furthermore there exist functorial isomorphisms
$$ (K*L)^\vee \cong K^\vee * L^\vee \ $$
since $DK*DL\cong D(K*L)$ and $(-id_X)^*(K*L)\cong (-id_X)^*(K)
* (-id_X)^*(L)$, which follows by basechange from the cartesian
diagram
$$ \xymatrix@+0.5cm{ X\times X \ar[r]^a\ar[d]_{-id_X \times -id_X} & X\ar[d]^{-id_X}  \cr
X\times X \ar[r]^a & X  \cr } \ .$$
 Since ${}^\vee$ is a
(contravariant) tensor functor with respect to the convolution, this usually allows to identify
$K$ and ${K^\vee}{}^\vee$ and $f^{\vee}{}^\vee$ and $f$, such that $(K^\vee*K)^\vee =
K^\vee{}^\vee *K^\vee = K*K^\vee$.

 %For later use we notice, that obviously the following holds. Given
%$$ f=f_1+f_2: K_1\oplus K_2 \to K, \quad \mbox{then} \quad f^\vee=(f_1^\vee,f_2^\vee): K^\vee
%\to K_1^\vee \oplus K_1^\vee \ .$$

\goodbreak
\section{The Grothendieck ring $K_*(X)$}
Let $K^0(Perv(X))$ denote the Grothendieck group  of the abelian category of perverse sheaves
$Perv(X)$. Let $t^{1/2}$ denote an indeterminate and consider the \lq{enlarged\rq}\
Grothendieck group
$$ K^0_*(X)=K^0(Perv(X))\otimes_{\Z} \Z[t^{1/2},t^{-1/2}] \ .$$
As a variant we may also consider $ K^0(Perv(F,X))\otimes_{\Z}
\Z[F]$ for the larger group $Z[F]$ containing $\Z[t^{1/2},t^{-1/2}]$
(see [KF], p. 189) or the Grothendieck of the abelian category
$Perv_m(X)$ of mixed perverse Weil sheaves. Notice that Verdier
duality \lq{acts\rq}\ on the Grothendieck group and this action is
considered to be extended to the ring $\Z[t^{1/2},t^{-1/2}]$ by
$Dt^{\pm 1/2}=t^{\mp 1/2}$. Pars pro toto we restrict ourselves to
the subcategory of mixed perverse Weil sheaves, in the case where
$k$ is the algebraic closure of a finite field. The other cases are
similar. Under the assumption made each irreducible object in
$Perv(X)$ is a pure complex of some weight $w\in\Z$, which defines a
weight graduation in the enlarged Grothendieck group. It is easy to
see, that the convolution product induces a ring structure on the
enlarged Grothendieck group. For simple mixed perverse sheaves $K$
and $L$ we have $K*L = \oplus_{\nu\in\Z} A_\nu[\nu](\nu/2)$ for
irreducible perverse sheaves $A_\nu$ using Gabber's theorem (for
simplicity we choose a square root of $q$ to define the half
integral Tate twists). We define the product of the classes of $K$
and $L$ in the extended Grothendieck group to be $\sum_\nu cl(A_\nu)
t^{-\nu/2}$. This makes the enlarged Grothendieck group into a
graded ring with respect to the weight filtration. The unit element
of this ring is the class of the perverse skyscraper sheaf
$\delta_{0}$, which is concentrated in the neutral element $0\in X$.
Assigning to a perverse sheaf $K$ its Betti polynomial $h(K,t)$
$$ h(K,t) = \sum_{\nu\in \Z} h^{\nu}(X,K) \cdot t^{\nu/2}\quad , \quad h^\nu(X,K)=dim \ H^\nu(X,K) $$
defines a ring homomorphism
$$ h: K^0_*(X) \to
\Z[t^{1/2},t^{-1/2}] \ $$ by the K\"{u}nneth formulas. Since $cl(\delta_X)*cl(K)=
cl(\delta_X)\cdot h(K,t)$ we now define
$$ K_*(X) = K^0_*(X) \Big/
\Z[t^{1/2},t^{-1/2}] \cdot cl(\delta_X) \ .$$ $K_*(X)$ with its weight graduation is a graded
ring, an enhanced version of the usual homology ring $(H_*(X,\Z),*)$ of $X$ endowed with the
$*$-product (dual to the cohomology ring endowed with the cup-product). In the following we
often do not distinguish between perverse sheaves $K$ on $X$ and their classes $cl(K)$ in
$K_*(X)$, writing $K$ in both cases for simplicity  whenever there is little danger of
confusion. $t$ is an element of weight $0$. $h$ induces a ring homomorphism
$$  K_*(X) \overset{h_*}{\longrightarrow} \Z[t^{1/2},t^{-1/2}]/h(\delta_X,t) \cong
\Z[t^{1/2},t^{-1/2}]/(t^{1/2}+2+t^{-1/2})^{g} \ ,$$ which under the specialization $t^{1/2}
\mapsto -1$ becomes  $K\mapsto \sum (-1)^\nu h^\nu(X,K)$.

\bigskip\noindent
The degree $deg_{t^{1/2}}(h(K,t))$ of $h(K,t)$ is called the degree of $K$, which is $\leq g-1$
for a perverse sheaf $K$ on $X$ without a nontrivial constant perverse quotient sheaf. Similar
degree $deg_{t^{-1/2}}( h(K,t))\geq g-1$ for a perverse sheaf $K$ on $X$ without a nontrivial
constant perverse subsheaf on $X$. If for $K\in Perv(X)$ the highest exponent of $h(K,t)$ in
$t^{1/2}$ and $t^{-1/2}$ is $\leq i$, we say $K$ is of type $i$. $K$ is of $0$-type, if
$h(K,t)$ is a constant polynomial. If $K$ is pure and of type $i$, subquotients of $K$ are
again of type $i$. Also the Verdier dual $D(K)$ and $K^\vee$. If $K$ and $L$ are of type $i$,
their direct sum is of type $i$ and their convolution product of types $i$ and $j$ is of type
$i+j$.

\bigskip\noindent
\begin{Lemma} \label{null}{\it For a perverse sheaf $K$ without constant quotient or constant perverse subsheaf,
such that $K*K=d\cdot K$ holds in $K_*(X)$ for some integer $d\neq
0$, the polynomial $h(K,t)$ is constant.}
\end{Lemma}

\bigskip\noindent
\underbar{Proof}: $R=\Z[X,X^{-1}]/(X+2+X^{-1})^g \cong \Z[X,X^{-1}]/(X+1)^{2g} \cong
\Z[X]/(X+1)^{2g} \cong \Z[\varepsilon]/\varepsilon^{2g}$. Therefore $K*K=d\cdot K$  for some
integer $d\neq 0$ implies $P(P-d)=0$ in $R$ for the residue class $P$ of $h(K,t)$ in $R$.
Either $P\in (\varepsilon)$ and $P=0$ using $P^2=dP$, or $P\not\in (\varepsilon)$. But then
$P\in R^*$ and $P(P-d)=0$ implies $P=d$.

\bigskip\noindent
\underbar{Functoriality}: A homomorphism $f:X\to Y$ between abelian
varieties induces a ring homomorphism $f_*: K_*^0(X)\to K_*^0(Y)$ by
$K\mapsto Rf_*(K)$.

\section{Translation invariance}\label{TI}
For the skyscraper sheaf $\delta_0$ concentrated at the neutral element $0\in X$ we have $
K*\delta_0 = K$. More generally $K*\delta_{\{x\}} = T^*_{-x}(K)$, where $x$ is a closed
$k$-valued point in $X$ and where $T_{x}(y)=y+x$ denotes the translation $T_x:X\to X$ by $x$.
In fact
$$ T^*_y(K*L) \cong T^*_y(K)*L \cong K*T^*_y(L)
$$ holds for all $y\in X(k) $ by the proper basechange theorem
using the diagram
$$ \xymatrix{ X\times X \ar[d]_{T_y\times id} \ar[r]^-a & X\ar[d]^{T_y} \cr X\times X \ar[r]^-a & X \cr } \ .$$
 For $K\in D_c^b(X,\overline\Q_l)$ let $Aut(K)$ be the abstract group of all closed $k$-valued points $x$ of
$X$, for which $T^*_x(K)\cong K$ holds.
 A complex $K$ is called translation-invariant, provided $Aut(K)=X(k)$.
 As a consequence of the formulas
above, the convolution of an arbitrary $K\in D_c^b(X,\overline\Q_l)$ with a
translation-invariant complex on $X$ is a translation-invariant complex. In fact $Aut(K)$ is
contained in $Aut(K*L)$.

\bigskip\noindent
Suppose $K$ is a translation-invariant perverse sheaf on $X$. Then each irreducible constituent
of $K$ is translation-invariant (e.g. use Lemma \ref{l9}). Since $K=\delta_E$ for some triple
$(E,U,Y)$, where $E$ is a smooth etale sheaf on $U$ with Zariski closure $Y=\overline U$,
translation invariance of $K$ implies $Y=X$, since the subvariety $Y$ is uniquely determined by
$K$. Hence $E$ must be a translation-invariant $\overline\Q_l$-sheaf over a Zariski open dense
subset $U$ of $X$. $E$ corresponds to some $\overline\Q_l$-adic representation
$$ \rho_E: \pi_1(X\setminus V_E, x_0) \to Gl(E_{x_0}) $$
for the etale fundamental group of some open dense subvariety
$U=X\setminus V_E$ of the abelian variety $X$ and some geometric
point $x_0$ in $ U$. Suppose $V_E$ is chosen minimal. Then $V_E$ is
the ramification divisor of the $l$-adic coefficient system $E$.
$V_E$ must either have dimension $g-1$, or otherwise it must be
empty by the purity of branched points. If $V_E$ is empty, then the
coefficient system $E$ is smooth on $X$, and $E$ corresponds to a
representation of the etale fundamental group $\pi_1(X,x_0)$ of $X$.
The group $\pi_1(X,0)$ is abelian for an abelian variety $X$ (this
also holds for $char(k)\neq 0$ by \cite{Mu}, p.167). The irreducible
$\overline\Q_l$-adic representations $E_\chi$ of
$\pi_1(X,0)=\pi_1(X)^{ab}$ correspond to continuous characters (for
this notion see \cite{KW}, p.9)
$$ \chi \in Hom_{cont}(\pi_1(X,0),\overline\Q_l^*) \ , $$
which correspond to projective limits of characters of finite quotient groups
$$\pi_1(X)^{ab}/n\cdot \pi_1(X)^{ab} \cong Pic^0(X)[n]^*$$ (see \cite{Mi}, p.132). Since the
automorphism $x\mapsto -x$ acts on the first etale cohomology $H^1(X,\Z_l)$ via multiplication
by $-1$, we obtain $ (-id_X)^*E_\chi[g]\cong E_{\chi^{-1}}[g] \cong D(E_{\chi}[g])$. If $K$ is
a translation-invariant perverse sheaf, the ramification divisor $V_E$ of the underlying smooth
sheaf $E$ must be translation-invariant. This gives a contradiction unless $V_E=\emptyset$.
Therefore an irreducible translation-invariant perverse sheaf on $X$ is of the form
$\delta_{E}=\delta_\chi$ for the rank one smooth $\overline\Q_l$-sheaf $E=E_\chi$ attached to
some continuous character $\chi: \pi_1(X,0)\to \overline\Q_l^*$. Conversely for any such
character the perverse sheaf $\delta_\chi$ is translation-invariant. It is enough to show this
for characters $\pi_1(X,0) \to \mu_n$. The group $X$ acts trivially on $Pic^0(X)$ via ${\cal
L}\mapsto T_x^*({\cal L})\otimes{\cal L}^{-1}$. Hence it acts trivially on $H^1(X,\mu_n)\cong
Pic(X)[n]=\{a\in Pic(X)\ \vert \ na\}=0$, since $Pic(X)[n]$ is contained in $Pic(X)^0$ (the
Neron-Severi group is torsionfree).

\bigskip\noindent
For the trivial character $\chi$ the sheaf complex $E_\chi[g]=\delta_X$ is
translation-invariant in the stronger sense, that $a^*(\delta_X) \cong pr_2^*(\delta_X)$ holds
for the morphisms
$$ \xymatrix{ X\times X\ \ar@<1ex>[r]^-{a} \ar@<-1ex>[r]_-{pr_2} & \ X \cr} $$
defined by the addition $a$ and the second projection $pr_2$. Conversely, suppose $K$ is an
irreducible perverse sheaf on $X$ such that $a^*(K) \cong pr_2^*(K)$ holds. Then $K$ can be
rigidified along the zero section ([KW], p.187), thus becomes $X$-equivariant ([KW],p.188).
Hence $K\cong \delta_X$ by [KW] lemma 15.5. Notice, this implies
$$ \delta_X*K \cong \delta_X \otimes H^*(X,K) $$
for all complexes $K$. In fact $a\circ \rho=pr_1$ holds for the
automorphism
$$\rho: X\times X \to X\times X $$ defined by $\rho(x,y)=(x-y,y)$. Now
$\rho^*(\delta_X\boxtimes K) = \rho^*(pr_1^*(\delta_X)\otimes pr_2^*(K)) =
(pr_1\rho)^*(\delta_X)\otimes (pr_2\rho)^*(K) = (pr_1^*(-id_X))^*(\delta_X)\otimes pr_2^*(K)
\cong \delta_X\boxtimes K$ since $(-id_X)^*\delta_X\cong\delta_X$. We have a commutative
diagram
$$ \xymatrix{  X\times X \ar[rr]^\rho\ar[rd]_{pr_1} &  & X\times X \ar[ld]^a \cr
& X  & \cr} \ .$$ Hence $\delta_X*K=Ra_*(\delta_X\boxtimes K)=
Ra_*R\rho_*(\rho^*(\delta_X\boxtimes K))\cong Rpr_{1*}(\delta_X\boxtimes K)=\delta_X \otimes
R\Gamma(X,K)$, which proves the assertion $\delta_X*K\cong H^*(X,K)\otimes\delta_X$.

\bigskip\noindent
Now we return to the case of translation-invariant irreducible perverse sheaves $E_\chi$.
Notice $\rho^*(E_\chi\boxtimes E_{\chi'})= E_{\chi}\boxtimes E_{\chi'/\chi}$; it is enough to
show this for etale torsion sheaves, where it suffices to compute this in $H^1(X\times
X,\mu_n)\cong Pic^0(X\times X)[n]\cong Pic^0(X)[n] \times Pic^0(X)[n]$. Using this one shows as
above $E_\chi*E_{\chi'}=Ra_*(E_\chi\boxtimes E_{\chi'})= Ra_*R\rho_*(\rho^*(E_{\chi}\boxtimes
E_{\chi'})\cong Rpr_{1*}(E_{\chi}\boxtimes E_{\chi'/\chi})= E_{\chi}\otimes
R\Gamma(X,E_{\chi'/\chi})$. Now $R\Gamma(X,E_\chi)=0$ for  a nontrivial $\chi$. To see this
again it  suffices  to consider etale torsion $\Lambda_X$-module sheaves $E$ for an Artin
$K$-Algebra $\Lambda$. Then  we can assume $l\cdot E=0$ for some integer $l$. Hence $E$ is a
direct summand of the direct image $Rl_*\Lambda_{X}$ of the multiplication map $l:X\to X$. Then
$H^\bullet(X,E)=0$ for all nonconstant summands, since $H^\bullet(X,\Lambda_{X}) =
H^\bullet(X,Rl_*(\Lambda_{X}))$. This proves the claim, since the direct image
$Rl_*(\Lambda_{X})$ also contains the constant sheaf $\Lambda_{X}$ as a direct summand. This
proves the vanishing property claimed above, and implies $E_\chi*E_{\chi'}=0$ for all
$\chi\neq\chi'$. Therefore
$$ \delta_\chi * \delta_{\chi}\ =\ \delta_\chi \otimes H^*(X,\delta_X)  \
$$
$$ \delta_\chi * \delta_{\chi'} = 0 \quad , \quad \chi\neq \chi' \
.$$ Finally, since for pure perverse sheaves $K$ the convolution $\delta_\chi*K$ is
translation-invariant, Gabber's decomposition theorem implies the formula $\delta_\chi*K =
\sum_{\chi'} m(\chi')\cdot \delta_\chi'$ (in the extended Grothendieck group $K_*^0(X)$) with
certain coefficients $m(\chi')\in \Z[t^{1/2},t^{-1/2}]$. If we multiply by $\delta_\chi$, we
obtain $h(X,\delta_X)\cdot \delta_\chi*K = h(X,\delta_X)m(\chi)\cdot \delta_{\chi}$. Hence
$\delta_\chi*K = m(\chi)\cdot \delta_\chi$ (the additive group of $K_*^0(X)$ is torsion free,
since the subcategory of $Perv(X)$ generated by pure perverse Weil sheaves is semisimple). By a
comparison of the stalk cohomology at the point 0 we then obtain
$$ \delta_\chi*K = m(\chi)\cdot \delta_\chi $$
$$m(\chi)=h(t,\chi^{-1},K)=\sum h^\nu(X,K\otimes E_{\chi^{-1}})\cdot t^{\nu/2}\ .$$
The subcategory $Perv(X)_{inv}$ of the abelian category $Perv(X)$ consists of all perverse
sheaves, whose irreducible constituents are translation-invariant, is a Serre subcategory. Let
denote
$$ \overline{Perv}(X) $$
the corresponding  quotient abelian category of $Perv(X)$. The subcategory $T(X)$ of all $K\in
D_c^b(X,\overline \Q_l)$, for which ${}^pH^\nu(K)\in Perv(X)_{inv}$, is a thick subcategory of
the triangulated category $D_c^b(X,\overline\Q_l)$. Let
$$ \overline{D}_c^b(X,\overline\Q_l)$$
denote the corresponding triangulated quotient category. Then the convolution
$$ *:\ \overline D_c^b(X,\overline\Q_l) \times \overline D_c^b(X,\overline\Q_l) \to \overline
D_c^b(X,\overline\Q_l) $$ remains well defined. Furthermore, for a surjective homomorphism
$f:X\to Y$ between abelian varieties the direct image $Rf_*(K)$ of a complex $K\in T(X)$ is in
$T(Y)$. Hence $Rf_*$ induces a functor $\overline{D}_c^b(X,\overline\Q_l) \to
\overline{D}_c^b(Y,\overline\Q_l)$.

\section{Stalks of convolution products}\label{stalks}
For  $K,L\in D_c^b(X,\overline\Q_l)$ the stalk cohomology groups at a point $z\in X(k)$ of the
convolution product $K*L$ are
$$ {\cal H}^{-\nu}(K*L)_z \ = \ H^{-\nu}\Bigl(a^{-1}(z),K\boxtimes L\vert_{a^{-1}}(z)\Bigr)\ .$$
$a^{-1}(z)$ is isomorphic to $X$ and consists of the points $(x,z-x), x\in X(k)$. Hence
$$ {\cal H}^{-\nu}(K*L)_z \ = \ H^{-\nu}(X,K\otimes^L \varphi_z^*(L))\ .$$
For fixed $z\in X(k)$ we abbreviate $M=\varphi_z^*(L)$. Then there exists a spectral sequence
converging to $H^{-\nu}(X,K\otimes^L M)$
$$ H^{i}(X,{\cal H}^j(K\otimes^L M)) \Longrightarrow
H^{i+j}(X,K\otimes^L M) \ .$$  For $N=K\otimes^L M$ notice ${\cal H}^j(N) = \oplus_{a+b=j}
{\cal H}^a(K)\otimes {\cal H}^b(M)$. The $d_2$-terms look like
$$ \xymatrix{  H^{i-1}(X,{\cal H}^{j+1}(N)) \ar[drr] &
H^{i}(X,{\cal H}^{j+1}(N)) & H^{i+1}(X,{\cal H}^{j+1}(N))\cr
H^{i-1}(X,{\cal H}^{j}(N))\ar[drr]  & H^{i}(X,{\cal H}^j(N)) &
H^{i+1}(X,{\cal H}^{j}(N))\cr H^{i-1}(X,{\cal H}^{j-1}(N)) &
H^{i}(X,{\cal H}^{j-1}(N)) & H^{i+1}(X,{\cal H}^{j-1}(N))\cr } $$

\section{A vanishing theorem}\label{VT}
For  $K$ and $L$ in $D_c^b(X,\overline \Q_l)$ notice
$$
R{\cal H}om(\varphi_z^*(L),DK)\cong D(\varphi_z^*(L)\otimes^L DDK) \cong
D(\varphi_z^*(L)\otimes^L K)\ .$$ By Poincare duality this induces an isomorphism
$$ H^\nu(X,R{\cal
H}om(\varphi_z^*(L),DK)) \ \cong \ H^{-\nu}(X,\varphi_z^*(L)\otimes^L K)^\vee   \ .$$ Here
$V^\vee=Hom_{\overline\Q_l}(V,\overline\Q_l)$ denotes the dual vector space. On the other hand,
as already shown,
$$ H^{-\nu}(X,\varphi_z^*(L)\otimes^L K)^\vee \ \cong\ {\cal
H}^{-\nu}(L*K)_z^\vee \ \cong\ {\cal H}^{-\nu}((L*K)_z)^\vee \ .$$ Hence
$$ H^\nu(X,R{\cal
H}om(\varphi_z^*(L),DK)) \ \cong \ {\cal H}^{-\nu}((L*K)_z)^\vee   \ .$$ Now assume $K,L\in
{}^p D^{\leq 0}(X)$. Then $DK\in {}^p D^{\geq 0}(X)$ and $\varphi_z^*(L)\in {}^p D^{\leq
0}(X)$. Hence the cohomology groups $ H^\nu(X,R{\cal H}om(\varphi_z^*(L),DK))$ vanish for all
$\nu<0$ (see e.g. [KW] lemma III.4.3). For $\nu=0$ we obtain
$$ Hom_{Perv(X)}({}^p H^0(\varphi_z^*(L)),{}^p H^0(DK)) =  Hom_{D_c^b(X,\overline\Q_l)}(\varphi_z^*(L)),DK) $$ $$ = H^0(X,R{\cal
H}om(\varphi_z^*(L),DK)) = H^0(X,R{\cal H}om(\varphi_z^*(L),DK))  \ .$$  This implies

\bigskip\noindent
\begin{Lemma}\label{van} {\it For semi-perverse $K,L\in {}^pD^{\leq 0}(X)$ the following holds
\begin{enumerate} \item For $\nu>0$ the stalk cohomology ${\cal H}^\nu(DK*L)_z = 0$ vanishes
for all $z\in X$.
\item For $\nu=0$ the dual of the stalk cohomology at $z\in X(k)$ is
$${\cal H}^0(L*DK)_z^\vee  \ \cong \ Hom_{Perv(X)}({}^p H^0(\varphi_z^*(L)),{}^p H^0(K)) \ . $$
This isomorphism is functorial in $K$ and $L$.
\end{enumerate}}
\end{Lemma}

\bigskip\noindent
\begin{Corollary} \label{cor1}{\it If $L$ and $K$ are irreducible perverse sheaves on $X$, then
$$\{ z\ \vert \ {\cal
H}^0(L*K)_z\neq 0 \}  = \{z\in X\ \vert \ \varphi_z^*(L)\cong DK \} \ ,$$ hence
$$supp\ \Bigl({\cal
H}^0(L*K) \Bigr) =\overline{\{z\in X\ \vert \ \varphi_z^*(L)\cong DK \}}\ .$$ }
\end{Corollary}

\bigskip\noindent
$\{z\in X\ \vert \ \varphi_z^*(L)\cong DK \}$ is a torsor under $Aut(K)$, and is contained in
the closed subset
$$ S(K,L) =\{ z\in X\ \vert \ z - supp(L) = supp(K)\} \ .$$
Under the assumption on $K$ and $L$ of the last corollary
$$ \fbox{$ supp\ \Bigl({\cal H}^0(K*L)\Bigr) \subseteq S(K,L) $} \ .$$

\bigskip\noindent
Notice $DK\cong \varphi_z^*(L)$ is equivalent to $K^\vee \cong T_z^*(L)$. For irreducible
perverse sheaves $K$ and $L$ therefore
\begin{itemize} \item ${\cal H}^0(K*L)_z = 0$  if $K^\vee \not\cong T_x^*(L)$. \item ${\cal H}^0(K*L)_z \cong \overline\Q_l$ if $K^\vee \cong T_x^*(L)$.
\end{itemize} %By Gabber's theorem this  carries over to the case where $K$ and $L$ are pure
%perverse sheaves: For $K\cong\oplus n_A A$ and $L\cong \oplus m_A A$ with irreducible perverse
%sheaves $A$ $$ dim({\cal H}^0(K*L)_0) = \sum_A n_Am_A\ .$$ Hence the $\overline\Q_l$-dimension
%of ${\cal H}^0(L*DK)_0$ counts (with these multiplicities) the number of common perverse
%irreducible constituents of $\varphi_z^*(L)$ and $DK$.

\begin{Corollary}\label{zuvor}
For $K,L\in Perv(X)$ we have functorial isomorphisms $$ Hom_{Perv(X)}(K,L) \ \cong\ {\cal
H}^0(K*L^\vee)^\vee_0 \ .$$
\end{Corollary}

\bigskip\noindent
Now suppose the convolution $K*L^\vee$ is still a perverse sheaf. Since the cohomology sheaf
${\cal H}^0(M)$ of a perverse sheaf $M$ is supported in a finite set of closed points, ${\cal
H}^0(K*L^\vee)^\vee$ must be a skyscraper sheaf. Hence $Hom_{Perv(X)}(K,L) \ \cong\
\Gamma_{\{0\}}(X,{\cal H}^0(K*L^\vee))$ (sections with support in $0$). For a
translation-invariant sheaf $G$ on the other hand $\Gamma_{\{0\}}(X,{\cal H}^0(K*L^\vee))=0$.

\bigskip\noindent
This being said, let $Perv(X)'$ be a subcategory of $Perv(X)$ closed under $K\mapsto K^\vee$,
such that for all $K,L\in Perv(X)'$ the complex $K*L^\vee\in D_c^b(X,\overline\Q_l)$ is a
direct sum of a perverse sheaf and a sum $\oplus_\nu T_\nu[\nu]$ of translation-invariant
perverse sheaves $T_\nu\in Perv(X)$. Then the image of $Perv(X)'$ in the quotient category
$\overline{Perv}(X)$ defines a category $\overline{Perv}(X)'$ for which  convolution
$$ *:\ \overline{Perv}(X)' \times \overline{Perv}(X)' \to \overline{Perv}(X) $$
is well defined.

\begin{Corollary}\label{zuvor2}
Suppose $Perv(X)'$ is a semisimple subcategory of $Perv(X)$ such that for all $K,L\in Perv(X)'$
the complex $K*L^\vee\in D_c^b(X,\overline\Q_l)$ is a direct sum of a perverse sheaf and a sum
$\oplus_\nu T_\nu[\nu]$ of translation-invariant perverse sheaves $T_\nu\in Perv(X)$. Then we
have functorial isomorphisms
$$ Hom_{\overline{Perv}(X)'}(K,L) \ \cong\ \Gamma_{\{0\}}(X,{\cal H}^0(K*L^\vee)^\vee) \ .$$
\end{Corollary}

\begin{Corollary}\label{Canonical}
For $K,L,T\in Perv(X)$ suppose $T*K\in Perv(X)$ and $K^\vee*L\in Perv(X)$. Then 1) there exists
a canonical isomorphism
$$ \nu_{T,K,L}:\  Hom_{Perv(X)}(T*K,L) \ \cong\ Hom_{Perv(X)}(T,K^\vee*L) \ .$$ For $T_i\in
Perv(X)$ now suppose $T_i*K\in Perv(X)$. Then 2) for a morphism $f\in Hom_{Perv(X)}(T_2,T_1)$
the induced diagram
$$ \xymatrix{ Hom_{Perv(X)}(T_1*K,L) \ \ar[d]\ar[r]^{\nu_{T_1,K,L}} & \ Hom_{Perv(X)}(T_1,K^\vee*L) \ar[d]\cr
Hom_{Perv(X)}(T_2*K,L) \ \ar[r]^{\nu_{T_2,K,L}} & \ Hom_{Perv(X)}(T_2,K^\vee*L) \cr  } $$ is
commutative. Furthermore, if the assumptions  $K_i\in Perv(X)$, $T*K_i\in Perv(X)$ and
$K_i^\vee*L\in Perv(X)$ hold for $i=1,2$, then 3) they hold for $K=K_1\oplus K_2$ such that $
\nu_{T,K_1\oplus K_2,L} = \nu_{T,K_1,L} \oplus \nu_{T,K_2,L}$.
\end{Corollary}

\bigskip\noindent
\underbar{Proof}: $Hom_{Perv(X)}(T*K,L)\cong {\cal H}^0((T*K)*L^\vee)^\vee_0 \cong {\cal
H}^0(T*(K*L^\vee))^\vee_0 \cong {\cal H}^0(T*(K^\vee*L)^\vee)^\vee_0 \cong
Hom_{Perv(X)}(T,K^\vee*L)$ by corollary \ref{zuvor}. These isomorphisms are functorial in $T$
provided the underlying convolutions $T*K$ are perverse. The isomorphisms obviously are
additive in $K$. $\quad \Box$

\bigskip\noindent
Assume $K,L,K^\vee*L\in Perv(X)$. Then for  $T=\delta_0$ (also denoted $T=1$ in the following)
the assumptions are fulfilled. Hence
$$ \nu_{1,K,L} : \ Hom_{Perv(X)}(K,L) = Hom_{Perv(X)}(1,K^\vee*L) \ .$$

 %Conversely by $\nu_{1,K,L}^{-1}$ to $\gamma: 1\to
%K^\vee$ one can associate $\gamma^\vee: K\to 1$ such that $${\beta^\vee}{}^{\vee}=\beta\ .$$

\bigskip\noindent
An important case is $T=K^\vee$ and $L=1$. Assume $K$ and $K^\vee*K\in Perv(X)$. Then the
isomorphism $\nu_{K^\vee,K,1}$ gives an evaluation map $ev_K:K^\vee *K\to 1$ via
$$ ev_K\in Hom_{Perv(X)}(K^\vee*K,1) \cong Hom_{Perv(X)}(K^\vee,K^\vee) \ni id_{K^\vee} \ .$$
By the partial functoriality stated in the last lemma one has under the assumptions $
T,K,L,T*K,K^\vee*L\in Perv(X)$ and $K^\vee*L*K\in Perv(X)$ a commutative diagram
$$ \xymatrix@+0,5cm{ g\in Hom_{Perv(X)}(T*K,L) \ar[r]^{\nu_{T,K,L}} &  Hom_{Perv(X)}(T,K^\vee*L) \ni f \cr
 ev\in Hom_{Perv(X)}(K^\vee*L*K,L)\ar[u]_{f^*} \ar[r]^{\nu_{K^\vee*L,K,L}} &  Hom_{Perv(X)}(K^\vee*L,K^\vee*L) \ni id \ar[u]^{f^*} \cr}
 $$
In other words, under the assumptions made, there exists for $g:T*K\to L$ a unique $f:T\to
K^\vee*L$ making the diagram
$$ \xymatrix@+0,5cm{ T*K \ar[d]_{f*id_K}\ar[dr]^g & \cr (K^\vee *L)*K \ar[r]^-{ev} & L \cr} $$
commutative. Compare \cite{DM} (1.6.1).

\bigskip\noindent
In special cases there exists a more direct way to describe $\nu_{T,K,L}$.

\begin{Lemma} \label{danach} For $K,L\in {}^p D^{\leq 0}(X)$ one has $ Hom_{D_c^b(X,\overline\Q_l)}(K*L^\vee,1)
\cong {\cal H}^0(K*L^\vee)^\vee_0$ via adjunction for the right arrow
$$ \xymatrix{ K*L^\vee \ar[rr]\ar[dr]_{{}^{st}\tau_{\geq 0}} & & \delta_0= i_{0*}(\overline\Q_l) \cr
& {\cal H}^0(K*L^\vee) \ar@{..>}[ur]_{\exists !} & \cr} $$
\end{Lemma}

\bigskip\noindent
\underbar{Proof}: By standard truncation
$Hom_{D_c^b(X,\overline\Q_l)}(K*L^\vee,\delta_0) = Hom_{Et(X)}({\cal
H}^0(K*L^\vee),\delta_0)$ using $\delta_0\in {}^{st}D^{\geq
0}(X,\overline\Q_l)$ and $K*L^\vee \in {}^{st}D^{\geq
0}(X,\overline\Q_l)$ (lemma \ref{van}.1). By definition the right
side is $Hom_{Et(X)}({\cal
H}^0(K*L^\vee),Ri_{0,*}\overline\Q_{l,0})$ where $i_0:\{0\}\to X$ is
the inclusion map. By adjunction this is $Hom_{Perv(k)}(i_0^*{\cal
H}^0(K*L^\vee),\overline\Q_l) = {\cal H}^0(K*L^\vee)_0^\vee$.

\bigskip\noindent
For $K\in Perv(X)$ this gives a more explicit description of the morphism
$$ ev_K: K^\vee*K \cong K*K^\vee \ \longrightarrow 1 \ ,$$
which  by $Hom_{Perv(X)}(K,L)\cong {\cal H}^0(K*L^\vee)^\vee_0 \cong
Hom_{D_c^b(X,\overline\Q_l)}(K*L^\vee,1)$ (combining corollary \ref{zuvor} and lemma
\ref{danach}) corresponds to the identity map $id_K\in Hom_{Perv(X)}(K,K)$.

\section{Rigidity}\label{rig}

\bigskip\noindent
In a monoidal category $P$ in the sense of \cite{WF}, p.113), as for example in
$(D_c^b(X,\overline\Q_l),*)$, an object $K$ is called rigid (or dualizable in the sense of
\cite{De} (0.1.4)), if there exists an object $L$ together with morphisms $\delta$ and $ev$
$$ \delta: 1\longrightarrow K*L \quad , \quad ev: L * K \longrightarrow 1 \ ,$$
such that the following two rigidity statements hold: The composite morphisms
$$\xymatrix{ K=1*K \ar[r]^-{\delta*id} & (K*L)*K = K*(L*K) \ar[r]^-{id*ev}  &
K*1=K \cr}$$
$$\xymatrix{ L=1*L  & (L*K)*L = L*(K*L) \ar[l]_-{ev*id} &
L*1=L\ar[l]_-{id*\delta} \cr}\ $$ are the identity  $id_K$ respectively the identity $id_{L}$.
Then $L$ is called the dual of $K$. If such a triple $(L,\delta,ev)$ exists, it is unique up to
a unique isomorphism (\cite{WF}, p.113). If all objects of $P$ are rigid, the category is
called rigid.

%To define $ev_K$ for $K=\oplus_{i=1}^n K_i$ and given $ev_{K_i}$ put $$ev_K:
%K^\vee*K =\bigoplus_{i,j=1}^n K_i^\vee*K_j \overset{pr}{\to} \bigoplus_{i=1}^n K_i^{\vee}*K_i
%\to 1\ ,$$ where the last morphism is the sum $\sum_{i=1}^n ev_{K_i}$ of the morphisms
%$ev_{K_i}$. If each of the $ev_{K_i}$ satisfies the first rigidity condition, then the same
%holds for $ev_K$. To show this reduces to consider diagonal maps only
%$$\bigoplus_i K_i \overset{\oplus \delta_{K_i} *id}{\longrightarrow} \bigoplus_i (K_i*K_i^\vee)
%* \bigoplus_j K_j
%$$ $$ = \bigoplus_i K_i* \bigoplus_j (K_i^\vee *K_j) \overset{id*pr}{\to} \bigoplus_i K_i*(K_i^\vee*K_i)
%\overset{id*\sum ev_{K_i} }{\longrightarrow} \bigoplus_i K_i \ .$$

\bigskip\noindent
Let us consider the category $P$ of pure complexes in $D_c^b(X,\overline\Q_l)$ of geometric
origin. Since the associativity constraints are strict, for some given $ev_K: K^\vee *K\to 1$
we may try to define $Y$ to be $K^\vee$ and $\delta_K$ to be $\mu_K\cdot (ev_K)^\vee$ for some
$\mu_K\in\overline\Q_l^*$ using the contravariant functor $(.)^\vee$ and the evaluation maps
$ev_K$ already defined at the end of section \ref{Begin}. We obtain a commutative diagram $\bf
(*)$ defining $\delta_K$ via the identification $\lambda_K:K\cong K^\vee{}^\vee$ and $1^\vee=1$
$$ \xymatrix@+0,3cm{ 1 \ar[r]^-{ev_K^\vee} & K^\vee{}^{\vee} * K^\vee \cr
1 \ar[u]^{\mu_K}_\sim\ar[r]^-{\delta_K} & K*K^\vee \ar[u]_{\lambda_K*id_K}^-\sim \cr} \ $$ If
the first rigidity formula holds for $(K^\vee,ev_K^\vee,ev_K)$, then also the second: The last
diagram $\bf (*)$ and its dual $\bf(*)^\vee$, the dual of the first rigidity formula together
with $(\lambda_K)^\vee=(\lambda_{K^\vee})^{-1}$ (lemma \ref{eas}) show, that the diagram
$$\xymatrix{
   & (K^\vee{}^\vee{}^\vee*K^\vee{}^\vee)*K^\vee \ar[dl]_{\mu_K\cdot (ev_K)^\vee{}^\vee * id\ \ }\ar[d]^{((\lambda_K)^\vee*id)*id}_{(*)^\vee\quad }  &   &    \cr
K^\vee=1*K^\vee  & (K^\vee*{K^\vee}{}^\vee)*K^\vee \ar@{=}[r]\ar[l]_-{\delta_K^\vee*id} &
K^\vee*({K^\vee}{}^\vee*K^\vee)
 & K^\vee*1=K^\vee\ar[l]_-{id*ev_K^\vee} \ar[dl]^{id*\delta_K\cdot \mu_K^{-1}}_{(*)\quad } \cr
   & (K^\vee*K)*K^\vee \ar[lu]^{\mu_K\cdot ev_K*1} \ar@{=}[r]\ar[u]_{(id*\lambda_K)*id}  & K^\vee *(K *K^\vee) \ar[u]^{id*(\lambda_K*id)}   & \cr
}\ $$ is commutative. The middle line composed gives the identity morphism of $K^\vee$. Hence
the lower morphisms composed also give the identity morphism. This shows, that it is enough to
find a morphism $ev_K$, which  together with $\delta_K =\mu_K\cdot (ev_K)^\vee$ satisfies the
first rigidity property.

\bigskip\noindent
Now we make the

\bigskip\noindent
{\bf Assumption}: $K\in Perv(X)$ is simple of geometric origin.

\bigskip\noindent
Then $End_{Perv(X)}(K)=\overline\Q_l$. By lemma \ref{zuvor} und \ref{danach}
$$dim\bigl( Hom_{D_c^b(X,\overline\Q_l)}(K^\vee*K,1)\bigr) =1\ .$$
Hence up to a scalar there is a unique possible choice for $ev_K$.

\bigskip\noindent
Step 1.  Therefore, by the irreducibility of $K$ the first rigidity property holds for a
suitable generator of
$$Hom_{D_c^b(X,\overline\Q_l)}(K^\vee*K,1)$$ if and only if the composite map $K\to K$, given in the first
rigidity formula, is not the zero map. Whether this holds can be checked from the induced
endomorphism $H^\bullet(ev_K): H^\bullet(X,K)\to H^\bullet(X,K)$  provided $H^\bullet(X,K)$
does not vanish. That this can happen is evident in the case $K=\delta_\chi$ for a nontrivial
character $\chi$. However in the case where $K^\vee*K\in Perv(X)$, which will be relevant for
us later, the decomposition theorem implies $H^\bullet(X,K)\neq 0$, since $1=\delta_0$ then is
a direct summand of the pure complex $K^\vee*K$ in $Perv(X)$. Hence $H^\bullet(X,K)\neq 0$ by
the K\"unneth formula. So we make the additional

\bigskip\noindent
{\bf Assumption}: $H^\bullet(X,K)\neq 0$.

\bigskip\noindent
We then want to determine the map
$$ H^\bullet(X,K^\vee*K) \overset{H^{\bullet(ev_K)}}{\longrightarrow} H^\bullet(X,\delta_0) \ $$
(up to a nonvanishing scalar in $\overline\Q_l$) induced by a fixed generator
$$ ev_K : K^\vee*K \to \delta_0 \ $$
of the one dimensional $\overline\Q_l$-vectorspace
$Hom_{D_c^b(\overline\Q_l)}(K^\vee*K,\delta_0)$.

\bigskip\noindent
Step 2. There exist functorial isomorphisms
$$ i_K: H^\bullet(X,K^\vee) \cong H^\bullet(X,K)^\vee \ $$
defined by the composite morphism
$$ \xymatrix{ H^\bullet(X,K^\vee)\ar@{=}[r] & H^\bullet(X,(-id_X)^*(DK)) & H^\bullet(X,DK)\ar[l]_-{(-id_X)^*}^-{\sim}
\ar[r]^{PD}_\sim & H^\bullet(X,K)^\vee\cr} \ $$ with the last
isomorphism given by Poincare duality.  For $K,L\in
D_c^b(X,\overline\Q_l)$ there exist isomorphisms $c_{K,L}$
$$c_{K,L}: H^\bullet(X,K*L)\cong H^\bullet(X,K)\otimes^s H^\bullet(X,L)$$ defined as the
composite of the following isomorphisms
$$ H^\bullet(X,K*L)=H^\bullet(X,Ra_*(K\boxtimes L)) \cong H^\bullet(X\times X,K\boxtimes L)
\overset{\sim}{\longleftarrow} H^\bullet(X,K)\otimes^s H^\bullet(X,L) \ .$$ The second morphism
comes from Leray, the third from the K\"unneth formula. These functorial isomorphisms make
$H^\bullet(X,-)$ into a tensor functor in the sense of \cite{D}.2.7.

\begin{Lemma}\label{ACU}  The isomorphisms $c_{K,L}$ define a tensor functor ($\otimes$-functor
ACU)
$$H^\bullet(X,-) :\ (D_c^b(X,\overline\Q_l),*) \to (Vec_{\overline\Q_l}^{\pm},\otimes^s) $$
from $(D_c^b(X,\overline\Q_l),*)$  to the tensor category
$(Vec_{\overline\Q_l}^{\pm},\otimes^s)$ of $\overline\Q_l$-super-vectorspaces. \end{Lemma}

\bigskip\noindent
\underbar{Proof}: The K\"unneth theorem implies $(c_{A, B}\otimes^s id)\circ c_{A*B,C} =
(id\otimes^s c_{B,C})\circ c_{A,B*C}$, i.e compatibility with associativity. Compatibility of
the commutativity constraints can be reduced to the commutativity of the diagram
$$ \xymatrix@+0.3cm{ H^\bullet(X\times X, K\boxtimes L) \ar[rr]^{H^\bullet(\psi_{pr_1^*(K),pr_2^*(L)})} & & H^\bullet(X\times X, L\boxtimes K)
\cr H^\bullet(X,K) \otimes^s H^\bullet(X,L)\ar[u]^\sim_{Kuenneth} \ar[rr]^{superconstraints} &
& H^\bullet(X,L) \otimes^s H^\bullet( X, K) \ar[u]_\sim^{Kuenneth} \cr} $$ For this, one can
assume $K$ and $L$ to be replaced by injective resolutions. Then this follows immediately from
the definition of the complex constraints $\psi_{pr_1^*(K),pr_2^*(L)})$ and the K\"unneth
isomorphism. The case of unit object is trivial.

\bigskip\noindent
In particular we now have the functorial isomorphisms
$$ (i_K\otimes id)\circ c_{K^\vee,K}: H^\bullet(X,K^\vee*K) \cong H^\bullet(X,K)^\vee
\otimes^s H^\bullet(X,K)\ .$$

\bigskip\noindent
Step 3. For $f(x,y)=y-x$ $$ \xymatrix{ X\times X \ar[d]_{(-id_X)\times id}\ar[r]^f & X
\ar@{=}[d]\cr X\times X \ar[r]^a & X \cr } $$ we have an associated commutative diagram
$$ \xymatrix{  & H^{\bullet}(X,K^\vee *L) \ar@{=}[dl] \ar@{=}[dr] & \cr
H^{\bullet}(X,Ra_*(K^\vee\boxtimes L))\ar[d]_{\sim}\ar[rr]^{(-id_X)^*\times id} & &
H^{\bullet}(X,Rf_*(D(K)\boxtimes L))\ar[d]_{\sim}^{Leray} \cr H^{\bullet}(X\times
X,K^\vee\boxtimes L) \ar[rr]^{(-id_X)^*\times id} & & H^{\bullet}(X\times X,D(K)\boxtimes L)
\cr H^{\bullet}(X,K^\vee)\otimes^s H^{\bullet}(X,L)\ar[rr]^{(-id_X)^*\otimes id} \ar[u]^{\sim}
&  & H^{\bullet}(X,D(K))\otimes^s H^{\bullet}(X,L) \ar[u]^{\sim}_{Kuenneth} \cr } $$ If we
compare $(i_K\otimes id)\circ c_{K^\vee,K}$ with the map, which is obtained from going down the
left side (which is $c_{K^\vee,K}$) and then along the bottom line  of the large diagram, we
obtain an alternative description of $(i_K\otimes id)\circ c_{K^\vee,K}$ as the composite
$$ H^\bullet(X,K^\vee*K)=H^\bullet(X,Rf_*(D(K)\boxtimes K))  $$ $$ \cong\ H^\bullet(X\times X,D(K)\boxtimes K) \overset{\sim}{\longleftarrow}
H^\bullet(X,D(K))\otimes^s H^\bullet(X,K) $$ $$ \cong\ H^\bullet(X,K)^\vee \otimes^s
H^\bullet(X,K) \ .$$ The second isomorphism comes from Leray, the third from the K\"unneth
formula, the last from Poincare duality.

\bigskip\noindent
{\bf Our aim}: We would like to know, whether the tensor functor $H^\bullet(X,-) :
D_c^b(X,\overline\Q_l) \to Vec_{\overline\Q_l}^\pm$ (defined using the above identification
isomorphism), makes the following diagram commutative $$ \xymatrix{ K^\vee * K
\ar[d]_{ev_K}\ar@{..>}[rr] & & H^\bullet(X,K)^\vee \otimes^s H^\bullet(X,K)
\ar[d]^{H^\bullet(ev_K)\circ c_{K^\vee,K}^{-1}\circ (i_K^{-1}\times id)} \cr \delta_0
\ar@{..>}[rr]& & H^\bullet(X,\delta_0) =\overline\Q_l \cr}
$$  i.e. maps $ev_K$  for a suitable choice of $ev_K$ to the canonical
evaluation map  $eval$ of $(Vec_{\overline\Q_l}^\pm,\otimes^s)$
$$ H^\bullet(X,ev_K) \ \ \overset{?}{=} \ \ eval \circ (i_K\times id) \circ c_{K^\vee,K} \ .$$

\bigskip\noindent
Step 4. Notice
$H^\bullet(X,\delta_0)=H^0(\{0\},i_0^*(\overline\Q_l)) =
\overline\Q_l$ for the closed immersion $i_0: \{0\}\hookrightarrow
X$. Hence  $H^\bullet(ev_K)$ factorizes over the cohomology in
degree zero
$$ \xymatrix{ H^\bullet(X,K^\vee * K) \ar[d]_{H^\bullet(ev_K)}\ar[rr]^{pr^0} &   &  H^0(X,K^\vee * K) \ar[d]^{H^0(ev_K)} \cr
H^\bullet(X,\delta_0) \ar@{=}[rr]& & H^0(X,\delta_0) =\overline\Q_l \cr}
$$

\bigskip\noindent
Step 5. For any $f:Y\to X$ and complex map $\varphi: L\to L'$ one has a commutative diagram
$$ \xymatrix{ H^\bullet(X,L)\ar[d] \ar[r] & H^\bullet(Y,f^*(L)) \ar[d] \cr H^\bullet(X,L')
\ar[r] & H^\bullet(Y,f^*(L')) \cr} $$ for the canonical morphisms $H^\bullet(X,L)\to
H^\bullet(Y,f^*(L))$, denoted $res_L$ in the case of closed immersions $f:Y\hookrightarrow X$.
Applied for a complex map $\varphi: L\to L'=\delta_0 = i_{0*}(\overline\Q_l)$ we get -- with
$L=K^\vee *K$ in mind -- the following commutative diagram
$$ \xymatrix{ H^0(X,L) \ar[d]_{H^0(\varphi)}\ar[r]^-{res} & H^0(\{0\}, i_0^*(L)) \ar[d]^{i^*_0(\varphi)}\cr
H^0(X,\delta_0) \ar[r]^-{res} & H^0(\{0\}, i_0^*(\delta_0)) = \overline\Q_l \cr}
$$ where $i_{0}^*(\delta_0)=i_{0}^*i_{0*}\overline\Q_l=\overline\Q_l$.

\bigskip\noindent
Step 6. We compute the restriction maps of step 5 for $L=Rf_*(DK\boxtimes K)$ via the proper
base change using the diagram
$$ \xymatrix{ X\times X \ar[r]^f &  X \cr
\Delta_X \ar[r]^{\overline f}\ar@{_{(}->}[u] & \{0\}.\ar@{_{(}->}[u]_{i_0} \cr } $$ For the
diagram of step 5 (see also step 3)
$$ \xymatrix{ H^0(X,K^\vee *K)= H^0(X,Rf_*(DK\boxtimes K))\ar[d]_{H^0(ev_K)} \ar[r]^-{res} &
H^0(\{0\},i_0^*(Rf_*(DK\boxtimes K))) \ar[d] \cr H^0(X,\delta_0) \ar[r]^{res} &
H^0(\{0\},\overline\Q_l) \cr } $$ this expresses $ i_{0}^*(Rf_*(D(K)\boxtimes K))$ by $
R\overline f_*(D(K)\boxtimes K\vert_{\Delta_X}) = R\overline f_*(D(K)\otimes^L K)$, where the
latter is a complex over $\{0\}=Spec(k)$, hence can be identified with $
H^\bullet(\Delta_X,D(K)\otimes^L K)$, making  the following  diagram commutative
$$ \xymatrix{ H^0(X\times X, D(K)\boxtimes K) \ar[r]^-{res} & H^0(\Delta_X, D(K)\otimes^L K) \cr
H^0(X, Rf_*(DK\boxtimes K))\ar[u]_\sim^{Leray} \ar[r]^-{res} &
H^0(\{0\},i_0^*(Rf_*(D(K)\boxtimes K))),\ar[u]^\sim_{Leray} \cr}
$$ with the vertical isomorphisms being of Leray type. Since the K\"unneth isomorphism is defined by cup-product, and since
cup-products are functorial with respect to pullbacks (this follows immediately from their
definition given e.g. in \cite{Mi2}, p.167ff), we have a commutative diagram
$$ \xymatrix{ \bigoplus_{i\in\Z} H^i(X,D(K)) \otimes^s H^{-i}(X,K) \ar[d]_{Kuenneth}^\sim\ar[r]^-{\cup} & H^0(\Delta_X, D(K)\otimes^L K)\ar@{=}[d] \cr
H^0(X\times X, D(K)\boxtimes K) \ar[r]^-{res} & H^0(\Delta_X, D(K)\otimes^L K) \cr}
$$
Notice, for the closed immersion $i_X: \Delta_X\to X\times X$ and the map $res=i_X^*$ this uses
$i_X^*(pr_1^*(x)\cup pr_2^*(y)) = i_X^*(pr_1^*(x))\cup i_X^*(pr_2^*(y)) = x\cup y$, and
$Kuenneth(x\otimes^\varepsilon y) = pr_1^*(x)\boxtimes pr_2^*(y)$.

%\footnote{\cite{V}, p. 210 states the conjecture to be true without proof. A proof should
%follow from \cite{SGA45}, section 1.2, p.155 together with \cite{SGA}, lemma 3.2.3, p.583 with
%$L=\Lambda$ (constant sheaf) combined with the functoriality of cup-products under pullback
%(see \cite{Mi2}, p.167ff).}
%
%\bigskip\noindent
%Step ?. In lemma \ref{} we have seen, that for perverse sheaves $K\in Perv(X)$ any morphism
%$ev_K$ factorizes uniquely
%$$ \xymatrix{ K^\vee * K \ar[rr]^{ev_K}\ar[dr]_{{}^{st}\tau_{\leq 0}} & & \delta_0= i_{0*}(\overline\Q_l) \cr
%& {\cal H}^0(K^\vee*K) \ar@{..>}[ur]_{\exists !} & \cr} $$

\bigskip\noindent
%Via the description of $Hom_{D_c^b(X,\overline\Q_l)}(K^\vee
%* K,\delta_0)$ at the end of section \ref{VT} it is clear, that the following morphism generates
%$Hom_{D_c^b(X,\overline\Q_l)}(K^\vee * K,\delta_0)$: Let $K=\delta_E$ with irreducible support
%$Y$ of dimension $d=dim(Y)$, and $j: U\hookrightarrow Y$ open dense such that $E\vert_U$ is an
%irreducible smooth etale $\overline\Q_l$-sheaf. Consider $f:X\times X\to X$ defined by
%$f(x,y)=y-x$. Then
%$$ K^\vee*K = Rf_*(DK\boxtimes K) \ .$$
%The space $Hom_{D_c^b(X,\overline\Q_l)}(K^\vee * K,\delta_0)$ was described in section \ref{VT}
%via the restriction $i_0^*$. This restriction can be computed by the proper basechange theorem.
%Since the $f^{-1}(\{0\})$ is the diagonal $\Delta_Y\cong Y $ defined by all $(y,y)\in Y\times
%Y$ we obtain
%$$ i_0^*(K^\vee*K) \ \cong \ H^\bullet(\Delta_Y,D(K)\otimes^L K) \ .$$
%Hence again, up to a scalar, there is no alternative for our choice
%$$ ev_K \in Hom_{D_c^b(X,\overline\Q_l)}(K^\vee * K,\delta_0) $$
%$$ = Hom_{\overline\Q_l}({\cal H}^0(i_0^*(K^\vee * K),\overline\Q_l)
% = Hom_{\overline\Q_l}(H^0(Y,D(K)\otimes^L K),\overline\Q_l) \ .$$
Let $K_X$ denote the dualizing sheaf on $X$. The Poincare duality map (\cite{SGA45}, section
1.2, p.155 together with \cite{SGA}, lemma 3.2.3, p.583 with $L=\Lambda$ the constant sheaf) it
is obtained by composing the cup product $H^{-i}(X,D(K)) \otimes^s H^{i}(X,K) \to H^0(\Delta_X,
D(K)\otimes^L K)$, with the map $S_K$ induced by  $DK \otimes^L K = R{\cal H}om(K,K_X)
\otimes^L K \to K_X $ on the zero-th cohomology group, and the canonical isomorphism
$H^0(X,K_X)=\overline\Q_l$. There are several equivalent possibilities to describe the cup
product (see \cite{Mi2}, V\S1, p.167ff and in particular prop.V.1.16). Composing all the recent
commutative diagrams, we obtain the following commutative diagram
$$ \xymatrix{ \bigoplus_{i\in\Z} H^i(X,K)^\vee \otimes^s H^{i}(X,K) \ar[r]^-{eval} &  \overline\Q_l\cr
\bigoplus_{i\in\Z} H^{-i}(X,D(K)) \otimes^s H^{i}(X,K) \ar[u]^{PD\otimes^s id}_\sim
\ar[d]_{Kuenneth}^\sim\ar[r]^-{\cup} & H^0(\Delta_X, D(K)\otimes^L
K)\ar@{->>}[u]^{S_K}\ar@{=}[d] \cr H^0(X\times X, D(K)\boxtimes K) \ar[r]^-{res} &
H^0(\Delta_X, D(K)\otimes^L K) \cr H^0(X, Rf_*(DK\boxtimes K))\ar[u]_\sim^{Leray} \ar[r]^-{res}
& H^0(\{0\},i_0^*(Rf_*(D(K)\boxtimes K))),\ar[u]^\sim_{Leray}\ar[dd]^{i_0^*(ev_K)} \cr
H^0(X,K^\vee *K)\ar@{=}[u]\ar[d]_{H^0(ev_K)} &
 \cr H^0(X,\delta_0) \ar[r]^{res}_\sim &
H^0(\{0\},\overline\Q_l) \cr } $$ By step 3 the left vertical map from top to bottom is
$H^\bullet(ev_K)\circ c_{K^\vee,K}^{-1}\circ (i_K^{-1}\times id)$. By lemma \ref{danach} the
map $i_0^*(ev_K)$ is nonzero, since $ev_K$ was chosen to be nonzero. So complete our aim it is
enough to show, that all $\overline\Q_l$-vectorspaces on the right side of the last diagram are
one dimensional. Since $S_K$ is surjective (Poincare duality is nondegenerate), and
$i_0^*(ev_K)$ is nonzero, this implies that all vertical maps on the right side are
isomorphisms, as well as the bottom horizontal map. Hence the composition of this map is given
by multiplication with a nonzero constant in $\overline\Q_l$. This constant therefore is 1, if
$ev_K$ was  chosen suitable.

\bigskip\noindent
Step 7. It remains to show $H^0(\Delta_X, D(K)\otimes^L K) \cong\overline\Q_l$. Let
$K=i_*(\delta_E)$ for $i:Y\hookrightarrow X$ be the support of $K$. Notice $K=i_*(\delta_E)$
and $D(K)=\delta_{\tilde E}\in Perv(Y)$, where $\tilde E\vert_U \cong (E\vert_U)^*$ holds on a
smooth open dense irreducible subset $U$ of the support $Y$ of $K$. $E\vert_U$ is an
irreducible smooth $\overline\Q_l$-sheaf on $U$. Let $d$ denote the dimension of $Y$. We claim
$$ H^0(\Delta_Y,D(K)\otimes^L K) \cong H^{0}(Y,\lambda_{\tilde E} \otimes \lambda_E) \cong
H^{2d}_c(U,E^*\otimes E) \overset{\sim}{\longrightarrow} H_c^{2d}(U,\overline \Q_{l,U})
\overset{Sp}{\longrightarrow} \overline\Q_l \ ,$$ where $(E\vert_U)^*\otimes E\vert_U \to
\overline\Q_{l,U}$ is the adjoint of $id\in Hom_{Et}(E\vert_U,E\vert_U)$, and where $Sp$
denotes the trace isomorphism. The first isomorphism either from section \ref{VT}, or from the
spectral sequence exploited in section \ref{help}, the second from excision and standard
vanishing theorems, the third from the vanishing theorems of etale cohomology implied by
Poincare duality on $U$, and the last from the trace isomorphism. We have the following
commutative diagram
$$ \xymatrix{  H^0(X, D(K) \otimes^L K) \ar[r] & H^0(X,K_X) \ar[r]^-\cong & \overline\Q_l \cr
 H^0(Y, D(\delta_E) \otimes^L \delta_E) \ar[u]^\sim\ar[r] & H^0(Y,K_Y)\ar[u] \ar[r]  &\overline\Q_l\ar@{=}[u] \cr
  H^0_c\bigl(U, (E\vert_U)^*[d] \otimes E\vert_U[d]\bigr) \ar[u]^\sim\ar[d]_\sim\ar[r] & H^0_c(U,K_U)\ar[u] \ar[r]^-\cong & \overline\Q_l \ar@{=}[u]\cr
   H^{2d}_c(U, \overline\Q_l ) \ar[rr]^-{Sp}_-\sim &  & \overline\Q_l \ar@{=}[u]\cr
 } $$
The morphisms on the right side are the trace isomorphisms. Let $j: U\hookrightarrow Y$ and
$i=i_Y:Y\to X$ denote the inclusions, and $s_X:X\to Spec(k)$ and $s_Y:Y\to Spec(k)$ the
structure morphisms. Then for a complex $L$ with support in $Y$, i.e. $L\in
i_*(D_c^b(Y,\overline\Q_l)$, we have $i^*(L)=i^!(L)$ and the adjunction morphisms $i_*i^!(L)
\to L$ is an isomorphism. Hence $H^\nu(X,L)= H^\nu(Y,i^!(L))=H^\nu(Y,i^*(L))$. Notice
$i^*(A\otimes^L B)=i^*(A)\otimes^L i^*(B)$ and $L=K\otimes^L D(K) \in
i_*(D_c^b(Y,\overline\Q_l)$. Hence $i^!(K\otimes^L D(K))=i^!(K\otimes^L R{\cal H}om(K,K_X)) =
i^*(K\otimes^L R{\cal H}om(K,K_X))= i^*(K)\otimes^L i^*(R{\cal H}om(K,K_X)) = \delta_E
\otimes^L i^!(R{\cal H}om(K,K_X))$, since $K=i_*(\delta_E)$ and $D(K)=i_*(D(\delta_E))$. Now
use $i^!(R{\cal H}om(A,B)) = R{\cal H}om(i^*(A),i^!(B))$ and $K_X=s_X^!(\overline\Q_l)$ and
$K_Y= s_Y^!(\overline\Q_l) = i^!s_X^!(\overline\Q_l) = i^!(K_X)$ to obtain $i^!(K\otimes^L
D(K))= \delta_E \otimes^L D(\delta_E)$ and $K_Y=i^!(K_X)$. Therefore, if we apply $i^!$ to the
evaluation $a: K\otimes^L D(K) \to K_X$, we obtain the evaluation $i^!(a): \delta_E \otimes^L
D(\delta_E) \to K_Y$. This gives the upper left commutative diagram induced by $i^!$ using
 $K,D(K)\in i_*(D_c^b(Y,\overline\Q_l)$. The second commutative diagram below
is obtained similarly by the restriction $j^*=j^!$ and the natural adjunction map
$Rj_!j^*(L)\to L$. Notice, as explained above, that this map induces an isomorphism
$H_c^0(U,L\vert_U)\to H^0(Y,L)$ on the first cohomology group for $L=\delta_E\otimes^L
D(\delta_E)$. Similarly $H_c^0(U,K_U) \to H^0(Y,K_Y)$ by $K_U=j^!(K_Y)=j^*(K_Y)$. Commutativity
of the lowest part of the last diagram is obvious. The composed map $H_c^0(U,K_U) \to
H^0(Y,K_Y)\to H^0(X,K_X)$ is an isomorphism. To show this, observe that one can find an open
dense subset $V\subseteq X$, such that $V\cap Y\subseteq U$. Hence without restriction of
generality $V\cap Y=U$, and $U \subseteq V$ is a smooth pair. Then
$H_c^0(U,K_U)=H_c^{2d}(U,\overline\Q_l)$ and $ H^0(X,K_X)=H^{2g}(X,\overline\Q_l)=
H_c^{2g}(V,\overline\Q_l)$, and the isomorphism comes from the Gysin map via the purity theorem
$(i\vert_U)^!\overline\Q_{l,V} = \overline\Q_{l,U}[2d-2g]$. See \cite{SGA}, XVI, theorem 3.7.

%The remaining arrows are induced
%
%
%The first diagram commutes, since $D$ commutes with $(i_Y)_*$ for the inclusion of the closed
%subvariety $i_Y:Y\hookrightarrow X$ and $K=(i_Y)_*(\delta_E)$ for $\delta_E\in Perv(Y)$.

\bigskip\noindent
Together the last steps imply

\begin{Corollary}\label{eva}
Suppose $K\in Perv(X)$ is an irreducible perverse sheaf, such that $H^\bullet(X,K)\neq 0$ (e.g.
$K^\vee*K\in Perv(X)$). Then for a suitable choice of $ev_K$ in the one dimensional
$\overline\Q_l$-vectorspace $Hom_{Perv(X)}(K^\vee*K,\delta_0)$
$$H^\bullet(X,ev_K) =  eval \circ (i_K\times id) \circ c_{K^\vee,K} $$ holds.
\end{Corollary}

\bigskip\noindent
By step 7 and the last part of step 6 we also obtained

\bigskip\noindent
\begin{Theorem}\label{th1}\label{Lefschetz} For  an irreducible perverse sheaf $K$ on a projective variety $Y$ the
diagram
$$ \xymatrix{ H^\bullet(Y\times Y,D(K) \boxtimes K) \ar[d]_{\mu}^\sim\ar[r]^-{res} & H^\bullet(\Delta_Y, D(K)\otimes^L K) \ar[d]^{S_K\ \circ\ pr^0}\cr
H^\bullet(Y,K)^\vee \otimes^s H^\bullet(Y,K) \ar[r]^-{eval} & \overline\Q_l \cr} $$ commutes,
where the left vertical isomorphism $\mu$ is composition $$ H^\bullet(Y\times Y,D(K) \boxtimes
K) \cong H^\bullet(Y,D(K)) \otimes^s H^\bullet(Y, K) \cong H^\bullet(Y,K)^\vee \otimes^s
H^\bullet(Y, K) $$ of the K\"unneth and the Poincare duality isomorphism.
\end{Theorem}

%\begin{Theorem}\label{th1}
%For projective smooth $Y$ and $K=\delta_Y$ the conjecture is true.
%\end{Theorem}

\bigskip\noindent
Suppose $Y$ projective and smooth and $K=\delta_Y$. Then in fact  this is a version of the
Grothendieck-Verdier fixed point formula \cite{SGA5}, p.101, \cite{SGA45}, cycle, section 2.3
and also \cite{Mi2}, theorem 12.3 or \cite{FK}, p.155. In section \ref{SP} we will use this
description of Poincare duality in terms of nondegeneracy of the the cup-product (theorem
\ref{th1}) for $K=\delta_Y$, where $Y$ is a smooth projective curve $Y$.

\bigskip\noindent
Step 8. We now show the dual version of corollary \ref{eva} by reducing it to corollary
\ref{eva}. For any complexes $A,B\in D_c^b(X,\overline\Q_l)$ and complex maps $\varphi: A\to B$
there
exists a commutative diagram %for the morphism (left vertical arrow of the
%diagram), which can be deduced from relative Poincare duality for the proper smooth morphism
%$f:X\times X\to X$ and $(-id_X)\circ f= f\circ (-id_X \times -id_X)$ and Poincare duality for
%$X$:
$$ \xymatrix{ H^\bullet(X,(A*B)^\vee) \ar[r]^-{i_{A*B}}_-\sim & H^\bullet(X,A*B)^\vee  & (H^\bullet(X,A) \otimes^s
H^\bullet(X,B))^\vee \ar[l]_-{(c_{A,B})^\vee}^-\sim \ar@{=}[d]\cr H^\bullet(X,A^\vee * B^\vee)
\ar[u]^-\sim \ar[r]^-{c_{A^\vee,B^\vee}}_-\sim & H^\bullet(X,A^\vee) \otimes
H^\bullet(X,B^\vee) \ar[r]^-{i_{A}\times i_{B}}_-\sim & H^\bullet(X,A)^\vee \otimes^s
H^\bullet(X,B)^\vee \cr}
$$
For the other isomorphisms see step 2. In fact this corresponds to diagram 4.3.3.3 of
\cite{Sav}, since $H^\bullet(X,-)$ is a tensor functor by lemma \ref{ACU}.

\bigskip\noindent
Now assume, that the diagram
$$ \xymatrix{ H(K^\vee*K) \ar[d]_-{H(ev_K)} \ar[r]^-{c_{K^\vee,K}} &
H(K^\vee) \otimes^s H(K) \ar[r]^-{i_K\otimes^s id} & H(K)^\vee \otimes H(K)\ar[dll]^-{eval} \cr
H(1) & & \cr} $$ commutes. We abbreviated $H^\bullet(X,-)$ by $H(-)$ for simplicity. Then the
dual of this diagram becomes the right lower sub-diagram of
$$ \xymatrix{ H(K*K^\vee)\ar[d]_-{H(\lambda_K*id)} \ar[r]^-{c_{K,K^\vee}} &
H(K)\otimes^s H(K^\vee)\ar[d]_-{H(\lambda_K)\otimes^s id} \ar[rr]^-{id\otimes^s i_K} & &
H(K)\otimes^s H(K)^\vee \ar[dd]^-{can\otimes^s id} \cr H(K^\vee{}^\vee*K^\vee)
\ar[r]^-{c_{K^\vee{}^\vee,K^\vee}} & H(K^\vee{}^\vee)\otimes^s H(K^\vee)
\ar[r]^-{i_{K^\vee}\times i_K} & H(K^\vee)^\vee \otimes^s H(K)^\vee  & \cr H((K^\vee*K)^\vee)
\ar[u]^\sim \ar[r]^-{i_{K^\vee*K}} & H(K^\vee*K)^\vee  & H(K^\vee)^\vee \otimes^s H(K)^\vee
\ar@{=}[u]\ar[l]_-{(c_{K^\vee,K})^\vee} & H(K)^\vee{}^\vee \otimes^s H(K)^\vee
\ar[l]_-{(i_{K})^\vee\otimes^s id } & \cr H(1^\vee) \ar[u]^{H(ev_K^\vee)} \ar@{=}[r] &
H(1)^\vee \ar[u]_{H(ev_K)^\vee} \ar[urr]_-{eval^\vee} & & \cr}
$$
Using the first diagram of lemma \ref{compdia} (left bottom), the compatibility diagram stated
above (middle) and the second compatibility diagram of lemma \ref{compdia} (right upper
sub-diagram) together with functoriality (upper left diagram) we obtain
$$ (eval)^\vee = (id\otimes i_K) \circ c_{K,K^\vee} \circ
H(\delta_K) \ .$$ Ignoring the identification isomorphisms, this can be simplified to the
statement: $H(ev_K)=eval$ implies $H(\delta_K) = eval^\vee$.

\bigskip\noindent Step 9. By functoriality
and $(c_{A, B}\otimes^s id)\circ c_{A*B,C} = (id\otimes^s c_{B,C})\circ c_{A,B*C}$

\bigskip\noindent
{\scriptsize
$$ \xymatrix{
H(1*K)\ar@{=}[d]^{c_{1,K}} \ar[r]^-{H(\delta_K)*id} & H((K*K^\vee)*K)\ar[d]_{c_{K*K^\vee,K}}
\ar@{=}[r] & H(K*(K^\vee *K)) \ar[d]_{c_{K,K^\vee*K}} \ar[r]^-{id*H(ev_K)} &
H(K*1)\ar@{=}[d]_{c_{K,1}} \cr H(1)\otimes^s H(K) \ar[r]^-{H(\delta_K)\otimes^s id} &
H(K*K^\vee)\otimes^s H(K) \ar[d]_{c_{K,K^\vee}*id}  & H(K)\otimes^s H(K^\vee *K) \ar[d]_{id*
c_{K^\vee,K}} \ar[r]^-{id \otimes^s H(ev_K)} & H(K)\otimes^s H(1) \cr  & (H(K)\otimes^s
H(K^\vee))\otimes^s H(K)\ar[d]_{id*i_K*id} \ar@{=}[r] & H(K)\otimes^s (H(K^\vee) \otimes^s
H(K))\ar[d]_{id*i_K*id} &  \cr  1 \otimes^s H(K) \ar@{=}[uu] \ar[r]^-{eval^\vee \otimes^s id} &
(H(K)\otimes^s H(K)^\vee)\otimes^s H(K) \ar@{=}[r] & H(K)\otimes^s (H(K)^\vee \otimes^s H(K))
\ar[r]^-{id\otimes^s eval} & H(K)\otimes^s 1 \ar@{=}[uu]\cr}
$$}
we get:
$$H((id_K*ev_K)\circ (\delta_K * id_K)) = (id\otimes^s
eval)\circ (eval^\vee \otimes^s id) = id$$ is the identity of $H^\bullet(X,K)$, since finite
dimensional super-vectorspaces are rigid in $Vec_{\overline\Q_l}^\pm$.

\bigskip\noindent
By step 1 this implies, that $K\in Perv(X)$ is rigid with dual $(K^\vee,ev_K^\vee,ev_K)$, if
the conjecture above holds for $K$ and if $H^\bullet(X,K)\neq 0$.

%\begin{Corollary} \label{sppo} For smooth proper irreducible subvarieties $Y$ of $X$ the perverse sheaves $\delta_Y$
%are rigid. \end{Corollary}

\begin{Corollary} \label{sppo}
If $K\in Perv(X)$ is simple of geometric origin with $H^\bullet(X,K)\neq 0$, then $K$ is rigid.
\end{Corollary}

\bigskip\noindent
Finally for an additive, $\overline\Q_l$-linear karoubienne category $P$ with the additional
properties listed in \cite{De}, 1.2. (examples will be the categories ${\cal BN}$ and ${\cal
T}(X)$ later) the following holds (\cite{De} 1.15)

\begin{Lemma}\label{thereafter} Direct sums, tensor products of rigid objects are rigid.
Direct summands of rigid objects are rigid.
\end{Lemma}

\bigskip\noindent
Notice, that for a rigid pure complex $K\in Perv(X)$ and a homomorphism $f:X\to Y$ between
abelian varieties over $k$ the direct image complex $Rf_*(K) \in D_c^b(X,\overline\Q_l)$ again
is rigid, hence all its summands by the decomposition theorem. Taking such direct images of
$\delta_Y$ as in corollary \ref{sppo}, and convolution products, direct sums and subfactors all
give rigid objects. In particular, this implies that the abelian categories ${\cal BN}$ and
${\cal T}(X)$ defined in section \ref{APP} and \ref{categor} are rigid.

\bigskip\noindent
\begin{Lemma}\label{triangel} Suppose $K_1\to K_2 \to K_3 \to K_1[1]$ is a distinguished
triangle in $D_c^b(X,\overline\Q_l)$, such that $K_1$ and $K_3$ are rigid. Then also $K_2$ is
rigid. Translates of rigid objects are rigid.
\end{Lemma}

\bigskip\noindent
\underbar{Proof}: For the first statement see \cite{Le}, Lemma IV.1.2.3. The second statement
is obvious. \qed

\bigskip\noindent
Notice, that $D_c^b(X,\overline\Q_l)$ is the derived category of the abelian category $Perv(X)$
by a quite deep result of Beilinson \cite{Be}. The complexes $K\in D_c^b(X,\overline\Q_l)$ with
$H^\bullet(X,K)=0$ define a tensor ideal $N$, i.e. $H^\bullet(X,K)=0$ implies
$H^\bullet(X,K*L)=0$ for any $L\in D_c^b(X,\overline\Q_l)$. For a distinguished triangle
$(f,g,h): K_1\to K_2 \to K_3 \to K_1[1]$, for which $K_3\in N$ and the morphism $f:K_1\to K_2$
factorizes over an element $K_3'\in N$, the objects $K_1,K_2$ are also in $N$. Obviously $f$
induces an isomorphism $H^\bullet(X,f)$, which must be zero. Hence $H^\bullet(X,K_i)=0$ for
$i=1,2$. This shows, that $N$ defines a thick subcategory. This allows to pass to the quotient
monoidal category $D_c^b(X,\overline\Q_l)/N$. The canonical functor $D_c^b(X,\overline\Q_l) \to
D_c^b(X,\overline\Q_l)/N$ is a tensor functor. Beilinson's result stated above  and corollary
\ref{sppo} implies, that the full triangulated monoidal subcategory of complexes $K\in
D_c^b(X,\overline\Q_l)$, for which all simple perverse constituents of all perverse cohomology
sheaves ${}^pH^i(K)$ are of geometric origin, defines a rigid monoidal subcategory of
$D_c^b(X,\overline\Q_l)/N$. If $k$ is the algebraic closure of a finite field, then by
fundamental results of Lafforgue on the Langlands program every complex in
$D_c^b(X,\overline\Q_l)$ is mixed. Using Beilinson's results stated above and the weight
filtration of mixed perverse sheaves, lemma \ref{triangel} implies

\bigskip\noindent
\begin{Theorem}\label{rigid-monoidal} {If $k$ is the algebraic closure of a finite field, then the monoidal
category $D_c^b(X,\overline\Q_l)/N$ is a rigid monoidal category.} \end{Theorem}

\goodbreak
\section{Hard Lefschetz}
Assume that $K$ and $L$ are pure perverse sheaves of weights $w$ respectively $w'$. Then $K$
and $L$ are semisimple by the decomposition theorem. The hard Lefschetz theorem for the
projective morphism $a:X\times X\to X$ implies
$$ K*L = \bigoplus_A \ \bigoplus_{\nu=-\nu(A)}^{\nu(A)} A[2\nu](\nu) \ , $$
where the sum ranges over finitely many irreducible pure perverse sheaves $A\in Perv(X)$ of
weight $w+w'$. Hence in the Grothendieck group the class of $K*L$ becomes
$$ \sum_{A}  A \cdot (t^{-\nu(A)/2} + \cdots + t^{\nu(A)/2}) \ .$$
For $[n]_t= (t^{n/2}-t^{-n/2})/(t^{1/2}-t^{-1/2})$ the class of
$K*L$ is $\sum_{A}   [1+\nu(A)]_t\cdot A$. A priori there are the
following bounds
$$ \nu(A) \leq g-1\ ,$$ unless both $K$ and $L$ contain nontrivial constant perverse subsheaves.
This follows from \cite{KW} III.11.3, since $a:X\times X\to X$ is a smooth morphism with
geometrically connected fibers of dimension $g$. On the other hand
$$\nu(A)\leq dim\ supp(A)\ $$ because ${\cal H}^\nu(K*L)=0$  by lemma \ref{van} for
$\nu\geq 1$.

\goodbreak
\section{Computation of ${\cal H}^{-1}(K*L)$}\label{help}

\bigskip\noindent
Suppose $K=\delta_E$ and $L=\delta_F$ are irreducible perverse sheaves on $X$.  Since $${\cal
H}^{-1}(K*L)_z = H^{-1}(X,K\otimes^L \varphi_z^*(L))\ ,$$ the spectral sequence of section
\ref{stalks} has $E_2$-terms
$$\bigoplus_{m} E_2^{m-1,-m}  \ =\ \bigoplus_{a,b\leq 0} \ \ H^{a+b-1}\Bigl(X,{\cal
H}^{-a}(K)\otimes {\cal H}^{-b}(\varphi_z^*(L))\Bigr)\ .$$ In the
$E_\infty$-limit  ${\cal H}^{-1}(K*L)_z $ becomes a subquotient of
$\bigoplus_{m} E_2^{m-1,-m}$. We consider now the single terms. For
this assume $a\leq b$, otherwise switch the roles of $K$ and $L$.
Then, unless $a=dim\ supp(K)$, we have $$dim\ supp({\cal H}^{-a}(K))
\leq a-1$$ by the $IC$-sheaf property [FK] prop.III.9.3(4). Hence $
2\cdot dim\ supp({\cal H}^{-a}(K)) \leq 2(a-1) < a+b-1$, hence the
corresponding $E_2$-term of the spectral sequence is zero. Thus,
without restriction of generality,
$$a=dim\ supp(K) \leq b\leq dim\ supp(\varphi_z^*(L))\ .$$ Then
$supp({\cal H}^{-a}(K)\otimes {\cal H}^{-b}(\varphi_z^*(L)))$ is contained in $$Y=supp(K) \cap
supp ({\cal H}^{-b}(\varphi_z^*(L)))\subseteq supp(K) \cap supp(\varphi_z^*(L))\ .$$ Suppose
$Y$ is a proper subset of $supp(K)$.  Since $K$ is an irreducible perverse sheaf, $supp(K)$ is
irreducible. Hence $dim(Y)\leq dim\ supp(K) - 1 =a-1$. Again this implies $2dim(Y)= 2(a-1) <
a+b-1$, and the corresponding term of the spectral sequence vanishes. Also for $b\geq a+2$ once
more $2dim(Y)\leq  2a < a+b-1$ and the corresponding term vanishes. Hence the stalk ${\cal
H}^{-1}(K*L)_z$ vanishes except if
\begin{enumerate} \item $supp(K)= supp({\cal H}^{-a}(K))\subseteq supp ({\cal
H}^{-b}(\varphi_z^*(L)))\subseteq supp(\varphi_z^*(L))$, and
\item Either $b=a=dim \ supp(K)$, or $b=a+1=dim\ supp(K) +1$.
\end{enumerate}

\bigskip\noindent
\underbar{The subcase $b=a$}: Then,  unless $b=dim\ supp(L)$, the dimension of $supp\ {\cal
H}^{-b}(\varphi_z^*(L))$ is $<b$, hence $<a$. Therefore we can assume $a=supp(K)=
supp(\varphi_z^*(L))$ by dimension reasons. Then
$$ {\cal H}^{-1}(K*L)_z\neq 0 \Longrightarrow z\in S(K,L) \ .$$
The relevant $E_2$-term in the spectral sequence is
$$ H^{2a-1}\Bigl(supp(K),{\cal H}^{-a}(K)\otimes {\cal
H}^{-a}(\varphi_z^*(L))\Bigr) = H^{2a-1}\Bigl(supp(K), E\otimes \varphi_z^*(F)\Bigr) \ $$
${\cal H}^{-1}(K*L)_z$ in the $E_\infty$-limit is a quotient group of the cohomology group
$$ H^{2a-1}\Bigl(supp(K), E\otimes \varphi_z^*(F)\Bigr)\quad , \quad a=dim(supp(K))\ .$$
To visualize the different $E_2$-terms: All nonvanishing terms are contained in the lower
triangle of the square with coordinates $(0,0),(0,2a),(-2a,0),(-2a,2a)$ visualized below with
the following notations
\begin{enumerate}
\item $M:=\varphi_z^*(L)$
\item $Y= supp({\cal H}^{-a}(K))$ is of dimension $a$
\item $Y= Supp({\cal H}^{-a}(M))$
\item ${\cal Z}=H^{2a}(Y,{\cal H}^{-a}(K) \otimes {\cal H}^{-a}(M))$
\item ${\cal Y} = H^{2a-1}(Y,{\cal H}^{-a}(K) \otimes {\cal H}^{-a}(M))$
\item ${\cal X}_r=H^{2r-2}(Y,{\cal H}^{-r}(K) \otimes {\cal H}^{-r}(M))$
\end{enumerate}

%0 & 0 & 0 & ..... & 0 & 0 & 0\cr X_1 & 0 & 0 & ..... & 0 & 0 &
%0\cr
%* & 0 & 0 & ..... & 0 & 0 & \cr
%* & * & X_2 & ..... & 0 & 0 & 0\cr
%* & * & * &  .... & 0 & 0 & 0\cr
%* & * & * &  .... & 0 & 0 & 0\cr
%* & * & * &  .... & X & 0 & 0\cr

$$ \xymatrix{ 0 & 0 & 0 & 0 &..... & 0 & 0 & 0 \cr
0 & 0 & 0 & 0 &..... & 0 & 0 & 0\cr {\cal X}_1 & 0 & 0 & 0 &..... & 0 & 0 & 0\cr \bullet & 0 &
0 & 0 &..... & 0 & 0 & 0\cr \bullet & \bullet & {\cal X}_2 & 0 & ..... & 0 & 0 & 0\cr \bullet &
\bullet & \bullet & 0 & .... & 0 & 0 & 0\cr
 \bullet & \bullet & \bullet & \bullet &  ..\  0 & 0 & 0 & 0\cr \bullet & \bullet & \bullet & \bullet &
.... & {\cal X}_{a} & {\cal Y} & {\cal Z}\cr} $$ The diagonal is isomorphic to its
$E_\infty$-limit ${\cal H}^0(K*L)_z$. The diagonal to the left of the diagonal describes the
$E_2$-term. The $E_\infty$-limit ${\cal H}^{-1}(K*L)_z$ is obtained by taking successive
kokernels with respect to the higher differentials
$$ d: {\cal X}_{a-1} \to {\cal Y} \quad , \quad \delta: {\cal X}_{a-2}\to {\cal Y}
 \quad \cdots \quad \delta: {\cal X}_{1}\to {\cal Y} \ .$$
Since all ${\cal X}_r$ have weights $\leq -2$, we get an exact sequence
$$    {\cal Y}_{w\leq -2} \to {\cal H}^{-1}(K*L)_z \to
Gr({\cal Y})_{w=-1} \to 0 \ $$ from the weight filtration of ${\cal Y}$.

\bigskip\noindent
 \underbar{Purity assumption}: Suppose $H^{2a-1}(Y,{\cal
H}^{-a}(K) \otimes {\cal H}^{-a}(M))$ is pure of weight $-1$.

\bigskip\noindent
Under this purity assumption  all higher differentials $d:{\cal  X}_r\to {\cal Y}$ are zero,
since all  ${\cal X}_r$ have weight $\leq -2$. Hence the ${\cal  X}_1$-term survives in the
spectral sequence
$$ {\cal  X}_1= H^0\Bigl(Y,{\cal H}^{-1}(K)\otimes {\cal H}^{-1}(L)\Bigr)
\hookrightarrow {\cal H}^{-2}(K*L)_z \ .
$$ Furthermore under this purity assumption  $$\fbox{$ {\cal H}^{-1}(K*L)_z \cong H^{2a-1}\Bigl(Y,{\cal
H}^{-a}(K) \otimes {\cal H}^{-a}(M)\Bigr) $}\ .$$ Theorem \ref{thm1} below is a special case of
this assertion.

\bigskip\noindent
\underbar{The subcase $dim\ supp(L) > dim\ supp(K)=a$}: Then only
the term for $b=a+1$ contributes to the stalk cohomology and
$$ {\cal H}^{-1}(K*L)_z\neq 0 \Longrightarrow supp({\cal
H}^{-a}(K)) \subseteq z - supp({\cal H}^{-a-1}(L))  \ .$$ Again there is a single $E_2$-term in
the spectral sequence, which contributes to the $E_\infty$-limit term ${\cal H}^{-1}(K*L)_z$
coming from $a=dim\ supp(K)$ and $b=a+1\leq dim\ supp(L)$. ${\cal H}^{-1}(K*L)_z$ is a
subquotient of the $E_2$-term
$$ H^{2a}\Bigl(supp(K),{\cal H}^{-a}(K)\otimes {\cal
H}^{-a-1}(\varphi_z^*(L))\Bigr)\ .$$ Notice $dim\ supp ({\cal H}^{-a-1}(\varphi_z^*(L))) \leq
a$. Hence $K$ is an irreducible component of $supp({\cal H}^{-a-1}(\varphi_z^*(L)))$ of highest
dimension  unless $$dim\ (supp(K)) +1=a+1= dim\ (supp(L))\ .$$

\bigskip\noindent
This discussion implies

% If $L=\delta_E$, then the support of
%${\cal H}^{-a-1}(\varphi_z^*(L))$ has dimension $<a$ unless

\bigskip\noindent
\begin{Lemma} {\it Suppose $K$ and $L$ are irreducible perverse sheaves on $X$ with $dim\
(supp(K)) \leq dim\ (supp(L))$.  Then
$$ \fbox{$ supp\ \Bigl({\cal H}^{-1}(K*L)\Bigr)\ \subseteq\ \{ z\in X \ \vert
\ supp(K) \subseteq z - supp(L)\} $} \ .$$ If both supports have the same dimension, then
$$ supp\ \Bigl({\cal H}^{-1}(K*L)\Bigr)\ \subseteq\ S(K,L)\ .$$} \end{Lemma}

\bigskip\noindent
\begin{Theorem}\label{thm1} {\it Suppose $K=\delta_Y$ and $L=\delta_Z$ for irreducible subvarieties $Y$ and $Z$ of $X$
of the same dimension $d$. Then ${\cal H}^{-1}(K*L)_z$ vanishes unless $z\in S(K,L)$. Suppose
the singularities of $Y$  have codimension $\geq 2$ in $Y$. Then for $z\in S(K,L)$
$$ \fbox{$ {\cal H}^{-1}(K*L)_z \cong H^{2d-1}(Y) \cong IH^{2d-1}(Y) $} \ .$$} \end{Theorem}

\bigskip\noindent
\underbar{Proof}: If $z\in S(K,L)$, then $supp(K)=supp(\varphi_z^*(L))=Y$.  We may assume for
simplicity $z=0$ and $Y=Z$ for simplicity of notation. Then ${\cal H}^{-a}(M)=\overline\Q_l$
and ${\cal H}^{-a}(K)=\overline\Q_l$ for $a=d=dim(Y)$ outside a subvariety $S$ of codimension
two in $Y$. Hence $H^{2a-1}(Y,{\cal H}^{-a}(K)\otimes {\cal H}^{-a}(M))$ is isomorphic to $$
H_c^{2d-1}(Y\setminus S, \overline\Q_l) \cong H^{2d-1}(Y)\ .$$ This follows from the long exact
sequence for  cohomology with compact supports attached to $(Y\setminus S,Y,S)$. To proof the
theorem following thew discussion above it only remains to prove the purity assertion for
$H^{2d-1}(Y)$, which will follow from the statement $$H^{2d-1}(Y)= IH^{2d-1}(Y)\ .$$ This last
assertion follows from the long exact sequence attached to the distinguished triangle
$\psi_Y\to \lambda_Y\to \delta_Y\to \psi_Y[1]$
$$ H^{d-1}(Y,\psi_Y) \to H^{d-1}(\lambda_Y) \to IH^{2d-1}(Y) \to
H^{d}(Y,\psi_Y)\ . $$ Since $\psi_Y\in {}^pD^{\leq 0}(X)$  and $dim\ supp({\cal
H}^{-\nu}(\psi_Y)) \leq max(\nu,d-2)$ the spectral sequence $ H^\mu(X,{\cal
H}^{-\nu}(\psi_Y))\Longrightarrow H^{\mu-\nu}(X,\psi_Y)$ has zero $E_2$-terms for the degrees
$\mu-\nu=d-1$ and $\mu-\nu=d$. Hence $H^{d-1}(Y,\psi_Y)=0$ and $ H^{d}(Y,\psi_Y)=0$. Since
$\lambda_Y\vert_{Y\setminus S}=\overline\Q_{l,Y\setminus S}$  outside the singular locus $S$,
the isomorphism $H^{d-1}(\lambda_Y)\cong H^{2d-1}(Y)$ is again obtained from the long exact
sequence for cohomology with compact supports.

%\bigskip\noindent
%{\bf Remark}: Since $Ext^1_{Perv(X)}(DK,\varphi_z^*(L))=0$ for
%irreducible perverse Weil sheaves of equal weight we may conclude,
%that the subspace of Frobenius invariant elements in ${\cal
%H}^{-1}(K*L)_z$ is trivial for all $z\in X$.

\section{The Theta divisor}
In this section we assume $char(k)=0$. Let $\Theta$ be a polarization of $X$ such that
$(X,\Theta)$ is irreducible as a polarized abelian variety. This is equivalent to the
condition, that the divisor $\Theta\subseteq X$ is an irreducible subvariety of $X$. For
basefields $k$ of characteristic zero then the theta divisor is a normal variety.  Using that
the codimension of the singular locus is $\geq 2$ (see \cite{EF}) this follows from the
Krull-Serre criterion, since $\Theta$ is the zero locus of a section of the line bundle defined
by the polarization, i.e defined by one equation in the smooth variety $X$. The perverse sheaf
$$K=\delta_\Theta \ $$ is a self dual $K=DK$. Since $\Theta$ is a
polarization, the map $X\to Pic(X)$ defined by $x\mapsto cl(\Theta +x)-cl(\Theta)$ is an
isogeny,  $Aut(\delta_\Theta)$ is finite. In particular $S(K,L)$ is finite for any perverse
sheaf $L$. $\Theta$ is a principal polarization if and only if the cardinality of $Aut(K)$ is
one. A principal polarization is a translate of a symmetric principal polarization. For a
symmetric principal polarization $\Theta=-\Theta$ one has $Aut(K)=\{0\}$. For a principal
polarization $\Theta$ there exists a unique $\kappa\in X$ such that $\varphi_\kappa^*(K)\cong
K$ for $K=\delta_\Theta$. Hence
$$ S(K,K) = \{\kappa\} \ .$$

\bigskip\noindent
\begin{Corollary} {\it If $K=\delta_\Theta$ for an irreducible principal polarization $\Theta$ of $X$, then
${\cal H}^0(K*L)_z$ is zero unless $L\cong \delta_{z-\kappa+\Theta}$. In particular ${\cal
H}^0(K*K)_z=0$ unless $z=\kappa$. Hence $$ \fbox{$ {\cal H}^0(\delta_\Theta*\delta_\Theta) =
\delta_{\{\kappa\}} $} \ .$$} \end{Corollary}

\bigskip\noindent
%Notice $Y=-Y$ for a symmetric theta divisor $Y=\Theta$. In this
%case $\varphi_0^*(\delta_Y)=\delta_Y$.  Once more using that the
%support of $\psi_\Theta$ has dimension $\leq (g-1)-2=g-3$ we see,
%that
%$$IH^{2g-3}(\Theta)\cong H^{2g-3}(\Theta)\ .$$ Here we used, that
%the cohomology $H^i(X,K)$ of any semi-perverse sheaf complex $K$
%in ${}^pD^{\leq}(X)$ vanishes, if $i> dim\ supp(K)$. So the
%intersection cohomology and the ordinary cohomology coincide in
%degree $2g-3$.

%
%\bigskip\noindent
%{\bf Corollary}: {\it For an irreducible perverse sheaf
%$K=\delta_Y$ attached to an irreducible subvariety $Y\subseteq X$
%and for $z\in S(K,K)$ we get (ignoring Tate twists) $$ {\cal
%H}^{-1}(\delta_Y
%* \varphi_0^*(\delta_Y))_0 \cong H^{dim(Y)-1}(Y,\delta_Y) \cong IH_1(Y) \ .$$}

\bigskip\noindent
If $X$ is the Jacobian $X$ of a curve $C$, then it will be shown in \ref{Cohomology groups}
corollary \ref{cor8}
$$IH_1(\Theta) \cong H_1(X)\ .$$ More generally
$IH_1(\Theta)\cong H_1(X)$ holds for arbitrary principal polarized abelian variety, provided
the singularities of $\Theta$ have codimension $\geq 3$ in $\Theta$ (see \cite{W2}). Hence

\bigskip\noindent
\begin{Corollary}\label{Thet} {\it For principally polarized abelian variety $(X,\Theta)$, where $(X,\Theta)$ is either a
Jacobian or where the singularities of $\Theta$ are of codimension $\geq 3$ in $\Theta$, then
the theta divisor $\Theta$ is irreducible and
$$ \fbox{$ {\cal H}^{-1}(\delta_\Theta*\delta_\Theta) = H^1(X) \otimes
\delta_{\{\kappa\}} $} \ .$$} \end{Corollary}

\chapter{Convolution with curves}\label{pingpong} Let $C\hookrightarrow X$ be an irreducible
smooth projective curve of genus $g\geq 2$ in $X$. Suppose $E$ is a smooth etale $\overline
\Q_l$-sheaf on an open dense subset of $C$. Let $\delta_E\in Perv(X)$ be the corresponding
perverse sheaf with support in $C$. Then $E$ extends to a (not necessarily smooth)
$\overline\Q_l$-sheaf on $C$ also denoted $E$, such that $\delta_E=E[1]$. Let $K$ be an
irreducible perverse sheaf on $X$ with support $Y$ of dimension $d_K$. The convolution $$ L=
\delta_E * K  $$ is the direct image complex of $\delta_E\boxtimes K$ under the the smooth
proper morphism
$$ f: C \times X \longrightarrow X \ $$
defined by the restriction of $a$, i.e by $f(x,y)=x+y$. Since $f$ is smooth with fibers of
dimension $1$ by [KW] III.7.1
$$ L= \delta_E * K \quad \in \ {}^p D^{[-1,1]}(X) \ .$$
Our aim is to show

\begin{Theorem}\label{thm2}
Suppose $C$ generates $X$ as an abelian variety, and suppose the dimension $d=d_K$ of the
support of $K$ is $< g=dim(X)$, then $L=\delta_E*K$ is a perverse sheaf.
\end{Theorem}

\bigskip\noindent
\begin{Corollary}\label{pin} {\it For $d\leq dim(X)$ iterated convolutions $\delta_{E_1}*\cdots *\delta_{E_d}$ of perverse
sheaves $\delta_{E_i}$ with support in smooth projective curves $C_i$ contained in $X$ are
perverse sheaves, if $d-1$ of the curves are generating curves for $X$}.
\end{Corollary}

\bigskip\noindent
A refined version of theorem will be stated below in theorem \ref{refined}.

\bigskip\noindent
\underbar{Spectral sequence}: The stalk cohomology sheaves ${\cal
H}^{-\nu}(L)_x$ at a point $x\in X$ can be computed from a
spectral sequence, where the following $E_2$-terms contribute
$$ \xymatrix@C=-1.2cm{H^0(x-C,\varphi_x^*(E)\otimes {\cal H}^{-\nu+1}(K)) \ar[drr]_{d_2} & \bullet   & \bullet  \cr
\ar[drr]_{d_2} \bullet & H^1(x-C,\varphi_x^*(E)\otimes {\cal H}^{-\nu}(K)) & \bullet \cr
\bullet & \bullet & H^2(x-C,\varphi_x^*(E)\otimes {\cal H}^{-\nu-1}(K)) \cr }
$$
All middle terms $H^1(x-C,\varphi_x^*(E)\otimes {\cal H}^{-\nu}(K))$ survive in the
$E_\infty$-term. The other contributions in $E_\infty$ are the kernel of the differential
$$d_2: H^0(x-C,\varphi_x^*(E)\otimes {\cal H}^{-\nu+1}(K))
\longrightarrow H^2(x-C,\varphi_x^*(E)\otimes {\cal H}^{-\nu}(K))$$ respectively the kokernel
of the differential
$$d_2: H^0(x-C,\varphi_x^*(E)\otimes {\cal H}^{-\nu}(K))
\longrightarrow H^2(x-C,\varphi_x^*(E)\otimes {\cal
H}^{-\nu-1}(K))\ .$$

\bigskip\noindent
\underbar{Notation}: $S\ominus T=\bigcap_{t\in T} (S-t)\subseteq S$. Notice $(S\ominus T)\oplus
T \subseteq S$.

\bigskip\noindent
\underbar{Claim}: If $T$ is irreducible and generates the abelian variety $X$, then
$$ dim (S\ominus T) < dim(S) \quad \mbox{ or } \quad S=X \ .$$
Similarly $S\oplus T \supseteq S$ and
$$ dim (S\oplus T) > dim(S) \quad \mbox{ or } \quad S=X \ .$$

\bigskip\noindent
\underbar{Proof}: Suppose $dim(S\ominus T)=dim(S)$. Then a component $S_1$ of $S$ of highest
dimension is contained in $S\ominus T$. Suppose $S_1,..S_r$ are the irreducible components of
$S$ of highest dimension, which are contained in $S\ominus T$. Let $\eta_i$ is a generic point
of $S_i$ and $t$ a generic point of $T$. Then $\eta_i - t \in S$ is one of the generic points
$\eta_j$ of $S$. Hence the $\eta_i$ are permuted, and there exists an integer $n_i$ such that
$\eta_i - n_i \cdot t = \eta_i$. Now, since $t$ generates $X$ if and only if $n_i\cdot t$
generates $X$ (consider the map $n_i:X\to X$), we get $S_i - x = S_i$ for all $x\in X$. Hence
$S_i=X$.

\bigskip\noindent
\underbar{Supports}: Let $S^\nu(K)=S^\nu$ and $S^\nu(L)$ denote the supports (see the
convention \pageref{supports}) of the stalk cohomology sheaves ${\cal H}^{-\nu}(K)$ and ${\cal
H}^{-\nu}(L)$ respectively. Then

\begin{Lemma}\label{lemma8} For $L=\delta_E*K$ we have
$$ S^\nu(L)\ \subseteq\ \bigl( S^{\nu-1}(K) \oplus C\bigr) \cup \bigl( S^{\nu}(K) \ominus
-C\bigr) \cup \bigl( S^{\nu+1}(K) \ominus -C\bigr)\ .$$ \end{Lemma}

\bigskip\noindent
\underbar{Proof}: This is almost obvious. The cohomology groups $H^i(x-C,...)$ vanish for
$i=1,2$ unless $x-C \subseteq S^\nu(K)$ respectively $x-C \subseteq S^{\nu+1}(K)$. Moreover, in
the second case the coefficient system $\varphi_x^*(E)\otimes {\cal H}^{-\nu-1}(K)$ must admit
a nontrivial constant quotient sheaf in order to contribute.

\bigskip\noindent
This being said, we make the

\bigskip\noindent
\underbar{Assumption}: $C$ generates $X$.

\bigskip\noindent

\bigskip\noindent
For the proof of theorem \ref{thm2} we can assume $K$ to be irreducible with dimension $d_K>0$.
If $d_K=0$, then $K=\delta_{x_0}$ for some closed point $x_0$ of $X$. Hence $\delta_E*K$ is a
translate of $\delta_E$. Hence we will assume $$1\leq d=d_K < g\ $$ for the rest of this
section. Then, since by assumption $K$ is irreducible perverse, we know $dim\ S^\nu(K) < \nu$
for $\nu < d$ and $dim\ S^{d}(K) = d$ for the dimensions of the support. Hence lemma
\ref{lemma8} implies

\begin{enumerate}
\item $S^\nu(L)=\emptyset  $ for $\nu <0$.
\item $dim\ S^\nu(L) < \nu $ for alle $0\leq \nu < d -1 $
\item $S^{d-1}(L) \subseteq S^{d}(K)\ominus -C$ (hence there exist
components of dimension $g-1$ only if $S^{d}(K)\ominus -C$ has
dimension $d-1$; otherwise $dim\ S^{d-1}(L) < d-1$)
\item $dim\ S^{d}(L) < d$
\item $S^{d+1}(L) \subseteq S^{d}(K)\oplus C$ has dimension $\leq
d+1$.
\end{enumerate}
For 4. we used $d_K< g$ to get the inequality  $dim\ (S^d(K)\ominus -C) < d$. The listed
properties show
$$ \delta_E * K \in {}^p D^{\leq 0}(X) \ .$$
The same applies for $DK$ and $D\delta_E$, hence $D(\delta_E*K) \in {}^p D^{\leq 0}(X)$. This
implies theorem \ref{thm2}.

\bigskip\noindent
We slightly improve on theorem \ref{thm2}:  In fact [KW] III.5.13 implies
%, that either
%$dim\ S^{d+1}(L) = d+1$ or $S^{d+1}(L) =\emptyset$. Hence
$ S_{d+1}(L)=\emptyset$  unless $ dim\ S^d(K)\oplus C = d+1$. We already used above, that
$H^0(x-C,E\otimes {\cal H}^{-d}(K))$ is nonzero only for $x\in C\oplus S^d$. If $x-C \nsubseteq
S^d$ or equivalently if $x\notin S^d\ominus -C$, then $x-C$ and $S^d$ intersect in finitely
many points $x_i$. If one of the points $x_i$ is contained in the open subset $U\subseteq
S^d(K)$, where ${\cal H}^{-d}(K)$ is smooth, then $H^0(x-C,\varphi_x^*(E)\otimes {\cal
H}^{-d}(K))$ is nonzero, if $\varphi_x^*(E)_{x_i}\neq 0$, since then $\varphi_x^*(E)\otimes
{\cal H}^{-d}(K)$ is a skyscraper sheaf with nonzero stalk at $x_i$. This is the case for
$$ x \in (U_E \oplus U) \quad , \quad x\not\in (S^d(K)\ominus -C) \ .$$
since then $\varphi_x(u') = x-u' = x_i \in U$ for $u'\in U_E$,
where $U_E\subseteq C$ is some dense open subset, where $E$ is
smooth. Such points $x$ exists, if the dimension of $U_E \oplus U$
or $C\oplus S^d(K)$ is $d+1$, since the dimension of
$S^d(K)\ominus -C$ is $<d$. So, to show the converse it is enough
to observe, that no further cancelation arises from the
differentials of the spectral sequence. However, this is obvious
in present case, since the only relevant differential is the
$d_2$-differential
$$ d_2: H^0(x-C,\varphi_x^*(E)\otimes {\cal H}^{-d}(K))
\longrightarrow H^2(x-C,\varphi_x^*(E)\otimes {\cal H}^{-d
-2}(K))=0 \ .$$ Furthermore for a  point $(u',y)\in U_E\times U$
in general position and $x=u'+y$, all other solutions $u''+y'=x$,
where $u''\in C$ and $y'\in S^d(K)$, are points where $u''$ is in
general position in $C$ and $y'$ is in general position in
$S^d(K)$. Hence we conclude for such $x$, that $ dim\
H^0(x-C,\varphi_x^*(E)\otimes {\cal H}^{-d}(K))$ is equal to $$
rank(E)\cdot rank({\cal H}^{-d}(K)) \cdot deg(a:C\times Y \to
Y\oplus C)\ .
$$

\bigskip\noindent
\underbar{Notation}: For curves $M$ we use the notation $S\circleddash M$, where $S\circleddash
M = S\ominus M$, if the dimension $of S\ominus M$ is $dim(S)-1$, and where $S
\circleddash M =\emptyset$ otherwise. %Then $(S\oplus M)\circleddash M =S$ for $dim(S) \leq
%g-2$, and $(S\circleddash M)\oplus M =S$ or $=\emptyset$.

%Notice
%$$ \xymatrix{ & S\oplus M \ar[dr]^{  \circleddash M}  & \cr S \ar[dr]_{  \circleddash M}\ar[ur]^{\oplus M} & & S \cr
% & S\circleddash M \ar[ur]_{\oplus M} & \cr} $$
%commutes for $dim(S) \leq g-2$.

\bigskip\noindent
\begin{Theorem} \label{refined}{\it Suppose the curve $C$ generates $X$ and $\delta_E$ is a perverse sheaf with support $C$.
Then, for an irreducible perverse sheaf of geometric origin $$K = K_d \in Perv(Y) \subseteq
Perv(X)$$ with support $Y$ of dimension $d < dim(X)$, the convolution $L=\delta_E
* K$ with  $\delta_E$ is a perverse sheaf on $X$
$$ L = L_{d-1} \oplus L_{d+1} \quad \in \ Perv(X)\
$$ with two summands
\begin{enumerate} \item $L_{d+1}$ is a
nonzero perverse sheaf, whose irreducible perverse components are supported in the
$d+1$-dimensional irreducible scheme $\ Y\oplus C$. The generic rank $rank ({\cal
H}^{-d-1}(L))$ of the corresponding coefficient system defined by $L$ on $\ Y\oplus C$ is
$$rank ({\cal H}^{-d-1}(L)) \ = \ rank(E)\cdot rank({\cal
H}^{-d}(K)) \cdot deg(a:C\times Y \to Y\oplus C)\ .
$$

\item $L_{d-1}$ is a perverse sheaf, whose irreducible
perverse components are supported in the $(d-1)$-dimensional
scheme $\ Y\circleddash -C$.
\end{enumerate}} \end{Theorem}

\bigskip\noindent
Notice, that even if $Y\circleddash -C$ is nonempty, the
\lq{exceptional}\rq\ perverse sheaf $L_{d-1}$ from the theorem
above might be trivial. For $L_{d-1}$ to be nontrivial it is
required, that $ E\otimes {\cal H}^{-d}(K)\big\vert_{Y\circleddash
-C} $ has a nontrivial constant etale quotient sheaf in a
neighborhood of a generic point $\eta$ of $Y\circleddash -C$. More
precisely, the kokernel
$$ H^0(\eta-C,E\otimes {\cal H}^{1-d}(K)\big\vert_{Y\circleddash -C})
\longrightarrow H^2(\eta-C,E\otimes {\cal H}^{-d}(K)\big\vert_{Y\circleddash -C})
$$
must not vanish.

\bigskip\noindent
\underbar{Example}: For two smooth projective curves $C$ and $C'$ with irreducible perverse
sheaves $\delta_E$ and $\delta_{E'}$ of support $C$ and $C'$ the convolution $L=\delta_E
*\delta_{E'}$ is a direct sum $L=L_0\oplus L_2$. The perverse sheaf $L_0$ vanishes except in
the case, where $T_x^*(\delta_{E'}) \cong (\delta_E)^\vee$ holds for some $x$ in $X$. Similarly
${\cal H}^{0}(\delta_C *\delta_{C'})$ and ${\cal H}^{-1}(\delta_C *\delta_{C'})$ vanish except
if $C=x-C'$ holds for some point $x\in X$. If $X$ is the Jacobian of $C$ this point $x$ is
uniquely determined, if it exists. Now suppose $C'=x- C$. By replacing $C'$ by a translate now
let us assume without restriction of generality, that $C'=-C$. Then
$$ \delta_C * \delta_{C'} =\delta_C *\delta_{-C} = \delta_0 \oplus L_2 $$
where $L_2$ is a perverse sheaf with support in $C-C$. Furthermore
$$ {\cal H}^{-1}(L_2) \cong H^1(C,\overline \Q_l)
\otimes \delta_0 \ .$$ The rank of the smooth coefficient system defining $L_2$ is the degree
of the map $$C\times C'\to C\oplus C' \subseteq X\ .$$ See \ref{CundA} for further details.

%\bigskip\noindent
%\underbar{Example}: If $C$ is a hyperelliptic curve and $X$ is its
%Jacobian, then there exists a unique point $e\in X$ such that
%$$ C = e - C \ .$$
%Hence $L=\delta_C*\delta_{-C}$  is  $L=L_0 \oplus L_2$, where
%$L_0$ is the Dirac sheaf $L_0= \delta_{0}$ and $L_2$ has generic
%rank two. Furthermore, since $\delta_C *\delta_C$ is a translate
%of $L$, it contains the Dirac sheaf $L_0=\delta_{0}$ translated by
%$e$ and similar the translate of $L_2$.
%
%\bigskip\noindent
%Before we discuss the remaining case, where the support of $K$ has
%dimension $d=dim(X)$ we make the following

\bigskip\noindent
\underbar{Remark}: If $char(K)=0$ part of the discussion above is a special case of the fact,
that for an irreducible perverse sheaf $K=\delta_E$ with support $Y_K$ of dimension $d=d_K$ one
can consider the class
$$ c(K)=rank({\cal H}^{-d}(K)) \cdot [Y_K] \ \in \ H_d(X,\Z) \ ,$$
where $rank({\cal H}^{-d}(K))$ is the rank of the coefficient system ${\cal H}^{-d}(K)$ at the
generic point of $Y_K$. For another $L=\delta_F$ with support $Y_L$ in $X$, the convolution
$K*L$ need not be perverse. Assume $d_K+d_L \leq g$. Then ${\cal H}^{-\nu}(K*L)=0$ for $\nu<
d_K + d_L$. If we define $c(K*L)\in H_{d_K+d_L}(X,\Z)$ to be $rank({\cal H}^{-d_K-d_L}(K*L))
\cdot [Y_K\oplus Y_L]$, where $[Y_K\oplus Y_L]$ is understood to be zero, if the dimension of
$Y_K\oplus Y_L$ is not equal to $d_K+d_L$, then
$$ c(K*L) = c(K)*c(L) \ ,$$
in the homology ring $H_*(X,\Z)$ endowed with the $*$-product,
which is induced on homology by the map $a:X\times X\to X$ via the
K\"{u}nneth theorem.

\bigskip\noindent
\section{The highest dimensional case}
We continue the discussion of convolutions and study $L=\delta_E*K$, where we now consider the
remaining case where $K$ is an irreducible perverse sheaf with support $X$. We retain the
notations and assumptions of the last section except, that now $d=d_K=dim(X)$.  Although the
discussion is similar to the case where $d<dim(X)$, we have to consider the cohomology degrees
$\nu=g-1,g$ and $\nu=g+1$ more carefully.
 To simplify the discussion we do not deal with
the most general case, but assume

\bigskip\noindent\underbar{Assumptions}: Let $X$ be generated by
$C$. Assume one of the following equivalent conditions holds:
\begin{enumerate}
\item $L=\delta_E*K$ is perverse. \item
 $L=\delta_E*K\in {}^p D^{\leq 0}(X)$.
 \item $S_{g+1}(L)=\emptyset$
\item $H^0\bigl(x-C,\varphi_x^*(E)\otimes {\cal H}^g(K)\bigr)= 0$
for all points $x\in X$. \end{enumerate} Notice [KW] III.5.13 for
3. Furthermore notice, that the spectral sequence collapses in
degree $\nu=g+1$, again by [KW] III.5.13.

\bigskip\noindent
\underbar{Remark}: Concerning this assumption recall that $\delta_E*K\in {}^p D^{[-1,1]}(X)$.
The perverse sheaves $f^*{}[1] {}^p H^{\pm 1}(L)$ are the maximal perverse subsheaves
respectively the maximal perverse quotient sheaves of $\delta_E\boxtimes K$ on $X$. Therefore,
if ${}^p H^{\pm 1}(L)$ is a constant sheaf on $X$ (a situation which will be typical later, see
also the curve lemma \ref{Enlarge} and lemma \ref{notra}) and $E$ is irreducible, then
$\delta_E\boxtimes K$ is an irreducible perverse sheaf on $C\times X$, and if ${}^p H^{\pm
1}(L)$ were nonzero, then $\delta_E\boxtimes K\cong f^*[1]({}^p H^1(\delta_E*K))$ must be the
trivial perverse sheaf on $C\times X$. This forces $K=\delta_X$ and $E=\overline \Q_l$ under
the assumptions made.

\bigskip\noindent
This being said, we now discuss the cohomology ${\cal H}^{-\nu}(\delta_E*K)$ under the
assumption above. We now skip the cases where $\nu<g-1$, since in these cases the discussion
remains the same as in the last section.

\bigskip\noindent
\underbar{The degree $\nu=g-1$}:  This case is similar to the previous case $\nu<d_K-1$. The
support $S_{g-1}(L)$ is contained in the union of $ S^{g-2} \oplus C$ and $ S^{g-1} \ominus -C$
and $S^{g} \ominus -C$, except that now $S^{g} \ominus -C=X\ominus -C=X$. Since the dimension
of $ S^{g-2} \oplus C$ and $ S^{g-1} \ominus -C$ is $< g-1$, we get an exceptional perverse
constituent $L_{d-1}$ in $L=\delta_E*K$ only if
$$ Z=\Bigl\{ x\in X\ \ \Big\vert H^2\bigl(x-C,\varphi_x^*(E)\otimes {\cal
H}^{-g}(K)\bigr) \Bigr\} $$ is nonempty and has dimension $g-1$.
The dimension of $Z$ is at most $g-1$ by our assumption $L$ to be
perverse. Notice, this does not suffice for $L_{g-1}$ to be
nonzero. In addition the differential $d_2$ from the spectral
sequence, which computes ${\cal H}^{1-g}(L)_x$
$$ d_2: H^0\bigl(x-C,\varphi_x^*(E)\otimes {\cal H}^{1-g}(K)\bigr)
\longrightarrow H^2\bigl(x-C,\varphi_x^*(E)\otimes {\cal
H}^{-g}(K)\bigr) \ ,
$$  must have nontrivial cokernel  for $x$ in a
$g-1$-dimensional open subset of $Z$. This now is a necessary and
sufficient condition for $L_{g-1}$ to be nontrivial. Notice, the
left side vanishes for $x\notin C\oplus S^{g-1}(K)$. Cancelations
might be possible on irreducible components, which are common
irreducible components both of $C\oplus S^{g-1}(K)$ and $Z$ of
dimension $g-1$.

\bigskip\noindent
\begin{Theorem} \label{thm4}{\it Suppose $X$ is generated by the curve $C$. Suppose $\delta_E$ has support $C$, suppose
$L=\delta_E*K$ is a perverse sheaf on $X$ and $K$ is an irreducible perverse sheaf of geometric
origin with support of dimension $g$. Then there exists a decomposition
$$ L = L_{g} \oplus L_{g-1} \ ,$$
where $L_{g-1}$ and $L_g$ are perverse sheaves, whose irreducible constituents  have supports
of dimension $g-1$ and $g$ respectively. Furthermore the generic rank of the coefficient system
defining $L_g$ is $$ dim\ H^1(x-C,\varphi^*_x(E)\otimes {\cal H}^{g}(K)) \ ,$$  where $x$ is a
generic point of $X$.} \end{Theorem}

\bigskip\noindent
By the assumptions of the theorem $H^0(x-C,\varphi^*_x(E)\otimes {\cal H}^{g}(K))$ vanishes for
all $x\in X$, as already explained. The algebraic set $Z$ of $x\in X$, where
$H^2(x-C,\varphi^*_x(E)\otimes {\cal H}^{g}(K))$ does not vanish, has dimension $\leq g-1$ by
the same reason. Hence for a generic point the the generic rank of $L_g$ can be computed by the
Euler characteristic. We state the formula for the Euler characteristic in the tame case, which
then gives
$$ rank(L_g) \ = \ - \chi^{tame}(x-C,\varphi^*_x(E)\otimes {\cal
H}^{g}(K)) $$
$$ = - rank(K) \cdot (2-2g)   + \sum_{c\in \Sigma} \Bigl( rank((\varphi^*_x(E)\otimes {\cal
H}^{g}(K)) \ - \ rank ((\varphi^*_x(E)\otimes {\cal H}^{g}(K))_c
\Bigr) \ .$$ Here $\Sigma$ is a finite set of closed point of $C$,
where the underlying etale sheaf $(\varphi^*_x(E)\otimes {\cal
H}^{g}(K)$ is smooth on $x-C$. Therefore
$$ rank(L_g) = 2(g-1)rank(K) + \sum_{c\in \Sigma}
\Bigl( rank(K)rank(E) - rank (E\otimes \varphi^*_x{\cal
H}^{g}(K))_c\Bigr) \ .$$

\section{The curve lemma} \label{Enlarge}
For an abelian variety $X$ let $\tilde C\hookrightarrow X$ be a irreducible projective curve
and let $\pi: C\to \tilde C$ be its normalization. We want to show, that convolution of a
perverse sheaf $K$ (of geometric origin) with $\delta_{\tilde C}$ again is a perverse sheaf.
Convolution with skyscraper sheaves preserves the category $Perv(X)$ for trivial reason.
Therefore, since $\pi_*(\delta_C)$ is the direct sum of $\delta_{\tilde C}$ and a skyscraper
sheaf, we may replace $\delta_{\tilde C}$ by the sheaf $\pi_*(\overline \Q_{l,C}[1])$ in the
following without loss of generality. The  morphism $f$ defined by $f(x,y)=x+\pi(y)$
$$ \xymatrix@+0,5cm{ X \times C \ar[r]^-{pr_1}\ar[d]^\sim_F &  X \cr X \times C \ar[ru]_f\ar[r]_{id\times\pi} & \ \ X\times \tilde C \ar[u]_a\cr } $$
 is smooth equidimensional with geometrically connected fibers of
dimension 1, since $F(x,y)=(x-\pi(y),y)$ is an isomorphism. Notice $f\circ \nu=id_X$ for the
inclusion $\nu:X\to X\times C$ defined by $x\mapsto (x,0)$, if $0\in C(k)$. Since $K\boxtimes
\overline\Q_{l,C}[1]$ is a perverse sheaf on $X\times C$ for $K\in Perv(X)$, therefore
$$ L= K * \pi_*(\delta_{C}) \ =\ Rf_*(K\boxtimes \overline\Q_{l,C})[1]
\ $$ satisfies (\cite{KW}, III.7.1)
$$ L \in {}^p D^{[-1,1]}(X,\overline\Q_l) \ .$$

\bigskip\noindent
Assume the base field is the algebraic closure $k$ of a finite field $\kappa$ with $q$
elements. For an algebraic scheme $Y_0$ over $\kappa$ let denote $Y=Y_0 \times_{Spec(\kappa)}
Spec(k)$. Assume that the embedding $\tilde C\hookrightarrow X$ comes from a closed immersion
$\tilde C_0 \hookrightarrow X_0$ defined over $\kappa$. Let $K_0$ be a perverse sheaf on $X_0$,
with extension $K$ on $X$. Let $\kappa_m$ be the unique extension field of $\kappa$ of degree
$m$ in $k$. Let $F$ respectively $F_m$ denote the geometric Frobenius of $Gal(k/\kappa_m)$. The
set of rational points $X_0(\kappa_m)$ can be identified with the set of fixed points
$X_0(k)^{F_m}$. For a closed point $x$ choose a geometric point $\overline x$ over $x$. For a
complex $K_0\in D_c^b(X_0,\overline\Q_l)$ let
$$ f^{K_0}: X_0(\kappa) \to \overline\Q_l $$
be the function defined by
$$ f^{K_0}(x) = \sum_\nu (-1)^\nu Tr(F; {\cal H}^\nu(K)_{\overline
x}) \ .$$ Similarly define $f_m^{K_0}: X_0(\kappa_m) \to \overline\Q_l$ for all $m\geq 1$.
Recall the following well known facts (e.g. [L], \cite{KW} III.12.1)
\begin{enumerate}
\item (Cebotarev) For semisimple perverse sheaves $K_0$ and $L_0$
on $X_0$, the equality $f_m^{K_0}(x)=f_m^{L_0}(x)$ for all $m$ and all $x\in X_0(\kappa_m)$
implies  $K_0\cong L_0$ in $D_c^b(X_0,\overline\Q_l)$.
\item $f_m^{K_0\oplus L_0}(x)= f_m^{K_0}(x) \ +\  f_m^{L_0}(x) $
for all $x\in X_0(\kappa_m)$
\item $f_m^{K_0\otimes^L L_0}(x)= f_m^{K_0}(x)\cdot f_m^{L_0}(x) $
for all $x\in X_0(\kappa_m)$
\item For a morphism $g_0: Y_0 \to X_0$ defined over $\kappa$ we
have $$f_m^{g_0^*(K_0)}(x)= f_m^{K_0}(g_0(x))$$ for all $x\in
Y_0(\kappa_m)$.
\end{enumerate}
For an abelian variety $$\pi_0:X_0\to Spec(\kappa)$$ over $\kappa$ the $\kappa_m$-rational
points $X_0(\kappa_m)$ define a finite group with respect to the group law of the abelian
variety. Suppose for all $m$ there exist subgroups $$ \Gamma_m \subseteq X_0(\kappa_m) \ $$
with index $[X_0(\kappa_m):\Gamma_m]$ bounded by a constant $C$ independently from $m$. Suppose
$K_0\in Perv(X_0)$ is a semisimple perverse sheaf on $X_0$ such that
$$  f_m^{K_0}(x) = f_m^{K_0}(0)\quad , \quad  \forall m\geq 1\ \ \forall x\in \Gamma_m \
.$$ Under these assumption we claim

\bigskip\noindent
\begin{Lemma}\label{l9} {\it The sheaf $K$ is a translation-invariant perverse sheaf on $X$.}
\end{Lemma}

\bigskip\noindent
\underbar{Proof}: For $N = C!$ and $x\in X_0(\kappa_m)$ we have $N\cdot x \in \Gamma_m$.
Multiplication by $N$ is an isogeny defined over $\kappa$
$$ g_0: \ X_0 \overset{N}{\longrightarrow} X_0 \ .$$
Since $g_0$ is a finite flat morphism $L_0=g_0^*(K_0)\in Perv(X_0)$,
$$ f_m^{L_0}(x) = f_m^{K_0}(N\cdot x) = f_m^{K_0}(0) $$
holds by the assumption above, since $N\cdot x \in \Gamma_m$. Let $i_0: \{0\} \to X_0$ be the
inclusion of the neutral element and $M_0=i_0^*(K_0)$. By definition
$f_m^{M_0}(0)=f_m^{K_0}(0)$ and $M_0\in D_c^b(Spec(\kappa),\overline\Q_l)$. Let $A_0$ be the
semisimplification of $\bigoplus_{\nu\equiv g\mod \ 2} {\cal H}^\nu(M_0)$ and $B_0$ be the
semisimplification of $\bigoplus_{\nu\equiv g+1\mod \ 2} {\cal H}^\nu(M_0)$. Then
$C_0=\pi_0^*(A_0)[g]$ and $D_0=\pi_0^*(B_0)[g]$ are geometrically constant perverse sheaves on
$X_0$, such that
$$ f_m^{L_0}(x) + f_m^{D_0}(x) = f_m^{C_0}(x) $$
holds for all $m$ and all $x\in X_0(\kappa_m)$. Therefore $L_0\oplus D_0 \cong C_0$ by
Cebotarev. Hence $L$ is a summand of $C$, thus a geometrically constant perverse sheaf on
$X_0$. Factorize the isogeny $g:X \to X$ into a purely inseparable isogeny $u$ (division by the
connected component of the group scheme $X[N]$) and an etale isogeny $v$
$$ g: X \overset{u}{\longrightarrow} X' \overset{v}{\longrightarrow} X \ .$$
The isogeny $u$ is purely inseparable. The functors $u_*u^*$ and $u^*u_*$ define an equivalence
of $Perv(X)$ and $Perv(X')$ and map constant sheaves to constant sheaves. This easily follows
e.g. from \cite{FK} I.3.12. We deduce, that $L'=v^*(K)$ is a constant perverse sheaf on $X'$.
Then, since $v:X'\to X$ is a finite etale Galois morphism,  \cite{KW} III.15.3d implies
$$ K \hookrightarrow v_*v^*(L') \ .$$
Since $v_*v^*(L')$ is a direct sum of translation-invariant sheaves $E_\chi$ for characters
$\chi\in Hom_{cont}(\pi_1(X,0),\overline\Q_l^*)$ of finite order, which become trivial on the
etale cover $X'$ of $X$, hence $K$ is a translation-invariant perverse sheaf on $X$. This
proves the claim.

\bigskip\noindent
{\bf Definition}: A perverse sheaf $\delta_{E}$ on $\tilde C$, attached to a smooth coefficient
system $E$ on some Zariski open dense subset of $\tilde C_{reg}$, for which $E$ becomes trivial
on some finite branched covering $C'\to C$ of the normalization $C$ of $\tilde C$, will be
called an admissible perverse sheaf on $X$.

\bigskip\noindent
\begin{Lemma}\label{curvel} {\it Suppose $k$ is the algebraic closure of a finite field and  suppose $X$ is generated by
$\tilde C$. Then for a semisimple perverse Weil sheaf  $K_0$ on $X_0$ the convolution
$K*\delta_{\tilde C}$ is a direct sum  $$K*\delta_{\tilde C} \cong P \oplus  T\ ,$$ where $P$
is a semisimple perverse Weil-sheaf $P$ on $X$ and $T$ is a sum of complex translates of
translation-invariant perverse sheaves on $X$. Furthermore $T$ is perverse for an irreducible
perverse sheaf $K_0$ on $X_0$, which is not translation-invariant. In addition: The same
statement holds for $\delta_{\tilde C}$ replaced by any admissible  perverse sheaf $\delta_E$
on $X$. }
\end{Lemma}

\bigskip\noindent
\underbar{Proof}:  $\tilde C_0$, $C_0$ and $X_0$ are defined over the finite field $\kappa$. By
a finite base field extension of $\kappa$ we may assume, that $C_0$ contains a
$\kappa$-rational point $P_0$. By a translation of $K_0$ and $C_0$ we may assume $P_0=0$. Then
the Jacobian $J(C_0)$ of $C_0$ is defined over $\kappa$ and there exists a morphism $C_0\to
J(C_0)$ (depending on $P_0$) defined over $\kappa$. Consider the Albanese morphism $f_0:
J(C_0)\to X_0$ induced by the morphism $\pi_0:C_0\to \tilde C_0 \to X_0$. Since $\tilde C$
generates $X$, the map $f_0$ is a surjection of abelian varieties with kernel, say $K_0$,
$$  0 \to K_0\to J(C_0) \overset{f_0}{\to}\ X_0 \to 0 \ .$$
Put $X'_0= J(C_0)/(K_0)^0$. Then the map $J(C_0)\to X_0$ factorizes $J(C_0)\to X_0' \to X_0$,
where the second map is an isogeny
$$ 0 \to \pi_0(K_0) \to X_0' \overset{\varphi_0}{\to} X_0 \to 0 \ .$$
Since $(K_0)^0$ is Zariski connected, the isogeny $id-Frob_m : (K_0)^0(k) \to (K_0)^0(k)$ is
surjective (Lang's theorem). Therefore the homomorphism of finite groups
$$ J(C_0)(\kappa_m) \to X'_0(\kappa_m)$$
is surjective for all $m$. Finally one has an exact sequence
$$ X'_0(\kappa_m) \to X_0(\kappa_m) \to H^1(\kappa_m, \pi_0(K_0))
\ .$$ The order of the finite Galois cohomology group $H^1(\kappa_m, \pi_0(K_0))$ on the right
is bounded by a constant $C_1$ independently from $m$.

%\bigskip\noindent
%Let $C'$ be the image of $C$ in $X'$. Since
%$$ K\hookrightarrow \varphi_*(\varphi^*(K)) = K \otimes
%\bigoplus_\chi E_\chi $$ and $\varphi_*(\delta_{C'})=\delta_C$ it
%is enough to prove that $(\varphi^*(K))*\delta_{C'}$ is perverse.
%In fact, notice $\varphi^*(K)$ is pure and perverse and
%$$ \varphi_*(\varphi^*(K)*\delta_{C'}) = (\varphi_*\varphi^*(K))
%* \varphi_*(\delta_{C'})$$
%by the functoriality of the convolution product. Furthermore
%notice, that $\varphi^*(K)$ is not translation-invariant, if $K$
%is not translation-invariant. Hence we can replace $X$ by $X'$ and
%$K$ by $K'=\varphi^*(K)$ (respectively its irreducible
%constituents) and thereby assume $\pi_0(K)=0$.

\bigskip\noindent
Now assume $K_0$ is an irreducible perverse sheaf on $X_0$. We can assume by enlarging the base
field, that $K_0$ is geometrically irreducible. The assertion of the last lemma is trivial, if
$K$ is a translation-invariant perverse sheaf. So let us assume, that $K_0$ is not
geometrically translation-invariant.

\bigskip\noindent
Recall $K*\delta_{\tilde C} = Rf_*(K\boxtimes \Q_{l,C}[1]) $, and $f:X\times C\to X$ is smooth
with geometrically connected fibers of dimension 1. Hence $Rf_*(K\boxtimes \Q_{l,C}[1]) \in
{}^p D^{[-1,1]}(X)$, as already explained.  If the assertion of the lemma were false, then
either ${}^p H^1(Rf_*(K\boxtimes \Q_{l,C}[1])$ or ${}^p H^{-1}(Rf_*(K\boxtimes \Q_{l,C})[1])$
must be nontrivial.  Without restriction of generality assume, that we are in the first case
(the other case being analogous). Then  e.g. by \cite{KW} III.11.3 there exists a surjective
morphism of perverse sheaves on $X\times C$
$$ K\boxtimes \delta_{C} \twoheadrightarrow f^*[1] \ {}^p H^1\bigl(Rf_*(K\boxtimes
\Q_{l,C}[1])\bigr) \ .$$ Since $K$ is irreducible, $ K\boxtimes \delta_{C}$ is an irreducible
perverse sheaf on $X\times C$. Hence the last surjection defines an isomorphism of perverse
sheaves. Since it is already defined over $\kappa$, it can be viewed as an isomorphism of Weil
sheaf complexes. So one obtains an isomorphism in $Perv(X_0)$
$$ K_0\boxtimes \delta_{C_0}\ \cong\ f_0^*[1]\  {}^p H^1\bigl(Rf_{0*}(K_0\boxtimes
\Q_{l,C_0}[1])\bigr) \ .$$  For $z\in X(k)$ the geometric points of the fiber $f^{-1}(z)$ are
of the form $(z-\pi(y),y)$, where $y$ runs over the geometric points $C(k)$. The function
$f_m^{K_0}(x)$ for $x\in X_0(\kappa_m)$ therefore satisfies by property 3 and 4 of the function
$f_m(.)$
$$ f^{K_0}(z-\pi_0(y)) \cdot f_m^{\delta_{C_0}}(y) = -f_m^{L_0}(z) $$
for all $z\in X_0(\kappa_m)$ and all $y\in C_0(\kappa_m)$, where $L_0$ denotes the perverse
sheaf ${}^p H^1(Rf_{0*}(K_0\boxtimes \Q_{l,C_0}[1])$ on $X_0$. Now $f_m^{\delta_{C_0}}(y)=-1$
does not depend on $y\in C_0(\kappa_m)$. This implies for all $y\in C_0(\kappa_m)$ and all
$z\in X_0(\kappa_m)$ the equality
$$ f_m^{K_0}(z-\pi_0(y)) = f_m^{K_0}(z) \ .$$
But then also $ f_m^{K_0}(z-\pi_0(y_1)-\pi(y_2))=f_m^{K_0}(z-\pi(y_1))=f_m^{K_0}(z)$ and so on.
Let $\Gamma_m$ denote the subgroup of $X_0(\kappa_m)$ generated by the image of the points of
$C_0(\kappa_m)$. Then this implies $f_m^{K_0}(z+\gamma)=f_m^{K_0}(z)$ for all $z\in
X_0(\kappa_m)$ and all $\gamma\in \Gamma_m$. If we can show, that the index
$[X_0(\kappa_m):\Gamma_m]$ is bounded independently from $m$, then lemma \ref{l9} implies that
$K$ is a translation-invariant perverse sheaf on $X$ contradicting our assumption on $K$, which
then implies that $K*\delta_{\tilde C}$ must have been a perverse sheaf.

%But then let $M_{\pm}$ be the constant sheaves defined by the $(-1)^g$
%times even/odd stalk cohomology of ${}^p H^1(Rf_*(K\boxtimes
%\Q_{C,l}[1])$ at the point zero. We then consider the geometric
%constant perverse sheaves $M_\pm[g]$ and we get an identity
%$$ f^K(z) + f^{M_-}(z) = f^{M_+}(z) \ .$$
%If this holds for all $z\in X(\kappa_m)$ and for almost all $m\geq
%0$, then the Cebotarev density theorem implies (see \cite{KW}
%III.12.1 (3))
%$$ K \oplus M_- \cong M_+ \ ,$$
%hence $K$ is a geometrically constant perverse sheaf on $X$.
%
%\bigskip\noindent
%In fact let $J(C)$ be the Jacobian of $C$ and consider the
%Albanese morphism $f$ induced by $\pi$. Since $\tilde C$ generates
%$X$, the map $f$ is a surjection of abelian varieties with kernel
%say $K$
%$$  0 \to K\to J(C) \overset{f}{\to}\ X \to 0 \ .$$
%Then
%$$ J(C)(\kappa_m) \to X(\kappa_m) \to H^1(\kappa_m,K) \to $$
%Since now $K=K^0$ is connected $ H^1(\kappa_m,K)=0$ vanishes by

\bigskip\noindent
To control the index of the group $\Gamma_m$ generated by $\tilde C_0(\kappa_m)$ in
$X_0(\kappa_m)$ it is enough to bound the index of  the group generated by $C_0(\kappa_m)$ in
the Jacobian $J(C_0)(\kappa_m)$. If the index in the Jacobian is estimated by $C_2$, then by
the first reduction step  the index can is bounded by $[X_0(\kappa_m):\Gamma_m]\leq C_1C_2$.
Therefore we  now may assume $X_0=J(C_0)$ without loss of generality.

\bigskip\noindent
The quotient map $C_0^{2g} \to C_0^{(2g)}$ is a finite morphism  of degree $(2g)!$. Hence the
image of
$$ C_0^{2g}(\kappa_m) \to C_0^{(2g)}(\kappa_m)$$
contains $$ \geq q^{2gm}/(2g)! + o(q^{2mg})\ $$ different points.
%On the other hand a point in $J(C)(\kappa)$ represented by a
%divisor class is $\kappa_m$ rational if and only if the divisor
%class is $\kappa_m$-rational. Since the Brauer group of of a
%finite field vanishes, this is equivalent to the fact that the
%divisor class contains a $\kappa_m$-rational representing divisor
%$D\in C^{(2g)}(\kappa_m)$. Anyway, there are approximately
%There are $q^{gm}$-rational points in $J(C_0)(\kappa_m)$.
Over each point $z\in
J(C_0)(\kappa_m)$ the fiber of the morphism
$$ C_0^{(2g)} \to J(C_0) $$
is a Severi-Brauer variety, i.e. becomes a projective space of dimension $g$ over $k$. By a
result of Chatelet a Severi-Brauer variety is isomorphic over $\kappa_m$ to a projective space
over $\kappa_m$, if and only if it contains a $\kappa_m$-rational point. If the fiber contains
a $\kappa_m$-rational point $x\in C^{(2g)}(\kappa_m)$, then it contains $$ q^{gm}+q^{(g-1)m} +
\cdots +1=q^{gm}+o(q^{gm})
$$
such $\kappa_m$-rational points. In other words, there remain at least
$$ \frac{1}{(2g)!} q^{2gm}/q^{gm} + o(q^{gm}) \ =\ \frac{1}{(2g)!} q^{gm} + o(q^{gm})$$
different points in the image of $C_0^{2g}(\kappa_m)\to J(C_0)(\kappa_m)$, whereas
$$ \# J(C_0)(\kappa_m) = q^{gm} + o(q^{gm}) \ .$$
 Hence the subgroup generated by the image of
$C_0(\kappa_m)$ in $J(C_0)(\kappa_m)$ has index bounded by $C_2= (2g)!$, at least if $m\geq
m_0$ is chosen large enough. Concerning this we may replace $\kappa$ by a finite field
extension, so that this bound holds for all $m$. This completes the proof of the first
assertion.

%\bigskip\noindent
%Now suppose in addition that $K\in Perv(X)$ is perverse irreducible, but not translation
%invariant. Then by the proof of the first assertion ${}^p H^1(Rf_*(K\boxtimes \overline\Q_l[1])
%= {}^p H^1(P\oplus T)={}^p H^1(T)$. This implies as in the proof of the first assertion
%$K\boxtimes \delta_C \cong f^*[1]({}^p H^1(T))$. By a translation we may assume $0\in
%C(\kappa)$. Then $\nu^*(K\boxtimes\delta_C) \cong \nu^*f^*({}^p H^1(T))$. Hence $K\cong {}^p
%H^1(T)$, since $f\circ \nu=id_X$.

\bigskip\noindent To prove the statement for an additional coefficient system
$E$ on $C$ we notice, that the above argument did only  use, that the map $C\to X$ has an image
that generates $X$. Hence we could replace $C$ by a finite branched covering $\pi':C'\to C$,
which trivializes $E$, and $\delta_E\hookrightarrow \pi'_*(\overline \Q_{l,C'}[1])$ by
$\delta_{C'}=\overline\Q_{l,C'}[1]$ and study the map $f:X\times C'\to X$ in a similar way.
This completes the proof of the lemma.

\goodbreak
\bigskip\noindent
\begin{Theorem}\label{Geomorigin} {\it Suppose $k$ is an algebraically closed field of characteristic zero. Suppose $X$ is an
abelian variety over $k$ and $\tilde C$ is an irreducible projective curve contained in $X$
defined over $k$. Suppose $X$ is generated as an abelian variety by $\tilde C$. Then for a
semisimple perverse sheaf $K$ on $X$ of geometric origin the convolution $K*\delta_{\tilde C}$
is a direct sum
$$K*\delta_{\tilde C} \cong P \oplus  T\ ,$$ where $P$ is a
semisimple perverse sheaf $P$ of geometric origin on $X$ and $T$ is a sum of complex translates
of translation-invariant perverse sheaves on $X$ of geometric origin. Furthermore $T$ is
perverse for an irreducible perverse sheaf $K$ of geometric origin on $X$, which is not
translation-invariant. Moreover the same holds with $\delta_{\tilde C}$ replaced by some
admissible perverse sheaf $\delta_{E}$ on $\tilde C$.} \end{Theorem}

\bigskip\noindent
\underbar{Proof}: We may assume $k=\C$. Then we can apply the technique of \cite{BBD}, section
6.2 (in particular lemma 6.2.6 of loc. cit.) to reduce the proof of the theorem to the
analogous result proved over the algebraic closure of a finite field. Notice that all the
sheaves $\delta_E$ are sheaves of geometric origin.

\chapter{Jacobians} Let $X$ be the Jacobian of a smooth curve
$C$ of genus $g\geq 3$. Then, after choosing a point $P_0\in C(k)$, one obtains the Abel-Jacobi
period map $C\to X$, which induces maps $C^r=C\times \cdots \times C \to X$ by repeated
addition in the abelian variety $X$. These maps factorize over the smooth projective varieties
$C^{(r)}=C^r/\Sigma_r$ (quotient by the symmetric group $\Sigma_r$) $$ C^r
\overset{f}{\longrightarrow} C^{(r)} \overset{p}{\longrightarrow} X \ .$$
\bigskip\noindent
For increasing $d$ we have natural embeddings $C^{(d)} \hookrightarrow C^{(d+1)}$ making the
diagram
$$ \xymatrix{ ...\quad  \ar@{^{(}->}[r] & \ C^{(d)}\ \ar@{^{(}->}[r]\ar[dr]_{p_d} & \ C^{(d+1)}\ \ar@{^{(}->}[r]\ar[d]^{\ p_{d+1}} & \ C^{(d+2)}\ \ar@{^{(}->}[r] \ar[dl]^{\ p_{d+2}} & \quad ... \cr & & X & & \cr
}
$$
commutative. For $d>2g$ the map $p=p_d$ $$p_d: C^{(d)} \to X$$ becomes a projective bundle
morphism, hence in particular a projective morphism.  Hence all $p_d$ for all $d\geq 1$ are
projective morphisms by the diagram above. For $d<g$ the image of $p_d$ is the subvariety
$W_d\hookrightarrow X$. For details see \cite{ACGH} and \cite{Mi} \S 5 in the case of positive
characteristics.

\bigskip\noindent
The direct image under $f$ respectively $p$ of the constant perverse sheaf
$\lambda_{C^r}=\delta_{C^r}$ therefore is an $\Sigma_r$-equivariant sheaf complex
$$ f_*(\Q_{l,C^r}[r]) \ =\ \bigoplus_\alpha \sigma_\alpha \otimes {\cal F}_\alpha \ , $$
where $\sigma_\alpha$ ranges over a representative system of equivalence classes of irreducible
$\overline \Q_l$-representations of the symmetric group $\Sigma_r$. Then up to a shift by $r$
the ${\cal F}_\alpha$ are sheaves on $C^{(r)}$, but they are indeed also in $Perv(C^{(r)})$.
Therefore
$$ R(p\circ f)_*(\delta_{C^r}) = Rp_*(\bigoplus_\alpha \sigma_\alpha \otimes {\cal F}_\alpha
) \ = \ \bigoplus_\alpha\ \sigma_\alpha \otimes \delta_\alpha $$
for
$$ \delta_\alpha = Rp_*({\cal F}_\alpha) \ .$$
Here $\delta_\alpha$ is a sheaf complex on $X$, which is pure of
weight $w=r$. Notice, that every irreducible representation of the
symmetric group is defined over $\Q$. This implies, that the
Verdier dual $D({\cal F}_\alpha)$ is isomorphic to ${\cal
F}_\alpha$. As a consequence we obtain

\bigskip\noindent
\begin{Lemma} {\it $\delta_\alpha$ is a self dual pure complex on $X$ of weight $r=deg(\alpha)$.}
\end{Lemma}

\section{The multiplication law}\label{ML}
From the diagram for $r=r_1+r_2$ and the relative K\"unneth formula
$$ \xymatrix{ & C^{(r)}\ar[r]^-p & X \cr
C^r \ar[ru]^f \ar[rd]_-{f\times f} & &  \cr
 & C^{(r_1)}\times
C^{(r_2)}\ar[r]^-{p\times p}\ar[uu]^-{\tau}  & X\times X \ar[uu]_a
\cr}
$$
it follows, that $\tau_*(\bigoplus \sigma_\alpha \otimes{\cal
F}_\alpha \boxtimes \bigoplus_\beta \sigma_\beta \otimes {\cal
F}_\beta) = \bigoplus_\gamma \sigma_\gamma \otimes {\cal
F}_\gamma$. Hence $\bigoplus \bigoplus_\beta
(\sigma_\alpha\boxtimes \sigma_\beta) \otimes \tau_*({\cal
F}_\alpha \boxtimes \otimes {\cal F}_\beta) = \bigoplus_\gamma
\sigma_\gamma\big\vert_{\Sigma_{r_1}\times \Sigma_{r_2}} \otimes
{\cal F}_\gamma$. Let $m_{\alpha\beta}^\gamma$ denote the
dimension of the space of intertwining operators
$$Hom_{\Sigma_r}(\sigma_\gamma, Ind_{\Sigma_{r_1}\times
\Sigma_{r_2}}^{\Sigma_r}(\sigma_\alpha\boxtimes \sigma_\beta)) =
Hom_{\Sigma_{r_1}\times
\Sigma_{r_2}}(\sigma_\gamma\big\vert_{\Sigma_{r_1}\times
\Sigma_{r_2}}, \sigma_\alpha\boxtimes \sigma_\beta)\ .$$ These
integers have a combinatorial description (this is the
Littlewood-Richardson rule described in section \ref{LRrule}). We
conclude

\bigskip\noindent
\begin{Lemma}\label{mult} {\it $\delta_\alpha * \delta_\beta = \bigoplus_\gamma m_{\alpha\beta}^\gamma \cdot
\delta_\gamma$, where $m_{\alpha\beta}^\gamma$ are the Littlewood-Richardson coefficients.}
\end{Lemma}

\bigskip\noindent
\section{Cohomology groups}\label{Cohomology groups}
Weyl's theorem describes, how the the $r$-th tensor representation  of the standard
representation of $Gl(V)$ on $V$ (for a finitely dimensional vector space $V$ over an
algebraically closed field) decomposes into irreducible constituents $V^\alpha$. The group
$\Sigma_r\times Gl(V)$ acts on the tensor product $V^{\otimes r}$ in a natural way, where the
symmetric group $\Sigma_r$ permutes the factors of the tensor product (we call this the Weyl
action). Then Weyl's theorem asserts a decomposition
$$ V^{\otimes r} = \bigoplus_\alpha \sigma_\alpha \boxtimes
V^\alpha \ ,$$ where summation is over a set of representatives for the irreducible
representations $\sigma_\alpha$ of the symmetric group $\Sigma_r$. The $V^\alpha$ are
irreducible algebraic representations of the group $Gl(V)$ of highest weight
$\alpha=(\alpha_1,...,\alpha_r)$ respectively, where the $\alpha_i$ are integers such that
$\alpha_1 \geq \alpha_2 \geq ... \alpha_r\geq 0$ such that $deg(\alpha)=\sum \alpha_i =r $. The
$\alpha$'s correspond to the partitions of the number $r$.

\bigskip\noindent
 Let denote
$$ H_+ = H^0(C,\overline \Q_l)[1] \oplus H^2(C,\overline \Q_l)[-1]
$$
$$ H_- = H^1(C,\overline \Q_l) \ .$$
There is a natural action of the group $\Sigma_r$ on the product
$C^r$. This induces a natural action of $\Sigma$ on the cohomology
of $C^r$. By the K\"{u}nneth theorem $$ \iota: \ \bigotimes^r
H^\bullet(C,\overline \Q_l)[1] \overset{\sim}{\longrightarrow}
H^\bullet(C^r,\overline \Q_l)[r] \ .$$ The isomorphism $\iota$ is
given by the cup-product
$$  x_1\otimes ...\otimes x_r \mapsto pr_1^*(x_1)\cup ...\cup
pr_r^*(x_r) \ .$$ Notice, that the natural action of the group $\Sigma_r$ on the tensor product
on the left (the Weyl action) does not coincide with the action, which is induced on cohomology
by the permutation of factors of $C^r$ (the cohomological action). This is due to the fact,
that the cup-product is alternating on $H_-$ and symmetric on $H_+$. Since cup-products of
elements in $H_+$ and $H_-$ commute, for both actions we get
$$ \bigotimes^r (H_+\oplus H_-) = \bigoplus_{a+b=r}
\bigoplus_{\tau\in \Sigma_r/\Sigma_a\times \Sigma_b}(\bigotimes^a H_+ \otimes \bigotimes^b H_+)
\ $$ For both actions we therefore get a decomposition of the representation of the symmetric
group $\Sigma_r$ into the form
$$ \bigoplus_{a+b=r} Ind_{\Sigma_{a}\times \Sigma_{b}}^{\Sigma_r}
(\bigotimes^a H_+ \otimes \bigotimes^b H_-) \ .$$ Written in this
way it is now obvious, how the two action differ.

\bigskip\noindent
\begin{Lemma} {\it The two representations of $\Sigma_b$ on $\bigotimes^b H_-$, induced by the cohomological
action respectively the Weyl action, differ by a twist with the signum character.} \end{Lemma}

\bigskip\noindent
\underbar{Proof}: It is enough to consider involutions $\sigma$ of
$C^r$, which permute two factors $i\neq j$. Then for 1-forms
$u,v\in H_-$ we have $\sigma_{Weyl}(...\otimes u\otimes \cdots
\otimes v \otimes ...) = (... \otimes v\otimes \cdots \otimes u
\otimes ...)$. To compare $\sigma_{cohom} (... \cup pr_i^*(u)\cup
\cdots \cup pr_j^*(v) \cup ...) = (... \cup pr_j^*(u)\cup \cdots
\cup pr_i^*(v) \cup ...) = - (... \cup pr_i^*(v)\cup \cdots \cup
pr_j^*(u) \cup ...) $. $\quad \quad\square$

\bigskip\noindent
This being said, we obtain from Weyl's formula the following expression for the cohomological
action of $\Sigma_r$
$$  H^\bullet(C^r,\overline \Q_l)[r] \ = \ \bigoplus_{a+b=r} Ind_{\Sigma_{a}\times \Sigma_{b}}^{\Sigma_r}
\Bigl( (\!\bigoplus_{deg(\alpha)=a}\! \sigma_\alpha\boxtimes
H_+^\alpha) \boxtimes (\!\bigoplus_{deg(\beta)=b}
(\sigma_\beta\otimes sign)\boxtimes H_-^\beta)\Bigr) \ .$$ The sum
ranges over all partitions $\alpha$ of $a$ respectively  over all
partitions $\beta$ of $b$. For an irreducible representation
$\sigma$ of the symmetric group $\Sigma_b$ let now $\sigma^*$
denote the representation $\sigma\otimes \mbox{sign}$. Then
$$ H^\bullet(C^r,\overline \Q_l)[r] = \bigoplus_{a+b=r\ }
\ \bigoplus_{deg(\alpha)=a,deg(\beta)=b} Ind_{\Sigma_{a}\times
\Sigma_{b}}^{\Sigma_r}(\sigma_\alpha\boxtimes \sigma_\beta^*)
\boxtimes (H_+^\alpha \otimes H_-^\beta) \ .$$ If we decompose
$H^\bullet(C^r,\overline \Q_l)[r]$ into its irreducible
constituents $\bigoplus_\gamma \sigma_\gamma \boxtimes
H^{\bullet}(X,\delta_\gamma)$ under the cohomological action of
the group $\Sigma_r$, i.e. the sum runs over all partitions of
$r$, we obtain

\bigskip\noindent
\begin{Lemma} \label{dimension}{\it For $deg(\gamma)=r$ we have
$$H^{\bullet}(X,\delta_\gamma) = \bigoplus_{a+b=r}
\bigoplus_{deg(\alpha)=a,deg(\beta)=b} m_{\alpha \beta^*}^\gamma \cdot  H_+^\alpha \otimes
H_-^\beta\ .$$} \end{Lemma}

\bigskip\noindent
\underbar{Remark}: Since $H_+$ has dimension $2$, in this sum an
$\alpha=(\alpha_1,\alpha_2,\alpha_3,..)$ contributes only if $\alpha_i=0$ holds for all $i\geq
3$.

\bigskip\noindent
\begin{Corollary} {$deg_{t^{1/2}}(h(\delta_\gamma,t))\leq \gamma_1$ and equality
$$deg_{t^{1/2}}(h(\delta_\gamma,t))\ = \ \gamma_1$$ holds for $\gamma_1\leq 2g$.}
\end{Corollary}

\bigskip\noindent
\underbar{Proof}: This is immediate consequence of the
Littlewood-Richardson rules, which forces $\alpha_1 \leq \gamma_1$
for $m_{\alpha\beta}^{\gamma}\neq 0$. For $\gamma_1 \leq 2g$ the
leading coefficient in $t$ is obtained from the summand
$\alpha=(\gamma_1)$ and $\beta = (\gamma_2,\gamma_3,...)$. It
gives
$$ h(\delta_\gamma,t) = m_{(\gamma_1)\ (\gamma_2,\gamma_3,..)}^\gamma \cdot
dim \bigl( H_-^{(\gamma_2,...)^*}\bigr) \cdot t^{\gamma_1/2} \ + O(t^{(\gamma_1-1)/2}) \ .$$
Since $\gamma_1\leq 2g$, the dimension of $H_-^{(\gamma_2,...)^*}$ is nonzero by the Weyl
formulas. In fact the coefficient $m_{(\gamma_1)\ (\gamma_2,\gamma_3,..)}^\gamma$ can shown to
be one.

\bigskip\noindent
\underbar{Notation}: For $\gamma=(r,0,...0)$ we also write $\delta_r $ instead of
$\delta_\gamma$.

\bigskip\noindent
\underbar{Example}:  Since for $\gamma=(r,0,...0)$ the coefficients  $m_{\alpha
\beta^*}^\gamma$ vanish except for $\alpha=(a,0,..,0)$ and $\beta^*=(b,0,..,0)$, we obtain
$$ H^{\bullet}(X,\delta_r) = \bigoplus_{a+b=r} S^a(H_+)\otimes
\Lambda^b(H_-) \ , $$ where $S^a$ and $\Lambda^b$ denotes the
$a$-th respectively $b$-th symmetric respectively antisymmetric
powers. Hence $ h(\delta_r,t) = \sum_{i=0}^r [i+1]_t \cdot
\Lambda^{r-i}(H_-)$.

\bigskip\noindent
\underbar{Degree $g-1$}: Suppose $C$ is not hyperelliptic. Then as a special case we obtain
$IH^{\bullet}(\Theta)\cong H^{\bullet}(X,\delta_{g-1})$. In the general formula only $a=2g-3$
and $b=1$ contributes to $IH^{2g-3}(\Theta)$. Hence (ignoring Tate twists) we obtain
$IH^{2g-3}(\Theta) \cong H^1(C)$. In particular, since $IH^1(\Theta)\cong IH^{2g-3}(\Theta)$

\bigskip\noindent
\begin{Corollary} \label{cor8} $IH^1(\Theta) \cong H^1(C) \cong H^1(X)$. \end{Corollary}

\bigskip\noindent
\underbar{Remark}: The statement of the corollary remains true for
hyperelliptic curves, although in this case $\delta_\Theta =
\delta_{g-1} - \delta_{g-3}$ as will be shown in section
\ref{Hyperelliptic}. More generally  $\delta_{W_r} = \delta_{r} -
\delta_{r-2}$ for all $r\leq g-1$ in the hyperelliptic case. Using
this, one can argue as above to show $IH^{2g-3}(\Theta)\cong
H^1(C)$. In fact $IH^1(W_r)\cong IH^{2r-1}(W_r) =
H^{r-1}(X,\delta_r)\cong H^1(C)$ for all $r\leq g-1$, since
$\delta_{W_r} \oplus \delta_{r-2} = \delta_r$ and
$H^{r-1}(X,\delta_{r-2})=0$ by the formulas above.

\bigskip\noindent
\underbar{Trivial cases}: If $\alpha=0$, then $\delta_\alpha=\delta_0=1$ is the unit element of
the convolution product. If $\alpha=(1,0,..)$, then $\delta_1=\delta_C$.

\goodbreak
\section{Theorem of Martens} \label{MT}
The fibers of the map $p_d:C^{(d)} \to X$ are projective spaces. For $r\geq 1$ the subvariety
$W^r_d\subseteq X$ is the locus of points $x$ in the image of $p_d$, where $p_d^{-1}(x)$ is a
projective space of dimension $\geq r-1$. For $r=1$ the image $W_d=W_d^1$ of $C^{(d)}$ has
dimension $d$, if $1\leq d\leq g$.

\bigskip\noindent
Suppose $r\geq 2$ and $2\leq d\leq g-1$. Then, for a smooth curve
of genus $g\geq 3$ the following holds

\bigskip\noindent
\begin{Theorem} (Martens): {\it
Suppose $ 2(r-1)\leq d$. Then \begin{enumerate} \item If $C$ is not a hyperelliptic curve $$
dim(W^r_d) +2(r-1) \leq d - 1 \ .$$ \item  If $C$ is hyperelliptic, then $W^r_d=W_{d-2r+2}$
holds up to translation. Therefore $$ dim(W^r_d) + 2(r -1) = d\ .$$\end{enumerate}}
\end{Theorem}

\bigskip\noindent
\underbar{Proof}: See [M] and \cite{Gunning}, p.55 and p.57. For hyperelliptic $C$ there exists
a point $e\in X$ such that $$e-W_1 = W_1\ ,$$ and more generally $ W^r_d = W_{d-2r+2} -
(r-1)e$.

\bigskip\noindent
\begin{Corollary} {\it Suppose $g\geq 3$. Then for $1\leq d\leq g-1$ the morphism $p_d: C^{(d)} \to W_d$ is
semi-small. If $C$ is not hyperelliptic $p_d$ is small.} \end{Corollary}

\bigskip\noindent
\underbar{Proof}: Obviously for $r\geq 2$    $$ 2\cdot dim(p_d^{-1}(x)) \ =\ 2(r-1)\ < \ d - (d
- 2r + 2) \ \leq \ dim(W_d) - dim(W_d^r)$$ for $x\in W_d^r \setminus W^{r+1}_d$. Notice we can
assume $2(r+1)\leq d$, since otherwise $W^{r}_d$ is empty. Hence $p_d:C^{(d)}\to W_d$ is
semi-small. Furthermore $p_d$ is an isomorphism over $W_d\setminus W^2_d$. Similarly, for $C$
not hyperelliptic, we have a strict inequality above. Hence $p_d$ is small. This improves
corollary \ref{pingpong}. In fact we obtain

\bigskip\noindent
\begin{Corollary} \label{MA} {\it For degree $0\leq deg(\alpha)\leq g-1$ the complexes $\delta_\alpha$ are pure perverse
sheaves on $X$ of weight $-deg(\alpha)$. Suppose $C$ is not hyperelliptic. Then
\begin{enumerate} \item[(i)] These sheaves are irreducible perverse
sheaves\footnote{One could go back from this to view Marten's theorem as a consequence of
\ref{pingpong} and of the main result of \ref{Rep}, whose proof is basically independent.} with
support in $W_{deg(\alpha)}$.
\item[(ii)] $\delta_\alpha \cong \delta_\beta$ for these sheaves if and only if $\alpha =
\beta$.
\item[(iii)] $\delta_d = \delta_{W_d}$ is the intersection cohomology sheaf
of $W_d$ for $1\leq d\leq g-1$, hence
$$ IH^\bullet(W_d,\overline\Q_l)[d] = \bigoplus_{a+b=d}
S^a(H_+)\otimes \Lambda^b(H_-) \ .$$ \end{enumerate}} \end{Corollary}

\bigskip\noindent
\underbar{Proof}: Concerning property (i) see [KW] III.7.4 and
7.5(ii). Notice that III.7.5 does not require ${\cal G}$ to be a
smooth sheaf. In fact it is enough, that ${\cal G}$ in the
notations of loc. cit. is an etale sheaf on $X$ in the ordinary
sense, such that ${\cal G}[n]$ is a pure perverse sheaf on $X$.
Then $Rf_!{\cal G}[n] \in {}^p D^{\leq 0}(Y)$ holds for projective
morphisms $f:X\to Y$ by [KW], III.7.4 and $Rf_!{\cal G}[n] \in
Perv(Y)$ by the hard Lefschetz theorem [KW] IV.4.1. Concerning
property (ii) see [KW] III.15.3(c).

\bigskip\noindent
Obviously for $d<g$
$$     supp({\cal H}^{-d+2r}(\delta_d)) \ = \ W_d^r  \ .$$
Similar, the supports of the cohomology sheaves of the $\delta_C^{*d}$ can be related to the
spaces $W_d^r$ in an obvious way. Moreover

\begin{Lemma} Suppose $C$ is a general curve of genus $g$ and suppose $char(k)=0$. Then for $d<g$
$$ {\cal H}^{-d+2r}(\delta_d) \ \cong \ \overline\Q_{l,W_d^r} \ .$$
\end{Lemma}

\bigskip\noindent
\underbar{Proof}: By \cite{ACGH} p.215 (.7) and p.163 (1.6) it follows, that $C_d^r$ is smooth
of the correct dimension in the sense of loc. cit. Also by \cite{ACGH} p. 190 (4.4) the variety
$W_d^r$ is Cohen-Macaulay, reduced and normal. Since $W_d$ is normal, the sheaf ${\cal
H}^{-d}(\delta_d)$ is $\overline\Q_{l,W_d}$. Hence the restriction to $W_d^r$ is ${\cal
H}^{-d}(\delta_d)\vert_{W_d^r} = \overline\Q_{l,W_d^r}$. The natural map $\pi: C^r_d \to W_d^r$
is projective, and it is smooth of relative dimension $d$ over the dense open subset
$W_d^{r}\setminus W_d^{r+1}$ with projective spaces of dimension $r$ as fibers. Therefore the
hard Lefschetz theorem implies ${\cal H}^{-d}(R\pi_*\overline\Q_{l,C^r_d}) \cong {\cal
H}^{-d+2}(R\pi_*\overline\Q_{l,C^r_d})\cong \cdots \cong {\cal
H}^{-d+2r}(R\pi_*\overline\Q_{l,C^r_d})\cong \overline\Q_{l,W_d^r}$. Since
$R\pi_*(\overline\Q_{l,C^r_d}) = Rp_{d,*}(\overline\Q_{l,C^{(d)}})\vert_{C_d^r}$, this proves
the claim.

\section{The adjoint perverse sheaf}\label{CundA}
The roles of $\delta_C\in Perv(X)$ and $\delta_C^\vee=\delta_{-C}\in Perv(X)$ are dual to each
other
$$ \fbox{$ \delta_C * \delta_{-C}\ =\ 1 \ \oplus \ \Omega $}\ , $$ where $1=\delta_{0}$.
For the difference map $f(x,y)=x-y$
$$ f: C\times C \to X \ $$
by definition
$$ Rf_*(\overline\Q_{l,C\times C})[2] \ \cong \ \delta_C *
\delta_{-C} \ .$$ To compute the stalks of the left side we have to consider the intersection
$C\cap (z+C)$. For $z=0$ this intersection is $C$. So it remains to consider  the case $z\neq
0$. For this we distinguish two cases (see \cite{Gunning}, p.128):

\bigskip\noindent
\underbar{Intersections}: For $z\neq 0$ suppose $x\in C$ and $x\in z+C$. Then $z=x-y$ for
$x,y\in C$, hence $z\in C-C$. Suppose we have two solutions $(x_1,y_1)\in C^2(k)$ and
$(x_2,y_2)\in C^2(k)$ such that
$$ z=x_1-y_1=x_2-y_2\ .$$ Then $x_1+y_2=y_1+x_2$.

\bigskip\noindent
\underbar{Non-hyperelliptic case}: If $C$ is not hyperelliptic, $H^0(C,O_C(x_1+y_2))$ has
dimension 1. Therefore $x_1=y_1$ or $x_1=x_2$. In the first case $z=0$, but we ssumed $z\neq
0$. Hence $x_1=x_2$ and $y_1=y_2$. It follows, that in the non-hyperelliptic case $C\cap (z+C)$
is nonempty if and only if $z\in C-C$, such that the intersection $C\cap (z+C)$ consists of one
unique point for $z\neq 0$, respectively  $C\cap (z+C)=C$ for $z=0$.

\bigskip\noindent
Hence $f: C\times C\to C-C$, defined by $f(x,y)=x-y$, is birational. The diagonal $\Delta_C$ is
blown down to $0\in X$. Outside the inverse image of $\{0\}$ the map $f$ is an isomorphism.
This easily implies
$$ \fbox{$ \Omega=\delta_{C-C} $}\ .$$

\bigskip\noindent
\underbar{Hyperelliptic curves}:  In this case the intersection $C\cap (C+z)$ for $z\neq 0$
consists of two points. In fact, there is a unique isomorphism class ${\cal L}$ of line bundles
of degree 2  so that $H^0(C,{\cal L})=2$. The divisors in this class have the form $D= x +
\theta(x)$, where $\theta:C\to C$ denotes the hyperelliptic involution. The quotient $C/\theta$
is the projective space $\Proj^1$.  $f(x_1,y_1)=f(x_2,y_2)$ either implies $y_2=\theta x_1$ and
$y_1=\theta x_2$, or $(x_1,y_1)=(x_2,y_2)$. Hence $(x_2,y_2)=(\theta y_1,\theta x_1)$ or
$(x_2,y_2)=(x_1,y_1)$. This implies:

\bigskip\noindent
The diagonal $\Delta_C\hookrightarrow C\times C$ is contracted by $f$ to the point zero.
Outside $f^{-1}(0)$  the morphism $f$ defines a two fold ramified covering with the covering
automorphism
$$\sigma(x,y)= (\theta(y), \theta(x))\ .$$ The covering
automorphism acts on the points of the diagonal $\Delta_C$ by $(x,x)\mapsto (\theta x,\theta
x)$. Its fixed points of the covering automorphism are of the form $(x,\theta(x))$, and the
locus of these points is the branch locus. Furthermore $f$ is a semi-small morphism. This
together with Gabber's theorem implies
$$ \fbox{$ \delta_C*\delta_{-C} \cong Rf_*(\delta_{C\times C}) \cong \delta_{C-C}\oplus \delta_E \oplus
S $}
 $$ for a coefficient system $E=E(\rho)$ on
$C-C$ of generic rank one, which corresponds to the nontrivial quadratic character $\rho$ of
the fundamental group $\pi_1((C-C)\setminus V)$, which is defined to the branched covering $f$.
The ramification locus of $E$ is contained in the image $V$ of the branch locus. Finally $S$ is
a a perverse sheaf with support in $\{0\}$, which is easily determined by considering the fiber
cohomology $H^\bullet(f^{-1}(0))= H^\bullet(C)$. In fact
$$ \fbox{$ S=\delta_{\{0\}} $} $$ coming from $H^2(C)=\overline\Q_l$. $R^{-1}f_*(\delta_{C\times C})$ is
a skyscraper sheaf with support in $\{0\}$. We claim $$R^{-1}f_*(\delta_{C\times C}) ={\cal
H}^{-1}(A)$$ and ${\cal H}^{-1}(\delta_{C-C})=0$. Proof: The group $\Z/2\Z$ generated by the
involution  $\sigma$ acts on $C\times C$, and $f$ is equivariant with respect to the trivial
action of $\Z/2\Z$ on $X$. The automorphism $\sigma$ induces the automorphism $\theta$ on the
fiber $f^{-1}(0)\cong C$
$$ \xymatrix@+0,5cm{C\ \ar[d]_\theta\ar[r]^-{id\times id} & \ C\times C \ar[d]_-\sigma\ar[r]^-f & C-C \subseteq X \ar[d]_{id_X}\cr
C\ \ar[r]^-{id\times id} & \ C\times C \ar[r]^-f & C-C \subseteq X \cr}
$$
Therefore the complex $Rf_*(\delta_{C\times C})$ is an $\Z/2\Z$-equivariant perverse sheaf on
$X$. Since $\Z/2\Z$ acts trivially on $X$, it decomposes into a $+$-eigenspace and a
$-$-eigenspace: The $+$-eigenspace is
$$ \delta_{0} \oplus \delta_{C-C} \ ,$$
the $(-1)$-eigenspace is $A$. To prove our claim it is therefore enough to show, that $\sigma$
acts nontrivially on the stalk
$$ (R^1f_*(\delta_{C\times C}))_0 \ \cong \ H^1(C) \ .$$
Since $f^{-1}(0)=\Delta_C\cong C$ is the diagonal  $\Delta_C\hookrightarrow C\times C$, and
since $\sigma$ restricts to the automorphism $\theta$ on the fiber $f^{-1}(0)$, this follows
since the hyperelliptic involution $\theta$ acts on $H^1(C)$ by $-id$. Hence
$$ {\cal H}^{-1}(A)\ \cong\ H^1(C)\otimes \delta_{0} \quad , \quad
 {\cal H}^{-1}(\delta_{C-C})\ \cong\ 0 \ .$$

\bigskip\noindent
\underbar{Resume}: As a final result in the hyperelliptic case $\Omega=\delta_{C-C} \oplus A $,
where $A=\delta_E$ for the nontrivial coefficient system $E$ on $C-C$ described above.
%
%\bigskip\noindent
%In any case, whether $C$ is hyperelliptic or not, we get
%$$ \xymatrix{ K \quad \ar@{|->}[r]^-{*\!\delta_C} & \quad  L=K\! *\! \delta_C\
%\ar@{|->}[r]^-{*\!\delta_{-C}} & \quad K\ \oplus\ (K\! *\! \Omega) \cr} \ .$$

\section{Hyperelliptic curves I} \label{Hyperelliptic}
Suppose $C$ is hyperelliptic and $X$ is its Jacobian. Then as in the last section
$$ \xymatrix@+0,5cm{C\ \ar[d]_\theta\ar[r]^-{id\times \theta} & \ C\times C \ar[d]_-\tau\ar[r]^-a & C+C \subseteq X \ar[d]_{id_X}\cr
C\ \ar[r]^-{id\times \theta} & \ C\times C \ar[r]^-a & C+C \subseteq X \cr}
$$
The addition morphism $a$ induces  a branched covering outside $a^{-1}(\{e\}) \cong C
=image(id\times \theta: C\to C\times C)$ with the covering automorphism defined by
$$\tau(x,y)=(y,x)\ .$$ Recall, $\theta(x)+x=e$ for all $x\in C$ implies
$$  C = e - C $$
in the hyperelliptic case. Therefore
$$ Ra_*(\delta_{C\times C}) = \delta_C*\delta_C = \delta_2\oplus \delta_{1,1} = \delta_{C}*\delta_{e-C} =\delta_{\{e\}} \oplus
\delta_{e+C-C} \oplus T_{-e}^*(A)  \ ,$$ which is the translate of $\delta_C*\delta_{-C} =
\delta_0 \oplus \Omega= \delta_{0}\oplus \delta_{C-C}\oplus A$ by $e$. Considered as an
equivariant sheaf with respect to the action of the group $\{id,\tau\}$ with trivial action on
$X$, the eigenspace decomposition gives
\begin{enumerate} \item{ $\oplus$-eigenspace:}  $\quad \delta_2 = \delta_{\{e\}} \oplus \delta_{e+C-C}\ =\
\delta_{\{e\}} \oplus \delta_{C+C}  $ \item{$\ominus$-eigenspace:} $ \quad \delta_{1,1} =
T_{-e}^*(A)$
\end{enumerate}
Hence $H^\bullet(X,A)\cong H^\bullet(X,\delta_{1,1})$, thus  $dim(H^{\pm 1}(X,A))=2g$ and
$dim(H^{0}(X,A))=dim(S^2(\overline\Q_l^{2g}))+1$. If we consider the symmetric quotient
$C^{(2)}=(C\times C)/\Sigma_2$ for $\Sigma=\{id,\tau\}$, then
$$ \xymatrix@+0,5cm{C\ \ar[d]\ar[r]^-{id\times \theta} & \ C\times C \ar[d]\ar[r]^-a & C+C \subseteq X \ar[d]_{id_X}\cr
\Proj^1\ \ar[r] & \ C^{(2)} \ar[r]^-p & C+C \subseteq X \cr}
$$
and $p$ is a birational map, which blows down the fiber $p^{-1}(e)=\Proj^1$. Hence
$$Rp_*(\delta_{C^{(2)}}) = \delta_{e+C-C} \oplus \delta_{\{e\}} =  \delta_{C+C} \oplus \delta_{\{e\}}
 =  \delta_{W_2} \oplus \delta_{\{e\}} \
.$$

%\section{Surfaces}
%Let $A$ be a surface, i.e. an irreducible variety of dimension
%two. Suppose $A$ has isolated singularities. Then the intersection
%cohomology sheaf $\delta_A$  has the following cohomology sheaves
%\begin{enumerate}
%\item ${\cal H}^{-2}(A)$ is the constant sheaf.
%\item ${\cal H}^{-1}(A)$ is a skyscraper sheaf with support
%in the finitely many point of $A_{sing}$.
%\end{enumerate}
%and all other cohomology sheaves vanish.
%
%\bigskip\noindent
%Then with the notations of \ref{} we have $ \psi_E[1] =
%\tau^{st}_{\geq -1}\delta_A \cong {\cal H}^{-1}(A) $. Since this
%is a skyscraper sheaf $$\psi_A = {\cal H}^{-1}(A)$$ and $ \psi_E
%\in Perv(A)$. Hence we get an exact sequence in $Perv(A)$
%$$ 0 \to \psi_A \to \lambda_A \to \delta_A \to 0 \ .$$
%Therefore from the long exact sequence
%$$ 0\to  H^1(A) \to IH^{1}(A) \to Kern(H^0(\psi_A) \to H^2(A)) \to
%0  $$ Let $\pi:X\to A$ be a minimal desingularisation of $A$. This
%is a semi-small morphism. By Gabber's decomposition theorem,
%therefore $$ R\pi_*\overline\Q_l[2] = \delta_A \oplus D \ ,$$
%where $D$ is a perverse skyscraper sheaf with support in the
%singular points $A_{sing}$. Hence
%$$ \psi_A \cong R^1\pi_*\overline\Q_l) \ .$$
%Let $Y$ be the exceptional divisor, and $U\subseteq X$ is
%complement. Then by [KW], III.10.7 there exists an isomorphism
%$$ H^1(X) \cong H^1(U) \ .$$
%By the Leray spectral sequence we again get an exact sequence
%$$ 0\to  H^1(A) \to H^{1}(X) \to Kern(H^0(R^1\pi_*\overline\Q_l)) \to H^2(A)) \to
%0 $$ From the results above we conclude $$ IH^1(A)\cong H^1(X)\ $$

\bigskip\noindent
More generally, by the theorem of Martens  $ W^r_d = W_{d-2r+2} - (r-1)e$. Since the map $p_d:
C^{(d)}\to X$ is semi-small for $1\leq d\leq g-1$, we get for  $1\leq d \leq g-1$
$$ \delta_d = \delta_{W_d} \oplus \delta_{W_{d-2}-e} \oplus
\delta_{W_{d-4}-2e} \oplus \cdots \ $$ hence

\bigskip\noindent
\begin{Lemma}\label{hy} {\it  For hyperelliptic curves the sheaves $\delta_d$ are perverse sheaves for $1\leq d\leq
g-1$ such that (in the Grothendieck group)
$$ \fbox{$ \delta_{W_{d}} = \delta_d - \delta_{d-2} $} \ .$$
Furthermore ${\cal H}^{-\nu}(\delta_{W_{d}})=0$ for $\nu\neq -d$, hence $$ \fbox{$ \delta_{W_d}
= \overline\Q_{l,W_d}[d] $} \ .$$} \end{Lemma}

\bigskip\noindent
\underbar{Singular supports}: Since $C=e-C$ now
$$ W_d = d\cdot e - W_d $$
for all $d\leq g-1$. For simplicity of notation we now assume
$e=0$, which can be achieved by replacing $C$ with an appropriate
translate of $C$ in $X$. This convention will be used from now on
in this section. Then  $\delta_C=\delta_{-C}$ and
$$ W_d = - W_d \quad , \quad 0\leq d\leq g-1 \ .$$

\bigskip\noindent
\begin{Lemma} \label{hypell}{\it For $0\leq d\leq g-1$
$$ {\cal H}^{0}(\delta_{W_d}*\delta_{W_d}) = \delta_0 $$
$$ {\cal H}^{-1}(\delta_{W_d}*\delta_{W_d})= IH^1(W_d)\otimes \delta_0 = H^1(C)\otimes \delta_0 \
.$$ Furthermore ${\cal H}^0(\delta_{W_r}*\delta_{W_s})={\cal H}^{-1}(\delta_{W_r}*\delta_{W_s})
= 0 $ for all $\vert r-s\vert \neq 1,0$ and ${\cal H}^0(\delta_{W_r}*\delta_{W_s})=0$ and
${\cal H}^{-1}(\delta_{W_r}*\delta_{W_s}) = \delta_C[-1]$ for $\vert r-s\vert =1$. }
\end{Lemma}

\bigskip\noindent
\underbar{Proof}: The first two statements are special cases of the computations made in
section \ref{help} implicitly using, that the singularities of $W_d$ have codimension $\geq 2$
\cite{ACGH}, and that $W_d + x = W_d$ implies $x=0$ (one easily reduces this to the case $d=1$
by the use of the operation $\ominus$). The second assertion uses theorem \ref{thm1} together
with $IH^1(W_d) = H^1(C)$ (see the remark at the end of section \ref{Cohomology groups}). The
remaining assertions follow from the spectral sequence exploited in section \ref{help} and from
${\cal H}^{-\nu}(\delta_{W_{d}})=0$ for $\nu\neq -d$ and ${\cal
H}^{-d}(\delta_{W_{d}})=\overline\Q_{l,W_d}$. We skip the details.

\bigskip\noindent
\goodbreak

\chapter{Basic results}\label{Sbr}

\section{Riemann-Roch}\label{RR}
There exists a canonical point $\kappa\in X$, which depends on the choice of the point $P_0$
defining the Albanese morphism $C\to X$,  such that $ \kappa - W_r^\nu = W_s^\mu $ holds for
all $r+s=2g-2$ and $\mu=g-1-r+\nu$. This is a reformulation of the Riemann-Roch theorem (see
\cite{Gunning}, p. 49). In other words
$$ \kappa - W^\nu_{g-1-\tau} = W^{\tau +\nu}_{g-1+\tau} $$
holds for $0\leq \tau\leq g-1$. Since $ W^1_{g-1+\tau} = \cdots = W^\tau_{g-1+\tau}=X $ for all
$1\leq \tau\leq g-1$, we may assume $\mu=\tau +\nu$ and $1\leq \nu$ for the study of
$W_{g-1+\tau}^\mu$ in the range $1\leq \tau\leq g-1$.

\bigskip\noindent
\underbar{Notation}: For $\varphi_\kappa(x)=\kappa-x$ and a complex $K\in
D_c^b(X,\overline\Q_l)$ let $K^-$ denote the complex $\varphi_\kappa^*(K)$.

\bigskip\noindent
\begin{Lemma} \label{RR} {\it For $d=g-1+\tau$ and $1\leq \tau \leq g-1$ the following holds
$$ \delta_{g-1+\tau} \ =\ \delta_{g-1-\tau}^-(-\tau) \  \oplus\
\bigoplus_{\nu=0 }^{\tau-1} \delta_X[\tau-1-2\nu](-\nu) \ $$ and
$$ \delta_{g-1+\tau}^- \ =\ \delta_{g-1-\tau}(-\tau) \  \oplus\
\bigoplus_{\nu=0 }^{\tau-1} \delta_X[\tau-1-2\nu](-\nu) \ .$$} \end{Lemma}

\bigskip\noindent
\underbar{Proof}: Consider the closed subscheme
$Y=W^{\tau+1}_{g-1+\tau} \hookrightarrow X$ and its open
complement $j:U=X\setminus Y \hookrightarrow X$. Let $E$ be a
sheaf of $C^{(g-1+\tau)}$ (not a sheaf complex).
%shifted to degree $-(g-1+\tau)$.
The map $p:C^{(g-1+\tau)}\to X$ is an smooth morphism over $U$. Its fibers over $U$ are
projective spaces of dimension $\tau -1$. Therefore the standard truncation of the direct image
complex $L=Rp_*(E[g-1+\tau])$ on $U$ vanishes in degree $\geq \tau - g +1$.Hence the standard
truncation $\tau^{st}_{\geq \tau-g+1}L$ is zero on $U$. The distinguished triangle
$$ \tau^{st}_{\leq \tau -g}L \to Rp_*(L) \to \tau_{\geq \tau-g+1}^{st}
L \to $$ implies $ \tau^{st}_{\leq \tau-g +1}L\big\vert_U \ \cong\
L\big\vert_U $ and $K=\tau_{\geq \tau-g+1}^{st} L$ is a complex
with cohomology supported in $Y$.

\bigskip\noindent
The cohomology stalks of $K$ at points $y\in Y^\alpha\setminus Y^{\alpha +1}$ for $Y^\alpha=
W^{\tau+\alpha}_{g-1+\tau}\subseteq Y$ are
$$ {\cal H}^\nu(K)_y = \bigoplus_{\nu\geq 2\tau}
H^\nu\Bigl(\Proj^{\tau-1+\alpha},E[g-1+\tau]\big\vert \Proj^{\tau-1+\alpha}\Bigr) \ .$$ These
cohomology stalks vanish on $Y^\alpha$ in degrees
$$\nu > 2(\tau-1+\alpha) -(g-1+\tau)= \tau - g +1 +2(\alpha-1)\ .$$
 By Riemann-Roch $Y=\kappa- W_{g-1-\tau}$ and
$Y^{\alpha}=\kappa-W^{\alpha}_{g-1-\tau}$.
%Hence $dim(Y^\alpha,Y)=codim(W^{\alpha}_{g-1-\tau},W_{g-1-\tau})$.
Hence $codim(Y^\alpha,Y) = codim(W^{\alpha}_{g-1-\tau},W_{g-1-\tau})$. By Martens theorem  this
is $\geq 2(\alpha -1)$ for $g\geq 3$. Therefore $-dim(Y^\alpha) \geq -dim(Y) + 2(\alpha -1) =
\tau - g + 1 + 2(\alpha-1) $. We have shown ${\cal H}^\nu(K)\big\vert_{Y^\alpha \setminus
Y^{\alpha+1}}=0$ for $\nu
> -dim(Y^\alpha)$. Hence the complex $K$ is semi-perverse
$$ K \in {}^p D^{\leq 0}(Y) \subseteq {}^p D^{\leq 0}(X) \ .$$
Furthermore, if $C$ is not hyperelliptic, the restriction of $K$ to $Y^\alpha$ for $\alpha\geq
2$ is contained in ${}^p D^{\leq -1}(Y^\alpha)$ by Martens theorem.

\bigskip\noindent
Now assume $E$ to be a selfdual pure perverse sheaf of weight $w$ on $C^{(g-1+\tau)}$. Since
$p$ is a proper morphism, $L=Rp_*(E)$ again is a selfdual pure complex on $X$. Hence $ L =
\bigoplus {}^p H^\nu L\langle \nu\rangle$ by Gabber's theorem. Then $M^\nu = j_{!*} ({}^p H^\nu
L)$ and we get
$$ L = \tilde L \ \oplus \ \bigoplus M^\nu \langle \nu\rangle \ ,$$
where $\tilde L$ again is a selfdual pure complex, whose cohomology support is contained in
$Y$.  Notice $\tau^{st}_{\tau-g+1}L \vert_Y \in {}^p D^{\leq 0}(Y)$. The same applies for any
direct summand, hence $\tau^{st}_{\geq \tau-g+1}\tilde L \in {}^p D^{\leq 0}(Y)$. Since $\tilde
L$ is selfdual and pure, $\tilde L$ is a direct sum of translates of perverse sheaves. If
$\tilde L$ is not a perverse sheaf on $Y$, being selfdual, ${}^p H^\mu(\tilde L)\neq 0$ holds
for some $\mu>0$. However the corresponding translated summand of $\tilde L$ remains unaffected
by the truncation $\tau^{st}_{\geq \tau - g+1}=\tau^{st}_{\geq-dim(Y)}$ by [KW] III.5.13. This
contradicts $\tau^{st}_{\geq \tau-g+1}\tilde L \in {}^p D^{\leq 0}(Y)$ unless $\tilde L$ is in
$Perv(Y)$. Hence
$$ \tilde L \in Perv(Y) \ .$$
Concerning the sheaves $M^\nu$: Let us now assume that $E$ was the constant sheaf on
$C^{(g-1+\tau)}$. We apply [FK] III.11.3 to the projective smooth morphism $p^{-1}(U)\to U$.
Therefore ${}^p H^{\tau -1}(Rp_*\delta_{C^{(r)}}\vert U)(\tau - 1)$ is the constant perverse
sheaf on $\delta_U\in Perv(U)$. By the hard Lefschetz theorem applied for the projective bundle
\footnote{See \cite{Gunning} theorem 16(b)} morphism $p^{-1}(U)\to U$, then all perverse
cohomology sheaves ${}^p H^\nu(Rp_*\delta_{C^{(r)}})$ are constant perverse sheaves up to
shifts and Tate twists. Therefore
$$ \bigoplus M^\nu = \bigoplus_{\nu=0}^{\tau -1}
\ \delta_X[\tau-1-2\nu](-\nu) \ ,$$ where the term for $\nu=\tau-1$ (restricted to $U$)
corresponds to $H^{\tau -1}(Rp_*\delta_{C^{(r)}}\vert U)$. Since $\tau^{st}_{\geq \tau -
g+1}(\bigoplus M^\nu) =0$ in this case, we see that $\tilde L = \tau^{st}_{\geq \tau - g+1} L
=K$.

\bigskip\noindent
Finally we claim, that for constant $E$ the perverse sheaf $\tilde L$ on $Y$ is the
intermediate extension of the constant perverse sheaf on $Y$, if $C$ is not hyperelliptic. For
this remains to show, that $\tilde L\big\vert_{Y\setminus Y^2}$ is constant up to a Tate twist.
However, over $Y\setminus Y^2$ we can argue as above to show that ${}^p H^\tau(\tilde L) = {}^p
H^\tau(K)$ is constant up to a $\tau$-fold Tate twist. Therefore $$ \tilde L \ =\
\delta_{g-1-\tau}^-(-\tau) \ .$$ This remains true for hyperelliptic $C$ with a similar
argument, although the perverse sheaves are not irreducible any longer.

\bigskip\noindent
By a similar argument we get

\bigskip\noindent
\begin{Lemma} {\it For $d\geq 2g-1$ we get $\delta_d = \bigoplus_{\nu=0 }^{d-g}\ \delta_X[d-g-2\nu](-\nu)$.}
\end{Lemma}

\bigskip\noindent
In particular the class of $\delta_d$ in $K_*(X)$ is zero for $d >
\chi=2g-2$.

\bigskip\noindent
\section{Triviality}

\bigskip\noindent
\begin{Lemma}\label{lastt} {\it For a partition $\gamma=(\gamma_1,\gamma_2,...)$ with $\gamma_1> \chi$ the class of $\delta_\gamma$ in
$K_*(X)$ is zero.} \end{Lemma}

\bigskip\noindent
\underline{Proof}: Notice $m_{\alpha \beta}^\gamma=1$ for $\alpha=(\gamma_1,0,..)$ and
$\beta=(\gamma_2,\gamma_3,...)$. Hence $\delta_\gamma$ is a direct summand of the complex
$\delta_\alpha * \delta_\beta$. Since $\delta_\alpha$ is a direct sum of translates of the
perverse sheaf $\delta_X$ (up to Tate twists) for $\gamma_1> \chi$, the same holds true for the
convolution $\delta_\alpha * \delta_\beta$, hence also for each direct summand.

\bigskip\noindent
\begin{Lemma} {\it For $\alpha=(2g-2,..,2g-2,0,..)$ of degree $deg(\alpha)= r(2g-2)$ the class of
$\delta_\alpha$ in $K_*(X)$ is equal to the class of the perverse skyscraper sheaf
$\delta_{\{r\cdot \kappa\}}(r\chi)$ concentrated in the point $r\cdot \kappa \in X$.}
\end{Lemma}

\bigskip\noindent
\underbar{Proof}: $(\delta_{2g-2})^{*r} = \delta_{2g-2}* ...
* \delta_{2g-2} = \delta_\alpha \oplus \bigoplus_\beta
\delta_\beta$ for $\beta=(\beta_1,..)$ with $\beta_1>2g-2$. $(\delta_{2g-2})^{r} \equiv
\delta_\alpha$ in $K_*(X)$ by lemma \ref{lastt}. But $\delta_{2g-2} \equiv
\delta_{\{\kappa\}}(\chi)$ by lemma \ref{RR}. This proves the claim.

\section{Duality}

\bigskip\noindent
%We start with some remarks on skyscraper sheaves. Convolution $K*\delta_{x_0}$ of an arbitrary
%complex $K$ with a skyscraper sheaf $\delta_{x_0}$ translates $K$.  Sometimes it will annoying,
In the following iterated convolutions with skyscraper sheaves at $r\cdot \kappa\in X$ occur.
To avoid them, one could consider $X$-equivariant perverse sheaves on $X\times X$ instead of
perverse sheaves on $X$ . Notice, there is a natural equivalence between perverse sheaves on
$X$ up to translation in $X$, and $X$-equivariant perverse sheaves on $X\times X$. Convolution
products of $X$-equivariant perverse sheaves on $X\times X$ can be described by the formalism
of [KW], bottom of p.192. In the following, instead of using equivariant perverse sheaves, we
often prefer to normalize $\kappa$ to become $0$ using an appropriate translate of $C$ in $X$
instead of the curve $C$ itself.

\bigskip\noindent
\underbar{Notation}: For $\alpha=(\alpha_1,...,\alpha_r,0,..,0)$ suppose $\alpha_1\leq
\chi=2g-2$. Then put
$$\beta=(\chi-\alpha_r,\chi-\alpha_{r-1},...,\chi-\alpha_1,0,..,0)\ .$$
$\delta_\alpha$ and $\delta_\beta$ are pure perverse sheaves of weight $deg(\alpha)$
respectively $deg(\beta)= r\chi -deg(\alpha) $. Set $\tau=r(g-1) - deg(\alpha)$. Then

\bigskip\noindent

\bigskip\noindent
\begin{Theorem} \label{dual}{\it $(\delta_\alpha)^- \equiv T^*_{r\cdot \kappa}(\delta_\beta)(-\tau)$ holds as an equality in
$K_*(X)$.} \end{Theorem}

\bigskip\noindent
\underbar{Remark}: If $\kappa$ is normalized to be zero, then ignoring Tate twists this can
also be written
$$ \fbox{$ \delta_\alpha^\vee \equiv \delta_\beta $} \ ,$$
since $D(\delta_\alpha)\cong \delta_\alpha$.

\bigskip\noindent

\bigskip\noindent
\underbar{Proof}: For the proof we ignore the Tate twists and also assume $\kappa=0$. We use
induction and the two formulas: \ \
$$ \mbox{\it Formula 1)} \quad \quad \delta_{(\alpha_1,..,\alpha_{r-1})} * \delta_{\alpha_r}
\ = \ \delta_{\alpha}\ \oplus\ \bigoplus_{\nu} m_+(\nu) \delta_\nu \ $$ with summation  over
all $\nu=(\nu_1,..,\nu_r,0,..,0)$ with $\nu_r < \alpha_r$ where we write
$$m_+(\nu):=m_{(\alpha_1,..,\alpha_{r-1})(\alpha_r)}^{(\nu_1,...,\nu_r)} \in \N$$ for the multiplicity. Similarly
$$ \delta_{(\chi-\alpha_{r-1},..,\chi-\alpha_{1})} * \delta_{\chi - \alpha_r}
\ = \ \delta_{\beta}\ \oplus\ \bigoplus_{\mu} m(\mu) \delta_\mu \ ,$$ where the sum now only
runs over the $\mu=(\mu_1,..,\mu_r,0,..,0)$ with $\mu_1 > \chi - \alpha_r$. Again
$m(\mu)=m_{(\chi-\alpha_{r-1},..,\chi-\alpha_{1})(\chi - \alpha_r)}^\mu \in \N$ denotes the
multiplicity. Now we read this formula in $K_*(X)$. This allows to ignore all $\mu$ with $\mu_1
> \chi=2(g-1)$. In particular $(\mu_1,...,\mu_r)=(\chi-\nu_r,...,\chi-\nu_1)$ is then well defined.
Thus, as a formula that holds in $K_*(X)$, this gives: \ \
$$ \mbox{\it Formula 2)} \quad \quad \delta_{(\chi-\alpha_{r-1},..,\chi-\alpha_{1})} * \delta_{\chi - \alpha_r}
\ \equiv \ \delta_{\beta}\ \oplus\ \bigoplus_{\nu} m_-(\nu)
\delta_{(\chi-\nu_r,...,\chi-\nu_1)} \ ,$$ where summation runs over all
$\nu=(\nu_1,..,\nu_r,0,..,0)$ for which $\chi\geq \chi- \nu_r > \chi - \alpha_r$ holds, and
where
$$m_-(\nu) := m_{(\chi-\alpha_{r-1},..,\chi-\alpha_{1})(\chi -
\alpha_r)}^{(\chi-\nu_r,...,\chi-\nu_1)}\ .$$

\bigskip\noindent
For the proof of the duality statement between $\delta_\alpha$ and $\delta_\beta$ we now use an
iterated induction on $r$ and on $\alpha_r$. By the induction assumption on $r$ we show that
the left side of formula 1) and formula 2) correspond under duality. Then, by induction on $r$
and $\alpha_r$ we show, that all terms except $\delta_\alpha$ respectively $\delta_\beta$ on
the right´side of formula 1) respectively formula 2) again correspond under duality. This
implies the desired duality between $\delta_\alpha$ and $\delta_\beta$ by comparison.

\bigskip\noindent
\underbar{Equality of the left hand sides}:  For $-\tau_1 = \alpha_r - (g-1)$ we have by
induction on $r$
$$\delta_{\alpha_r}^- =
T_\kappa^*(\delta_{\chi-\alpha_r})(-\tau_1)\ ,$$  and similarly for $-\tau_{r-1}=
\alpha_1+..+\alpha_{r-1} - (r-1)(g-1)$
$$\delta_{(\alpha_1,..,\alpha_{r-1})}^- =
T_{(r-1)\kappa}^*(\delta_{(\chi-\alpha_{r-1},..,\chi-\alpha_1)})(-\tau_{r-1})\ .$$ Here we also
used $T_{x_0}^*(K)
* T_{x_1}^*(L) = T_{x_0+x_1}^*(K*L)$.

\bigskip\noindent
\underbar{Equality of the sums on the right side}: For a comparison of terms in the direct sums
on the right sides of formula 1) and 2) by induction on $r$ and on  $\alpha_r$ we have
$\delta_{(\nu_1,..,\nu_r)}^- \equiv T_{r\kappa}^*(\delta_{\chi-\nu_r,..,\chi-\nu_1)}$ (as an
equality in $K_*(X)$).

\bigskip\noindent
This being said, we remark that the Littlewood-Richardson coefficients
$$m_+(\nu)=m_-(\nu)$$ coincide. This will be shown in the next section. Taking this for granted
now this implies, that the sums on the right sides of formula 1) and 2) agree with
multiplicities. Hence $\delta_\alpha^- \equiv T^*_{r\kappa}(\delta_\beta)-(\tau)$ in $K_*(X)$,
as desired.

\bigskip\noindent
\underbar{Induction start}: Concerning the induction start
$\alpha_r=0$ we remark, that for $\alpha_r=0$ we are in principle
in the case of the induction used for the left sides. The only
subtlety is, that we now have to use, that
$\delta_{(\chi-\alpha_r,...,\chi-\alpha_1)} =\delta_{(
\chi,\chi-\alpha_{r-1}...,\chi-\alpha_1)} \equiv \delta_\chi *
\delta_{(\chi-\alpha_{r-1},...,\chi-\alpha_1)} \equiv
\delta_{(\chi-\alpha_{r-1},...,\chi-\alpha_1)}$ holds in $K_*(X)$.
In fact, this turns out to be a special case of the formulas
stated above.

%For instance in case $r=s$ we get for $1\leq r\leq g-1$
%$$ Rf_* (\delta_{W_r\times W_r}) \ = \
%Rf_* (\delta_{W_{r-1}\times W_{r-1}})

\section{The Littlewood-Richardson rule}
\label{LRrule}

A partition $\gamma=(\gamma_1,\gamma_2,\cdots)$ has an associated  Young diagram. Consider it
as lying in the first plane-quadrant touching the x-axis. The first column - touching the
y-axis - has height $\gamma_1$, the second column to the right has height $\gamma_2$, and so
on.

\bigskip\noindent
To describe the rule Littlewood-Richardson notice, that $m_{\alpha\beta}^\gamma$ vanishes
unless $\alpha_i\leq \gamma_i$ holds for all $i$. In other words: the Young diagram of $\alpha$
has to lie inside the Young diagram of $\gamma$. This defines the complementary skew tableau
$\gamma/\alpha$, which is the complement of the Young diagram of $\alpha$ in the Young diagram
of $\gamma$ (it is a subdiagram of the Young diagram of $\gamma$). Fill the entries of the skew
tableau with numbers $1,2,3,..$,  having the multiplicities $\beta_1,\beta_2,..$ respectively,
according to the filling rules 1-3) described below.

\bigskip\noindent
Then the Littlewood-Richardson rule says, that the multiplicity $m_{\alpha\beta}^\gamma$ is the
number of different possible fillings of the skew diagram $\gamma/\alpha$ satisfying the

\goodbreak
\bigskip\noindent
\underbar{Filling rules}: \begin{enumerate} \item[1)]  Entries in each column must be weakly
increasing in the direction away from the $x$-axis. \item[2)] Entries on the same horizontal
line of the skew diagram must be  strictly increasing in the direction away from the $y$-axis.
\item[3)] Finally, starting from the highest entry of the first column of the skew diagram,
then descending in direction of the $x$-axis, then starting from the highest entry in the
second column of the skew diagram descending in direction of the $x$-axis, and so on, each
filling defines a sequence of integers. Say $n_\nu$ counts, how often the integer $n$ occurs
among the first $\nu$ terms of the sequence. Then the filling rule 3) requires, that $n\leq m$
always has to imply $n_\mu \geq m_\mu$ (for all $\mu\geq 1$). In other words: The integer $n$
occurs not less often then $m$ among the first $\mu$ terms (no matter how $\mu$ is chosen).

\end{enumerate}

\bigskip\noindent
\underbar{1. example}:  $\alpha=(\alpha_1,...,\alpha_{r-1},0,...\ )$ and
$\beta=(\alpha_r,0,...\ )$

\bigskip\noindent
In this case $\gamma/\alpha$ is filled with the single number 1. Hence rules 1) and 3) become
meaningless and only rule 2) has a consequence: in each horizontal line of the skew diagram at
most one entry is possible. Hence in our example the total number of filling positions is
$\alpha_r$. Hence the possible diagrams $\gamma$, which allow a filling, have the form (*)
$$ \gamma_1 \geq \alpha_1 \geq \gamma_2 \geq ...  \geq
\gamma_{r-1} \geq \alpha_{r-1} \geq \gamma_r \geq 0 \ .$$  In particular, rule 2) excludes the
nonvanishing of $\gamma_{r+1}$. Every $\gamma=(\gamma_1,...,\gamma_r)$ satisfying condition (*)
and the following trivial condition (**) $$deg(\gamma) =\sum_{i=1}^r \alpha_i$$ (this is the
condition, that the total number of fillings is $\alpha_r$) can occur, and it occurs with
multiplicity
$$m_+(\gamma)=1\ .$$

\bigskip\noindent
\underbar{Additional remark}: For $\alpha_{r-1}\geq \alpha_r$, the lexicographic minimal
possible choice for such $\gamma$ gives $\gamma=(\alpha_1,..,\alpha_r)$. All the other
$\gamma$'s, which allow fillings as above, necessarily satisfy
$$ \alpha_r > \gamma_r  \geq 0 \ $$
where $\gamma_r = \alpha_r - \sum_{i=1}^{r-1} (\gamma_i -
\alpha_i)$ is uniquely determined by $\gamma_1,..,\gamma_{r-1}$.

\goodbreak
\bigskip\noindent
\underbar{2. example}: $\alpha=(\chi-\alpha_{r-1},...,\chi-\alpha_1,0,...\ )$ and $\beta=(\chi
- \alpha_r,0,...\ )$.

\bigskip\noindent
Again filling is done with the single number 1 with only rule 2) being relevant. Possible
$\gamma$, which allow a filling,  are now of the form
$$ \gamma_1 \geq \chi-\alpha_{r-1} \geq \gamma_2 \geq \chi -\alpha_{r-2} \geq ...  \geq
\gamma_{r-1} \geq \chi -\alpha_{1} \geq \gamma_r \geq 0 \ .$$ Again rule 2) excludes
$\gamma_{r+1}>0$. Again there exists a unique filling, if this condition together with the
degree condition $deg(\gamma) = r\chi -\sum_{i=1}^r \alpha_i$ holds. Hence the multiplicity is
$$ m_-(\gamma)=1
$$ in this case, and is zero otherwise.

\bigskip\noindent
\underbar{For the proof of duality}: Comparing the first and the second example we complete the
proof of theorem \ref{dual}. For this parameterize the possible $\gamma=(\gamma_1,..,\gamma_r)$
by their dual parameters $\chi-\check\gamma = (\chi-\check\gamma_r,...,\chi-\check\gamma_1)$.
This is possible without  problem for the parameters $\gamma_i$ if $i=2,..,r$. This gives a
true partition $\check \gamma$ only for $\gamma_1\leq \chi$. However we will see immediately,
that we can ignore the other cases by lemma \ref{lastt}. Ignoring these cases, the conditions
above  for the dual parameter become
$$ \chi -\check\gamma_r \geq \chi-\alpha_{r-1} \geq \chi-\check\gamma_{r-1} \geq \chi -\alpha_{r-2} \geq ...  \geq
\chi -\check\gamma_{2} \geq \chi -\alpha_{1} \geq \chi-\check\gamma_1 \geq 0 \ .$$ Equivalently
this means (*)'
$$ \check\gamma_r \leq \alpha_{r-1} \leq \check\gamma_{r-1} \leq
\alpha_{r-2} \leq ...  \leq \check\gamma_{2} \leq \alpha_{1} \leq \check\gamma_1 \leq \chi \
.$$ The additional trivial degree condition (counting the total number of fillings) is $[(\chi
-\check \gamma_r) - (\chi-\alpha_{r-1})] + ... + [(\chi - \check \gamma_2) - (\chi - \alpha_1)]
+ (\chi -\check \gamma_1) = \chi - \alpha_r$ or (**)'
$$ deg(\check \gamma) = \sum_{i=1}^r \alpha_i \ .$$

\bigskip\noindent
\underbar{Additional remark}: Suppose $\alpha_1\geq \alpha_2$. Then
$\gamma_i=\chi-\alpha_{r-i}$ for $i=1,..,r$ defines the lexicographical minimal choice for
$\gamma$, subject to the conditions (*)' and (**)' above. In the other
 cases $ \gamma_1 = \chi - \check \gamma_r > \chi -\alpha_1 $, i.e. $$ \alpha_1 > \check \gamma_r \ .$$

\bigskip\noindent
\underbar{Comparison}: The discussion above shows, that the solutions $\gamma$ of the
conditions (*) and (**) of the first example precisely match with the solutions $\check \gamma$
of the conditions (*)' and (**)' of the second example, except for the \lq{missing}\rq\
boundary condition
$$ 0 \leq \check \gamma_r $$
in the second example (which is equivalent to $\gamma \leq \chi$), respectively the
\lq{missing}\rq\ boundary condition in the first example
$$ \gamma_1 \leq \chi \ .$$
If we discard such values of $\gamma_1$ (which is possible in the setting of the proof of the
duality theorem in the last section) this establishes a bijection between solutions defined by
$$ \gamma \longleftrightarrow \chi - \gamma \ $$
This duality was already used for the proof of the duality theorem together with the fact, that
in these special cases all multiplicities are one.

\bigskip\noindent

\section{Perverse depth}\label{PD}

\bigskip\noindent
For all $\alpha$ the complex $\delta_\alpha$ is a direct sum of a perverse sheaf and a direct
sum $\bigoplus T_\nu[\nu]$ of translation-invariant perverse sheaves $T_\nu$ by theorem
\ref{Geomorigin}. In this section we prove a slightly stronger result. The proof uses (an
elementary special case) of the Littlewood-Richardson rule and affine vanishing theorems (see
the following lemma \ref{post}).

\bigskip\noindent We abbreviate $\varepsilon_r =
\delta_{\alpha}$ for $\alpha=([1]^r)=(1,...,1,0,..,)$. For the moment we assume the following
lemma later proved in section \ref{PoA}

\bigskip\noindent
\begin{Lemma}\label{post} For all $r$ there exists a constant $C_0$ such that $${}^p
H^\nu(\varepsilon_r)=0 \quad  \mbox{ holds for all } \quad \nu\notin [-C_0,C_0]\ .$$
\end{Lemma}

\bigskip\noindent
Using this lemma we claim a similar vanishing statement holds for all $\delta_\alpha$ which
satisfy $\alpha_1\leq \chi$ for some new constant $C_1$ instead of $C_0$:

\bigskip\noindent
To show this notice, that for $\alpha_1\leq \chi$ the complex $\delta_\alpha$ is a direct
summand of the complex
$$ \varepsilon_{\beta_1} * \cdots * \varepsilon_{\beta_s} \ ,$$
where $\beta=(\beta_1,...,\beta_s)$ is defined by $$\beta = \alpha^*\ .$$  This statement is an
elementary special case of the Littlewood-Richardson rule corresponding to the fact that every
irreducible representation of the linear group is contained in a tensor product of the
fundamental representations. Notice $s\leq \alpha_1\leq \chi$. Since the perverse depth of the
complex $\varepsilon_{\beta_1} \boxtimes \cdots \boxtimes \varepsilon_{\beta_s} \in D_c^b(X^s)$
is at most $s\cdot C_0 \leq \chi\cdot C_0$, and since the direct image under the iterated
smooth addition map $ a:X^s \longrightarrow X $ shifts the perverse depth at most by the
relative dimension $g\cdot (s-1)$, we obtain the estimate $s\cdot C_0 + g\cdot (s-1) \leq
\chi\cdot (C_0+g)$ for the perverse depth of $\varepsilon_{\beta_1}
* \cdots * \varepsilon_s$, hence the same estimate for the perverse depth of any $\delta_\alpha$.
This proves the claim.

\bigskip\noindent
By triviality $\delta_\alpha \equiv 0 $, which holds for all $\alpha$ with $\alpha_1 >\chi$, we
therefore get the following weaker result (now for all $\alpha$)

\bigskip\noindent
\underbar{Conclusion}: {\it There exists a constant $d$ such that for all partitions $\alpha$
the perverse cohomology sheaves ${}^p H^\nu(\delta_\alpha)$ are constant perverse sheaves for
$\nu\notin [-d,d]$.}

\bigskip\noindent
Choose $d$ minimal with this property.%, i.e  choose $\alpha$ with $\alpha_1\leq \chi$ such that
%$$ {}^p H^d(\delta_\alpha) \not\equiv 0 \mbox { in } K_*(X)\quad \mbox{ and }$$
%$$ {}^p H^\nu(\delta_\beta) \equiv 0 \mbox{ in } K_*(X) \mbox{ for all }
%\nu > d \mbox{ and } \beta \mbox{ with } \beta_1\leq \chi\ . $$
We claim $d=0$. Suppose $d>0$. Then there exists $\delta_\alpha$ such that
${}^pH^d(\delta_\alpha)$ is not constant. For certain irreducible perverse sheaves $A=A_i,i\in
I$ we write
$$ \delta_\alpha = \bigoplus_A \bigoplus_{\nu=-\nu_{A}}^{\nu_A} A[\nu] \quad , \quad (\delta_\alpha)^\vee = \bigoplus_A \bigoplus_{\nu=-\nu_A}^{\nu_A} A^\vee[\nu]  $$
using the Hard Lefschetz theorem and Gabber's theorem. For some nonconstant $A$ then $\nu_A =
\nu_{A^-}=\nu_{DA} = d$. Let $J \subseteq I$ parameterizing all $A_j$ with the same
property. %By duality $\delta_\alpha^\vee \equiv T^*_y(\delta_\beta)$ holds for some $y=r\cdot
%\kappa\in X(k)$.
Then
$$ \delta_\alpha * \delta_\alpha^\vee \ =\ \ \bigoplus_{i\in I}\
\bigoplus_{j\in I}\ \bigoplus_{\nu=-\nu_{A_i}}^{\nu_{A_i}}
\bigoplus_{\mu=-\nu_{A_j}}^{\nu_{A_j}} A[\nu] * A^\vee[\mu]  =\ \bigoplus_{i,j\in J} (A_i *
A_j^\vee)[-2d] \ + \ \mbox{lower shifts} \ .$$ By lemma \ref{van}
$$ {\cal H}^{2d}(\delta_\alpha * \delta_\alpha^\vee) = \bigoplus_{i\in J} {\cal
H}^0(A_i
* A_i^\vee) \neq 0 \ .$$
Duality $\delta_\alpha^\vee \equiv T^*_y(\delta_\beta)$ for some $y=r\cdot \kappa\in X(k)$
implies
$$ {\cal H}^{2d}(\delta_\alpha * \delta_\alpha^\vee) =
{\cal H}^{2d}( \delta_X-\mbox{translates}) \ \oplus \ T^*_y{\cal
H}^{2d}(\delta_\alpha*\delta_\beta) $$ $$ = \ {\cal H}^{2d}( \delta_X-\mbox{translates}) \
\oplus \ \bigoplus_{\gamma} m_{\alpha\beta}^\gamma \cdot T^*_y{\cal H}^{2d}(\delta_\gamma) \
.$$ Since all ${}^p H^{\nu}(\delta_{\gamma})$ are constant perverse sheaves for $\nu\notin
[-d,d]$, all ${\cal H}^{\nu}(\delta_\gamma)$ are constant sheaves for $\nu >d$ and therefore $
{\cal H}^{2d}(\delta_\alpha * \delta_\alpha^\vee) \cong \oplus^l\ \overline\Q_{l,X} $ for some
integer $l>0$. Hence for $l>0$
$$ \bigoplus_{i\in J} {\cal H}^{0}(A_i * A_i^\vee) \ \cong \ \bigoplus^l\ \overline\Q_{l,X}\ .$$ But each support of ${\cal H}^{0}(A_i * A_i^\vee)$ is contained in $S(A_i)=\{x\in
X\ \vert \ T^*_x(A_i)\cong A_i\}$. From corollary \ref{cor1} we conclude:  $T^*_x(A_i)\cong
A_i$ holds for all $x\in X$, if $j\in J$. This forces $A_i\cong \delta_{\chi_i}$ for some
unramified character $\chi_i\in Hom(\pi_1(X,0),\overline \Q_l^*)$ by section \ref{TI}. The next
lemma implies that all $\chi_i$ are trivial. This gives a contradiction proving $d=0$.

\begin{Lemma}\label{notra}
Suppose $\delta_\chi \hookrightarrow \delta_\alpha$ for some partition $\alpha$ and some
character $\chi$, then $\chi$ is the trivial character.
\end{Lemma}

\bigskip\noindent
\underbar{Proof}: If follows $\delta_\chi^{*r} \hookrightarrow \delta_\alpha^{*r}$. Since $
\delta_\chi^{*r} = t^{(r-1)g/2} \delta_\chi $ plus terms of lower order in the Grothendieck
group, the conclusion above implies that $\delta_\chi$ is constant by choosing $r$ large enough
so that $(r-1)g/2
> d$ holds. Hence $\chi$ must be constant.

\bigskip\noindent
We conclude

\bigskip\noindent
\begin{Theorem}\label{thm8} {\it For $\nu\neq 0$ and arbitrary $\alpha$ there exists an integer $m(\alpha,\nu)$ such that
$$ {}^p H^\nu(\delta_\alpha) \cong
 m(\alpha,\nu)\cdot \delta_X \quad \mbox{ for } \quad \nu\neq 0\ .$$
Furthermore ${}^p H^0(\delta_\alpha)$ does not contain irreducible constituents $\delta_\chi$
except for the trivial character $\chi$.} \end{Theorem}

\bigskip\noindent
\begin{Corollary} \label{Zerlegung}{\it There exist pure perverse sheaves ${}^p\delta_\alpha$,  which do not contain any
translation-invariant constituent, such that
$$ \delta_\alpha = {}^p \delta_\alpha \ +\  P_\alpha(t)\cdot \delta_X \quad , \quad P_\alpha(t) \in \Z[t^{1/2},t^{-1/2}]$$
holds in the Grothendieck ring $K^0(Perv_m(X))\otimes \Z[t^{1/2},t^{-1/2}]$.} \end{Corollary}

\bigskip\noindent
\underbar{Example}: In case $\alpha_1=g$ for instance $$ \delta_\alpha = {}^p\delta_\alpha \
\oplus \ dim(H_-^{(\alpha_2,...,\alpha_r,..)^*})\cdot \delta_X \ .$$ The multiplicity
$dim(H_-^{(\alpha_2,...,\alpha_r,..)^*})$ is a nonzero constant!

\section{Computation of $P_\alpha(t)$} \label{Compu}
For the perverse sheaves ${}^p\delta_\alpha$ the hypercohomology groups $IH^\nu(X,{}^p
\delta_\alpha)$ vanish in degrees $\nu \geq g$, since they do not contain constant summands.
Since the degree of $h(\delta_X,t)$ in $t^{1/2}$ is $g$, $P_\alpha(t)$ can be computed from
$h(\delta_\alpha,t)$ via
$$ \frac{h(\delta_\alpha,t)}{h(\delta_X,t)} = P_{\alpha}(t)
+ \frac{r(\alpha,t)}{h(\delta_X,t)} \quad , \quad
deg_{t^{1/2}}(r(\alpha,t)) < g \ .$$

\bigskip\noindent
\begin{Corollary}\label{estimate} {\it $\delta_\alpha = {}^p \delta_\alpha$ if $\alpha_1\leq
g-1$. More generally
$$deg_{t^{1/2}} \ (P_\alpha(t)) \ \leq\ \alpha_1-g\ .$$} \end{Corollary}

\bigskip\noindent
\underbar{Proof}: Since $h(\delta_X,t) = t^{-g/2}(1+t^{1/2})^{2g} = (t^{1/4} + t^{-1/4})^{2g}$
has degree $deg_{t^{1/2}}\ h(\delta_X,t) = g$, this follows from corollary \ref{Cohomology
groups} by comparing the dimension of the hypercohomology groups. Notice
$$ deg_{t^{1/2}} \ (h(\delta_\alpha,t)) = \alpha_1\ .$$
This also implies

\bigskip\noindent
\begin{Corollary} {\it For $\alpha_1\leq g-1$ the dimension $d_\alpha$ of the support $Y_\alpha$ of
$\delta_\alpha$ is
$$ d_\alpha \geq \alpha_1 \ .$$} \end{Corollary}

\goodbreak

\section{Proof of Lemma \ref{post}}\label{PoA}

Recall  $ \varepsilon_r = Rp_*({\cal F}_\alpha)$, where $\alpha=(1^r)$ is of degree
$deg(\alpha)=r$ and $p: C^{(r)}\to X$ was the Abel-Jacobi map. ${\cal F}_\alpha$ is a perverse
sheaf on $C^{(r)}$ of the form $ {\cal F}_\alpha = F_\varepsilon[r] $ for an etale $\overline
\Q_l$ sheaf $F_\varepsilon$ on $C^{(r)}$ obtained as a direct factor in $f_*\lambda_{C^r}$ for
the finite morphism $f: C^r \to C^{(r)}$. Its sheaf cohomology is concentrated in degree $-r$.
The morphism $f: C^r\to C^{(r)}$ defining $F_\varepsilon$ is etale on the complement of the
diagonal divisor
$$ \bigcup_{i\neq j} \Delta_{ij} \subseteq C^r \quad  , \quad \Delta_{ij} = \{(x_1,..,x_r)\in C\ \vert \ x_i=x_j\}\
.$$
 The image of $\bigcup_{i\neq j} \Delta_{ij}$ in $C^{(r)}$ defines an irreducible divisor $Y_r$
 of
$C^{(r)}$ equal to the image of the irreducible variety  $C^{(r-2)} \times \Delta(C)$ under the
natural map $$ C^{(r-2)} \times C^{(2)} \overset{\tau}{\longrightarrow} C^{(r)} \ .$$ Let
$U=U_r\overset{j}{\hookrightarrow}\ C^{(r)} $ denote the open complement of $Y=Y_r$
$$ U\ \overset{j}{\hookrightarrow}\ C^{(r)}\ \hookleftarrow\ Y \ .$$
The multiplication law \ref{ML} implies, that ${\cal F}_\varepsilon[r]$ is a direct summand of
$$ \tau_* \Bigl({\cal F}_{\varepsilon}[r-2]\boxtimes {\cal
F}_{\varepsilon}[2]\Bigr) $$ defined by the corresponding sheaves ${\cal F}_{\varepsilon}$ on
$C^{(r-2)}$ respectively $C^{(2)}$. The  stalk of ${\cal F}_\varepsilon$ at a point of $Y$
vanishes. For this it is enough to consider the case $r=2$. For $r=2$ the direct image of the
constant sheaf under the morphism $C\times C \to C^{(2)}$  decomposes into the direct sum of
the constant sheaf and the sheaf $(C^{(2)},{\cal F}_\varepsilon[2])$. The stalks have dimension
1 over the diagonal $\Delta \subseteq C^{(2)}$ and dimension 2 outside the diagonal. Hence the
stalks of $(C^{(2)},{\cal F}_\varepsilon[2])$ have dimension zero over the diagonal and
dimension 1 outside the diagonal.

\bigskip\noindent
Now return to the general case $r\geq 2$. Since $F_\varepsilon$ corresponds to the
sign-character, it becomes a smooth etale sheaf $E=F_{\varepsilon}\vert_{U_r}$ on $U_r$ of rank
one. Its stalks vanish on $Y$. Therefore
$$ F_\varepsilon = j_!(E)\quad , \quad j: U=U_r \hookrightarrow C^{(r)} \ .$$

\bigskip\noindent
There are special interesting cases: First let us consider the
case $r\leq g-1$. In these cases the morphism $p: C^{(r)} \to W_r$
is a birational map. Hence

\bigskip\noindent
\begin{Lemma} {\it For $2\leq r\leq g-1$ the perverse sheaf $\varepsilon_r$ in $Perv(W_r) \subseteq Perv(X)$
is the intermediate extension $$ \varepsilon_r = {\tilde j}_{!*}(E[r]) \quad , \quad {\tilde
j}: U_r \setminus p^{-1}(W_r^2) \hookrightarrow W_r $$ from $U_r\setminus p^{-1}(W^2_r)$ to
$W_r$ for the smooth rank 1 etale sheaf $E=F_\varepsilon\big\vert{(U_r\setminus
p^{-1}(W^2_r))}$.}
\end{Lemma}

\bigskip\noindent
Now consider the case $r>2g-2$. Then $p: C^{(r)} \to X$ is a smooth projective bundle morphism,
whose fibers are projective spaces of dimension $\mu=r-g$. Hence, by the proper basechange
theorem, the cohomology stalk of $\varepsilon_r$ is
$$ {\cal H}^\nu(\varepsilon_r)_x = H^\nu(\Proj^\mu, {\cal F}_\varepsilon)
= H^{\nu+r}_c(\Proj^\mu \cap U_r,F_\varepsilon) \ $$ since ${\cal F}_\varepsilon =
j_!(F_\varepsilon)[r]$. There are two possible cases. Either $\Proj^\mu \cap U_r = \Proj^\mu$
if $x\in X$ is not in the image of the divisor $Y\subseteq C^{(r)}$. Or otherwise $\Proj^\mu
\cap U_r$, which is the complement in $\Proj^\mu$ of the effective divisor $Y\cap \Proj^\mu
\subseteq \Proj^\mu$. Hence, since any effective divisor of the projective space $\Proj^\mu$ is
ample, the complement $\Proj^\mu \cap U_r$ is an affine variety. By the vanishing theorem for
affine varieties by Artin-Grothendieck therefore
$$ H^{\nu+r}_c(\Proj^\mu \cap U_r,F_\varepsilon) = 0 \quad , \quad
\nu+r < \mu = r-g \ .$$ This implies the vanishing of ${\cal
H}^\nu(\varepsilon_r)_x$ for $\nu<-g$.

\bigskip\noindent
\begin{Lemma} {\it For $r> 2g-2$ and $\nu< -g$ the cohomology sheaves ${\cal H}^\nu(\varepsilon_r)$
vanish.} \end{Lemma}

\bigskip\noindent
\underbar{Proof}: By the discussion above it is enough to show, that the composite morphism
$Y\hookrightarrow C^{(r)} \to X$ is surjective. Since $Y$ contains $C^{(r-2)}$, this is clear
for $r-2\geq g$, hence for $g\geq 3$ and for $g=2$ and $r>3$. For  $g=2$ and $r=3$ this is also
true by a direct inspection.

\bigskip\noindent
\begin{Corollary} {\it $\varepsilon_r \in {}^{st}D^{\geq -g}(X)$ for $r> 2g-2$. Hence ${}^p
H^\nu(\varepsilon_r)=0$ for $\nu \geq g$.} \end{Corollary}

\bigskip\noindent
Since $\varepsilon_r$ is a self dual complex on $X$ this implies, that the perverse cohomology
${}^p H^\nu(\varepsilon_r)$ vanishes for $\nu\notin [-g,g]$ provided $r> 2g-2$. Since any
complex $K\in D_c^b(X)$ has finite perverse depth, we conclude that there exists a constant
$C_0$, so that ${}^p H^\nu(\varepsilon_r)$ vanishes for $\nu \notin [-C_0,C_0]$. This proves
lemma \ref{post}.

\bigskip\noindent
\section{BN-sheaves}\label{APP}

\bigskip\noindent
\underbar{Notation}:   A nonconstant irreducible perverse constituent $K$ of some $
\delta_\alpha$, or equivalently an irreducible perverse constituent of some ${}^p
\delta_\alpha$, is called a {Brill-Noether sheaf} ({\bf BN-sheaf}).

\bigskip\noindent
{\bf Definition.} The category ${\cal BN}'$ of BN-sheaves is the full subcategory of $Perv(X)$
generated by direct sums of BN-sheaves on $X$. It is a semisimple abelian category,  closed
under $K\mapsto K^\vee = D(K^-)$. The category is not closed under convolution.

\bigskip\noindent
\begin{Lemma} \label{Ffaith}{\it  The cohomology functor $H^\bullet(X,L)$ is
faithful on the abelian category ${\cal BN}'$.} \end{Lemma}

\bigskip\noindent
\underbar{Proof}: It is enough to show that $H^\bullet(X,L)=0$ implies $L=0$ for $L\in {\cal
BN}$. Let ${\cal N}$ denote the full abelian subcategory generated by all perverse BN-sheaves
$L$, whose cohomology groups $H^\bullet(X,L)=0$. This category ${\cal N}$ is closed under
direct sums,  Verdier duality, duality and the convolution product. By semisimplicity (Gabber's
theorem) it is also closed under taking subquotients! But any subcategory ${\cal N}$ of the
category $Perv(X)$ of perverse sheaves with these properties, which contains a nonzero object,
must contain the unit $1=\delta_0$, since $1 \hookrightarrow L* DL^-$ for $L\neq 0$. Since the
cohomology $H^\bullet(X,\delta_0)=\overline\Q_l$ of $1=\delta_0$ does not vanish, this implies
${\cal N}=0$ and proves the claim.

\bigskip\noindent
Let ${\cal BN}_0'$ be the full subcategory of perverse sheaves in ${\cal BN}'$ of 0-type.
${\cal BN}_0'$ is closed under subquotients, direct sums, convolution and $K\mapsto K^\vee$.

\bigskip\noindent
As aready remarked ${\cal BN}'$ is not closed under convolution. From lemma \ref{mult} and
theorem \ref{thm8} we only have the weaker property, that the convolution of two BN-sheaves is
a direct sum of BN-sheaves and translates of constant perverse sheaves. So, by corollary
\ref{zuvor2}, we may pass from ${\cal BN}'$ to its image ${{\cal BN}}$ in the quotient category
$\overline{Perv}(X)$ (see section \ref{TI}) to obtain

\begin{Lemma} For BN-sheaves $K,L$ we have
$$ Hom_{{\cal BN}}(\overline K,\overline L) \ \cong\ \Gamma_{\{0\}}(X,{\cal H}^0(K*L^\vee)^\vee) \ .$$
\end{Lemma}

\bigskip\noindent
In fact the quotient functor $K\to \overline K$ from $Perv(X)$ to $\overline{Perv}(X)$ induces
an equivalence of ${\cal BN}'$ with its image category ${\cal BN}$. On the quotient category
convolution defines a tensor functor $$ *: {\cal BN} \times {\cal BN} \to {\cal BN} \ .$$

\bigskip\noindent
\begin{Lemma} \label{stable}{\it Let $K$ be a BN-sheaf. Then
\begin{enumerate}
\item $K*\delta_C$ and $K*\delta_{-C}$ are perverse sheaves.
\item If $K$ is irreducible and $\delta_X$ is a summand
of $K*\delta_C$ or $K*\delta_{-C}$, then the dimension of the support of $K$ is $d=g-1$ or
$d=g$.
\item The group of automorphisms of $K$ is finite.
\item $ K\neq 0 \Longrightarrow K*\delta_C\neq 0$.
\item The Verdier dual $D(K)$ is a BN-sheaf. Similarly $K^-$ and $K^\vee= D(K^-)$ are BN-sheaves.
\end{enumerate}} \end{Lemma}

\bigskip\noindent
\underbar{Remark}: These statements carry over for translates $T_x^*(K)$, $x\in X(k)$ and $K$ a
BN-sheaf.

\bigskip\noindent
\underbar{Proof}:  Since $(K*\delta_C)^- = K^-*\delta_{-C}$ the first claim is a special case
of theorem \ref{Geomorigin}. The second is a special case of theorem \ref{thm4}.  By \ref{cor1}
the support of ${\cal H}^0(DK*K^-)$ is a torsor under the automorphism group of $Aut(K)$. Since
$K$ is irreducible perverse nonconstant, by lemma \ref{notra} this support is a proper subset
of $X$. Since $DK*K^- \hookrightarrow \delta_\gamma \oplus $ some translates of constant
perverse sheaves, this implies ${\cal H}^0(DK*K^-) \hookrightarrow {\cal
H}^{0}({}^p\delta_\gamma)$. Since ${}^p\delta_\gamma$ is a perverse sheaf, its cohomology sheaf
${\cal H}^0({}^p\delta_\gamma)$ is isomorphic to the direct sum of all perverse summands with
zero-dimensional support. Hence this is a skyscraper sheaf on $X$. Therefore $Aut(K)$ is
finite. Concerning property 4). ${\cal H}^0(DK*K^-)$ is nonzero for any perverse sheaf $K\neq
0$. $K*\delta_C=0$ implies $DK*\delta_C=0$. Since $K\hookrightarrow \delta_\alpha$ for some
$\alpha$,  we get $K^- \hookrightarrow \delta_\alpha^- \equiv
T^*_{x_0}\delta_{\beta}\hookrightarrow \delta_C *\cdots \delta_{C+x_0}$ for some $\beta$ (using
duality). Hence $DK*K^- \hookrightarrow (DK*\delta_C)* \cdots \delta_{C+x_0} $, where the
convolution on the right side now is a sum of translates of constant perverse sheaves. Hence
${\cal H}^0(DK*K^-)$ is a constant sheaf contradicting 3). The last properties follow from
Verdier duality $D\delta_\alpha\cong \delta_\alpha$ respectively from Riemann-Roch duality
${}^p \delta_\alpha^\vee \cong {}^p\delta_\beta$.

\bigskip\noindent
By abuse of notation we will identify the two categories ${\cal BN}'$ and ${\cal BN}$ in the
following.

%\bigskip\noindent
%\begin{Lemma} \label{abg}{\it The convolution of two BN-sheaves is the direct sum of a perverse sheaf and a translate
%of constant perverse sheaves}. \end{Lemma}

%\bigskip\noindent
%\underbar{Proof}: This can be immediately reduced to the
%corresponding statement, that the convolution $\delta_\alpha
%*\delta_\beta$ is the direct sum of a perverse sheaf and a
%translate of constant perverse sheaves. This was shown in the last
%two sections.

\bigskip\noindent
\section{Preliminary multiplicity formulas}\label{MF}

\bigskip\noindent
Let $Perv(X)'$ be a subcategory of $Perv(X)$ such that for all $K,L\in Perv(X)'$ the complex
$K*L^\vee\in D_c^b(X,\overline\Q_l)$ is a direct sum of a perverse sheaf in $Perv(X)'$ and a
sum $\oplus_\nu T_\nu[\nu]$ of translation-invariant perverse sheaves $T_\nu\in Perv(X)$. Then
the  quotient category $\overline{Perv}(X)'$  in $\overline{Perv}(X)$ is closed under
convolution
$$ *:\ \overline{Perv}(X)' \times \overline{Perv}(X)' \to \overline{Perv}(X)' \ ,$$
a typical example being the category $Perv(X)'={\cal BN}$ of BN-sheaves.

\bigskip\noindent
For $\delta_C\in Perv(X)'$ and $K\in Perv(X)'$ put $L=K*\delta_C$. $L*L^\vee\cong K*K^\vee
* (1 \oplus \Omega)$ for
 $\delta_C*\delta_C^\vee=1\oplus \Omega$ and $\Omega\cong \Omega^\vee$ by section \ref{CundA}. Now
$ Hom_{\overline{Perv}(X)'}(A,B) \ \cong\ \Gamma_{\{0\}}(X,{\cal H}^0(A*B^\vee)^\vee)$ using
corollary \ref{zuvor2} implies
$$ Hom_{\overline{Perv}(X)'}(L,L) \ \cong\ Hom_{\overline{Perv}(X)'}(K,K) \oplus
Hom_{\overline{Perv}(X)'}(\Omega,K*K^\vee) \ .$$ For  semisimple $K\in Perv(X)'$ of geometric
origin let $K=\bigoplus_{i=1}^r n_i\cdot K_i$ be the decomposition into irreducible components
$K_i$. Assume none of the irreducible perverse sheaves $K_i$ to be translation-invariant. Then
$$ {}^p m(K) := \sum_{i=1}^r n_i^2 = dim(Hom_{\overline{Perv}(X)'}(K,K)) \ .$$
Furthermore $$L=\delta_C*K = P\oplus T_K$$ is a semisimple perverse sheaf on $X$ by theorem
\ref{Geomorigin}. Without restriction of generality  $P$ does not contain a
translation-invariant irreducible constituent, and
$$ {}^p m(L) := {}^p m(P) = dim(Hom_{\overline{Perv}(X)'}(L,L)) \ .$$%we have
%$$ m(L) = dim\ {\cal H}^0\Bigl(L* L^\vee\Bigr)_0 \quad , \quad
% m_X(L) = dim\ {\cal H}^0\Bigl(L* L^\vee\Bigr)_\eta \ ,$$ where $\eta$ is a point of general
%position in $X$.
%\bigskip\noindent
%\underbar{Convolution}: Let $K=L*\delta_C$.

\bigskip\noindent
\begin{Lemma} \label{L}{\it Suppose $K\in Perv(X)'$ is a semisimple perverse sheaf on $X$ of geometric
origin. Then for $L=K*\delta_C$
$$ {}^p m(L)\ = \ {}^p m(K)\ +\  dim\ Hom_{Perv(X)}\bigl(\Omega,
{}^pH^0(K*K^\vee)\bigr) \ . $$  If $C$ is not hyperelliptic, then
$$  dim\ Hom_{Perv(X)}\bigl(\Omega, {}^pH^0(K*K^\vee)\bigr) \ \leqq \
\frac{1}{2g}\cdot dim\ {\cal H}^{-1}\bigl(K*K^\vee\bigr)_0 \ .$$ } \end{Lemma}

%\bigskip\noindent
%\underbar{Remark}: For $L=\delta_\alpha$ and $\alpha_1\leq g-1$
%the assumptions of the last lemma are satisfied. Also notice,
%$L^\vee \oplus T_\alpha = L_{\beta}$, where $\beta =
%(\chi-\alpha_r,...,\chi-\alpha_1)$ and
%$\alpha=(\alpha_1,..,\alpha_r)$. .

\bigskip\noindent
\begin{Corollary}\label{Mul} {\it If $K$ in the last lemma is of 0-type, then
$ {}^p m(K) = {}^p m(L)$.} \end{Corollary}

\bigskip\noindent
\underbar{Proof}: ${}^pH^0(K*K^\vee)$ is again of $0$-type. Either $\Omega$ is irreducible with
$H^{-2}(X,\Omega)\cong \overline\Q_l$, or it is a direct sum of $\delta_{C-C}$ (with the same
property $H^{-2}(X,\delta_{C-C})\cong \overline\Q_l$) and $A$ with $H^{-1}(X,A) \cong
\overline\Q_l^{2g}$. Hence any map from $\Omega$ to ${}^pH^0(K*K^\vee)$ must be zero by
Gabber's theorem. Now apply lemma \ref{L}.

\bigskip\noindent
\begin{Corollary} \label{swip}{\it
$ {}^p m(K*\delta_C)={}^p m(K*\delta_{-C})$ under the assumptions of lemma \ref{L}.}
\end{Corollary}

\bigskip\noindent
Using the argument of section \ref{Compu} one also obtains from the K\"unneth formula

\bigskip\noindent
\begin{Lemma}\label{tri} {\it Suppose $K\in Perv(X)'$ is a semisimple perverse sheaf on $X$. Then
$$ dim(Hom_{Perv(X)}(\delta_X,K*\delta_{C}))\ = \  dim(
H^{g-1}(X,K))\ .$$} \end{Lemma}

\bigskip\noindent
Applied for $L^\vee =K^\vee *\delta_{-C}$ instead of $L=K*\delta_C$ gives
$$ dim(
H^{g-1}(X,K)) = dim( H^{g-1}(X,K^\vee)) \ .$$

%
%Here we used $\Omega = \delta_{C-C}\oplus A$ with $H^{\nu}(X,A)=0$ for $\nu\leq -2$,
%$H^{-2}(X,\delta_{C-C})=0$ for $\nu<-2$ and $H^{-2}(X,\delta_{C-C})\cong\overline\Q_l$. In our
%case $n(2)= h^{g-1}(X,L_i)h^{g-1}(X,L_j^\vee)$ by our computation of ${\cal H}^{-2}(T_{ij})$
%above. This proves the claim.
%
%\bigskip\noindent
%\underbar{For generic}  $x\in X$ this is
%$$ {\cal H}^0(\Omega * L*L^\vee)_x = \ \bigoplus_{i,j}\ n_in_j\cdot \Bigl(\
% {\cal H}^0(\Omega * {}^pH^0\bigl(L_i*L^\vee_j\bigr)_{nc})_x\ \oplus\
%{\cal H}^0(\Omega * T_{i,j})_x \ \Bigr) \
% $$ $$ =  \ \bigoplus_{i,j}\ n_in_j\cdot {\cal H}^0(\Omega * T_{i,j})_x  \ = \ \bigoplus \ \
%H^{g-1}(X,L)\otimes H^{g-1}(X,L^\vee) \ ,$$ where the last
%equality follows from the computation made above. The first
%equality follows from the fact, that ${\cal H}^0(\Omega * M)_x$
%vanishes for an arbitrary perverse sheaf $M$, if $x$ is in general
%position.
%
%
%\bigskip\noindent
%\underbar{$K$ is perverse}: This follows from the assumption, that
%$L$ is a sum of BN-sheaves. Under this assumption $K=L*\delta_C$
%and $K=\delta_{-C}$ are automatically perverse by \ref{APP}. In
%other words, we have shown
%
%
%
%
%%\bigskip\noindent
%%\underbar{Remark}: For $L=\delta_\alpha$ and $\alpha_1\leq g-1$
%%the assumptions of the last lemma are satisfied. Also notice,
%%$L^\vee \oplus T_\alpha = L_{\beta}$, where $\beta =
%%(\chi-\alpha_r,...,\chi-\alpha_1)$ and
%%$\alpha=(\alpha_1,..,\alpha_r)$. .

\bigskip\noindent

\section{Dimensions}
A collection of numbers $d(\alpha)\in \Z$ for the partitions $\alpha$, such that
$$d(\alpha)\cdot d(\beta) = \sum m_{\alpha\beta}^\gamma \cdot d(\gamma)$$ holds for all partitions
$\alpha$ and $\beta$ will be called a system of dimensions. Notice $d(0)=1$. For partitions
$\alpha$ of degree $deg(\alpha)=r$ parameterizing the classes of irreducible representations
$\sigma_\alpha$ of the symmetric group $\Sigma_r$, recall the involution $\alpha\mapsto
\alpha^*$. If $\alpha$ parameterizes the class of $\sigma_\alpha$, then $\alpha^*$
parameterizes the class of $\sigma_\alpha\otimes sign$. In particular
$$ (r,0,..,0)^*= ([r])^*=([1]^r)\ .$$

\bigskip\noindent
For  a system of dimensions $d(\alpha)$ also  $\tilde d(\alpha)=d(\alpha^*)$ is a system of
dimensions, which is an immediate consequence of
$$m_{\alpha\beta}^\gamma = m_{\alpha^*\beta^*}^{\gamma^*}\ .$$
This holds, since the sign character on $\Sigma_a\times \Sigma_b\hookrightarrow \Sigma_{a+b}$
restricted to $\Sigma_a\times \Sigma_b$  is the tensor product of the sign characters of $S_a$
and $S_b$. By induction on the lexicographic ordering, a system of dimensions is uniquely
determined by the numbers $d(\alpha)$ for $\alpha$ running over the partitions $(r,0,..,0)^*$
for $r\in \N$. Hence

\bigskip\noindent
\begin{Lemma} {\it Two systems of dimensions $d_1(\alpha)$ and $d_2(\alpha)$ coincide, if
$d_1([r]) = d_2([r])$ holds for all integers $r$.}  \end{Lemma}

\bigskip\noindent
In our case the Euler characteristics $$d_1(\alpha)= \sum_\nu  (-1)^\nu dim(
H^\nu(X,\delta_\alpha))$$ are a system of dimensions (depending only on the genus g). This
follows from the K\"{u}nneth formula. $\chi(S^a(H_+)) = (-1)^a \cdot (a+1)$ implies $d_1([r]) =
\sum_{a+b=r} (-1)^a \sigma^a \cdot \lambda^b(H_-)$ for $\lambda^b = dim(\Lambda^b(H_-))$ and
$\sigma^a = dim( S^a(H_+)) = a+1$ by lemma \ref{dimension}. For a vectorspace $H$ of dimension
$2g-2$ choose some isomorphism $H_+\oplus H\cong H_-$ (notice $g\geq 2$). Then
$\lambda_t(H_-)=\lambda_t(H_+)\lambda_t(H)$ and $\lambda_t(H_+)^{-1} = \sigma_{-t}(H_+)$ for
the power series $\lambda_t(.) = \sum_{i=0}^\infty \lambda^i(.)t^i$ and
$\sigma_t(.)=\sum_{i=0}^\infty \sigma^i(.)t^i $. Hence
$$ d_1\bigl((r,0,..,0)\bigr) \ =\ \lambda^r(H) = {2g-2 \choose r} \ ,$$
the binomial coefficient being zero for $r>2g-2$.

\bigskip\noindent
\begin{Corollary} {\it The systems of dimensions defined by
\begin{enumerate} \item $d_1(\alpha) = \sum_\nu (-1)^\nu  dim(H^\nu(X,\delta_\alpha))$
, and \item $d_2(\alpha) = dim(H^{\alpha^*})$
\end{enumerate} coincide. Hence for all partitions $\alpha$, for which $\alpha_1< 2g-2$ holds, the dimensions
$$d_1(\alpha)= \sum_\nu  (-1)^\nu dim(H^\nu(X,\delta_\alpha)) >0 \quad , \quad \alpha_1 <\chi $$ are
positive.} \end{Corollary}

\bigskip\noindent
\underbar{Proof}: The last assertion follows from $dim(H^{\alpha^*})>0$, which in turn follows
from the well known classification of irreducible representations of the group $Sl(2g-2)$.

\chapter{Some special cases}

\goodbreak
\section{Hyperelliptic curves II}

\bigskip\noindent
By the multiplication law (section \ref{ML}) we always have
$$ \delta_{g-1} * \delta_{g-1} = \delta_{2g-2}
\oplus \delta_{2g-3,1} \oplus \cdots \oplus \delta_{g-1,g-1} \ .
$$  Also
by the multiplication law $ \delta_{2g-3}*\delta_1 = \delta_{2g-2} \oplus \delta_{2g-3,1}$.
Therefore $\delta_{2g-3}*\delta_1$ is a direct summand of $\delta_{g-1}*\delta_{g-1}$. Combined
with $\delta_{2g-3} \equiv \delta_{\kappa -C}$ (lemma \ref{RR}) this implies
$\delta_{2g-3}*\delta_1 \equiv \delta_{\kappa - C}*\delta_C = \delta_{\kappa +C-C} \oplus ..$,
where the latter follows from \ref{CundA}. Hence $ \delta_{\kappa+C-C} \hookrightarrow
\delta_{g-1}*\delta_{g-1}$. By corollary \ref{MA} (iii) we have $\delta_{g-1}=\delta_\Theta$ in
the non-hyperelliptic case. So for non-hyperelliptic curves the above identities imply
$$\delta_{\kappa+C-C} \hookrightarrow \delta_\Theta*\delta_\Theta \ .$$
As a consequence the perverse sheaf $\delta_{\kappa+C-C}$ satisfies the properties of theorem
\ref{Key} below in the non-hyperelliptic case. This follows from section \ref{CundA} and
corollary \ref{Thet}. See also \cite{W}.

\bigskip\noindent
{\it The hyperelliptic case}. For hyperelliptic curves $C$ this argument brakes down. Instead
for $1\leq r\leq g-1$ consider the proper morphisms
$$ f_r : W_r \times W_r \to X $$
$$ (x,y) \mapsto x-y \ .$$
We have $C=e-C$, where for simplicity we assume the point $e$ to be $0$. This can be achieved
by a suitable translation of the curve $C$. Then $W_r = -W_r$. The inverse image
$f_r^{-1}(C-C)$ in $W_r\times W_r$ of the surface $C-C\subseteq X$
$$ \xymatrix{ (x,y) + W_{r-1}\ar[d] \ar[r] & f_r^{-1}(C-C) \ar@{->>}[d]\ar@{^{(}->}[r]
& W_r \times W_r \ar[d]^{f_r}\cr \{x-y\}\ar@{^{(}->}[r] & C-C
\ar@{^{(}->}[r] & X \cr }
$$
admits a morphism
$$ g: \ C\times C \times W_{r-1} \longrightarrow f_r^{-1}(C-C) $$
defined by $(x,y,D)\mapsto (x+D,y+D)\in W_r\times W_r$, such that $f_r\circ g$ surjects onto
$C-C$. Notice $C+W_{r-1}\subseteq W_r$.

\bigskip\noindent
For later reference observe, that if $\theta$ denotes the hyperelliptic involution, then
$\nu_r(a,b)= (\theta(a),\theta(b))$ acts on $W_r\times W_r$, such that $f_r\circ \nu_r= f_r$
and the following diagram commutes
$$ \xymatrix{ (C\times C)\times W_{r-1} \ar[d]_{\nu_1\times \theta}\ar[r]^-g & W_r\times W_r \ar[d]^{\nu_r}\cr
(C\times C)\times W_{r-1} \ar[r]^-g & W_r\times W_r \cr } $$ {\bf Claim}: {\it $g(C\times
C\times W_{r-1})$ is an irreducible component of $f_r^{-1}(C-C)$ of the highest possible
dimension $r+1$, and for some open dense subset $U\subseteq C-C$ the  fibers of $f_r$ are of
dimension $\leq r-1$ over $U$.}

\bigskip\noindent
Taking this for granted (see lemma \ref{35} below) let us consider
$$ \delta_{W_r}* \delta_{W_r} = \delta_{W_r}* \delta_{-W_r} = Rf_{r,*}(\delta_{W_r}\boxtimes
\delta_{W_r}) = Rf_{r,*}(\overline\Q_{l,W_r\times W_r}[2r]) \ .
$$
For the last equality recall, that $\delta_{W_r}=\overline\Q_{l,W_r}[r]$ holds in the
hyperelliptic case by lemma \ref{hy}. Choose $U\subseteq C-C$ as above. Consider the cohomology
sheaf of $\delta_{W_r}*\delta_{W_r}$ in degree -2 over the open dense subset $U$ of $C-C$
$$ {\cal H}^{-2}(\delta_{W_r}* \delta_{W_r})\big\vert_U =
R^{2r-2}f_{r,*}(\overline\Q_{l,W_r\times W_r})\big\vert_U \ .$$ The proof of lemma \ref{35}
shows, that $\delta_{W_r}* \delta_{W_r}$ is a direct sum of a perverse sheaf $K$ on $X$ and a
sum of translates of constant perverse sheaves, such that $${\cal H}^{-2}(\delta_{W_r}*
\delta_{W_r}) ={\cal H}^{-2}(K) \ .$$ Perversity implies, that the support of ${\cal
H}^{-2}(K)$ has dimension $\leq 2$. Furthermore any contribution of ${\cal H}^{-2}(K)\vert_U$
over the generic point of $U$ necessarily comes from a perverse direct summand $L\in Perx(X)$
of $K$ with support in $C-C$. This now allows to detect perverse direct summands of
$\delta_{W_r}* \delta_{W_r}$ by considering the higher direct image sheaf
$R^{2r-2}f_{r,*}(\overline\Q_{l,W_r\times W_r})\big\vert_U$. In fact by the usual excision
arguments this higher direct image contains the subsheaf
$$ R^{2r-2}(f_r\circ g)_*(\overline\Q_l)\big\vert_U \neq 0 \ ,$$
since as stated above the image of $g$ is an irreducible component of highest dimension $r+1$
in $f_r^{-1}(C-C)$ with relative dimension $r-1$ over $U$.

\bigskip\noindent
The morphism $(f_r\circ g)(x,y,D)=x-y$ can be identified with the cartesian product $ f_r \circ
g = f_1 \times s: \ (C\times C)\times W_{r-1} \to (C-C)\times Spec(k)=C-C$ , where $f_1=f_r$
for $r=1$, and where $s:W_{r-1}\to Spec(k)$ is the structure morphism. Since $f_1$ is finite
over a suitable chosen $U$, $R^if_{1,*}(\overline\Q_l)\big\vert_U =0$ holds for $i\geq 1$ over
a suitably small open dense $U\subseteq C-C$. For such a choice of $U$ from the K\"unneth
formulas, applied for the proper morphisms $f_1\times s$, we derive
$$R^{2r-2}(f_r\circ g)_*(\overline\Q_l)\vert_U =
R^0f_{1,*}(\overline\Q_l)\big\vert_U \otimes R^{2r-2}s_*(\overline\Q_l)\big\vert_U\ .
$$  Since $W_{r-1}$ is irreducible of dimension $r-1$, the trace map gives an isomorphism
$R^{2r-2}s_*(\overline\Q_l)\big\vert_U \cong \overline\Q_{l,U}$ ignoring Tate twists. Hence we
obtain an injective sheaf morphism
$$ f_{1,*}(\overline\Q_l)\vert_U \hookrightarrow {\cal H}^{-2}(\delta_{W_r}*
\delta_{W_r})\big\vert_U \ .$$ They same holds for $r=1$, where this map is an isomorphism with
$f_{1,*}(\overline\Q_l)\vert U = \Q_{l,U}\oplus E\vert_U$. See section \ref{CundA}. Hence
$$ \Q_{l,U} \oplus E\vert_U \cong {\cal H}^{-2}(\delta_1*\delta_1)\big\vert_U \hookrightarrow
{\cal H}^{-2}(\delta_{W_r}*\delta_{W_r})\big\vert_U \ .$$ By the support dimension $dim(U)=2$
and the definition of middle perversity one concludes, using perverse continuation from $U$ to
$C-C$ as already mentioned above, that
$$ {}^p H^0(\delta_{W_r} * \delta_{W_r}) = \delta_{C-C} \oplus A \oplus \cdots
\ $$ holds, where $A=\delta_E=\delta_{1,1}$ as defined in section \ref{CundA}.

\bigskip\noindent
\begin{Lemma}\label{35} {\it Suppe $1\leq r\leq g-1$. The fibers of $f_r$ have dimension $\leq r$. There exists an open dense subset $U\subseteq C-C$, such that the fibers of $f_r$
over $U$ have dimension $r-1$. In particular $dim(f_r^{-1}(C-C))\leq r+1$. Hence $g(C\times
C\times W_{r-1})$ is an irreducible component of $f_r^{-1}(C-C)$ of highest dimension.}
\end{Lemma}

\bigskip\noindent
\underbar{Proof}:  $f_r:W_r\times W_r \to X$ is proper. Hence for a closed point $x$ of $X$ we
have $dim(f_r^{-1}(x)) > r-1 \Longleftrightarrow \exists n > 2r-2$ such that $R^n
f_{r,*}(\overline\Q_l)_x \neq 0$ by the proper base change theorem. Since
$Rf_{r,*}(\overline\Q_l)[2r]=\delta_{W_r}*\delta_{W_r}$ as shown above, this is equivalent to
$\exists m>-2$ such that ${\cal H}^m(\delta_{W_r}*\delta_{W_r})_x\neq 0$. Since
$\delta_{W_r}*\delta_{W_r} \hookrightarrow \delta_r*\delta_r$ (lemma \ref{hy}) and since
$\delta_r*\delta_r$ is a direct sum of a perverse sheaf
$$P \hookrightarrow \bigoplus_{\alpha=(\alpha_1,\alpha_2),deg(\alpha)=2r} \delta_{\alpha_1,\alpha_2}^p$$ and a
sum $T$ of translates of the constant perverse sheaves $\delta_X$ (corollary \ref{Zerlegung}),
we obtain the following estimates: By perversity ${\cal H}^m(P) =0$ for $m>0$ and  ${\cal
H}^m(P)\vert_U =0$ for $m\geq -1$ for some suitable chosen dense open subset $U\subseteq C-C$.
Concerning the complex $T$. By corollary \ref{estimate} the maximal depth of translation of
$\delta_X$ within $T$ can be estimated by the maximum of the degrees  $deg_t(P_\alpha(t))$ with
$\alpha=(\alpha_1,\alpha_2)$, such that $\delta_\alpha \hookrightarrow \delta_r*\delta_r$.
Since $deg(\alpha)=2r$, we get $\alpha_1\leq 2r$. Hence $\alpha_1-g \leq 2r-g \leq r-1$ with
equality only for $r=g-1$. This proves $\alpha_1-g \leq g-2$, with equality only for $r=g-1$.
By corollary \ref{estimate} therefore $deg_t(P_\alpha(t)) \leq \alpha_1 - g \leq g-2$. Hence
${\cal H}^m(T)$ vanishes for $m\geq -1$. We conclude, that for $f_r$ the maximal fiber
dimension is $r$, and the fiber dimension over $U$ is $\leq r-1$. \qed

%
%It is enough to show, that $\delta_r*\delta_r$ is a perverse sheaf and a direct sum $T$ of
%translates of constant sheaves $\delta_X$, such that ${\cal H}^\nu(T)=0$ for $\nu\geq -2$. Then
%$R^{2r+\nu}s_*(\overline\Q_l)\big\vert_U$ vanishes for $\nu > -2$ at the generic point $\eta$
%of $U$, which implies that the dimension of $f_r^{-1}(\eta)$ is $\leq r-1$. However
%$R^{2r+\nu}f_{r,*}(\overline\Q_l) = \delta_{W_r}*\delta_{-W_r} = \delta_{W_r}*\delta_{W_r}
%\hookrightarrow \delta_r*\delta_r$. If we use the general result, proved in the later section
%\ref{PD}, that $\delta_{r}*\delta_r$ is perverse up to some $T$, where $T$ is a sum of
%translates of constant sheaves on $X$, it only remains to control the cohomology of $T$. In the
%notation of section \ref{PD} it is enough to estimate the degree of the polynomial
%$P_\alpha(t)$ for $\alpha =(\alpha_1,\alpha_2)$ of degree $deg(\alpha)=2r$. It follows from
%\ref{PD}, that this degree is $\leq \alpha_1 - g \leq 2r - g$. But $2r - g \leq r-1$ for $r\leq
%g-1$ with equality only for $r=g-1$. This shows, that ${\cal H}^{\nu}(T)= 0$ for $\nu \geq -2$
%except in the case $r=g-1$, where this only holds for $\nu > -2$. In fact, in this case
%$\delta_{g-1}*\delta_{g-1}$ contains a translate $\delta_X[-g+2]$. In fact, this translate
%restricts to the sheaf $\delta_{C-C}$ on $C-C$. In any case, however, this implies that the
%fibers of $f_r$ over an open dense subset $U$ of $C-C$ have dimensions $\leq r-1$.

%In fact, for this we can consider the stalks at the generic point
%$\eta_X$ of $X$. Clearly, the fiber dimension at $\eta_X$ is $\leq
%$
\bigskip\noindent
We remark, that in fact $\delta_{g-1}*\delta_{g-1}$ contains a translate $\delta_X[-g+2]$,
which restricts to the sheaf $\delta_{C-C}$ on $C-C$.

\bigskip\noindent
Recall from section \ref{CundA}, that the support of the irreducible perverse sheaf
$\delta_{1,1}=A$ is a translate of $C-C$, such that ${\cal H}^{-1}(A) = H^1(C)\otimes
\delta_{0}$. Since by our conventions $\Theta=-\Theta$, we obtain from the preceding discussion
now also in the hyperelliptic case

\bigskip\noindent
\begin{Theorem} \label{Key}{\it There exists an unique irreducible perverse sheaf  $A\hookrightarrow \delta_\Theta * \delta_\Theta$
for $\delta_{\Theta}=\delta_{W_{g-1}}$ characterized by one of the following equivalent
properties
\begin{enumerate}
\item ${\cal H}^{-1}(A)$ is non-zero, but not a constant sheaf.
\item ${\cal H}^{-1}(A)\cong H^1(C)\otimes \delta_{\{\kappa\}}$ is a skyscraper
sheaf concentrated in the point $\kappa\in X$ defined by
$\Theta=\kappa -\Theta$.
\end{enumerate}
Furthermore the support of $A$ is a translate of the subvariety $C-C\subseteq X$.}
\end{Theorem}

\bigskip\noindent
\begin{Corollary}\label{triple} {\it Suppose $C$ is hyperelliptic. Then $\delta_C
* \delta_C *\delta_C$ has seven irreducible perverse constituents.} \end{Corollary}

\bigskip\noindent
\underbar{Proof}:  $\delta_C * \delta_C *\delta_C = \delta_3 \oplus 2\delta_{2,1} \oplus
\varepsilon_3$, and $\delta_3=\delta_{W_3}\oplus \delta_C$. We will see, that $\varepsilon_3$
is irreducible.  From section \ref{Hyperelliptic} we know  $\delta_2=\delta_{W_2}\oplus
\delta_0$. Hence $\delta_3\oplus \delta_{2,1}=\delta_2*\delta_1= (\delta_{W_2}\oplus
\delta_0)*\delta_C= (\delta_{W_2}*\delta_{W_1})\oplus \delta_C$. Lemma \ref{hypell} implies,
that $\delta_C$ is a summand of $\delta_{W_2}*\delta_{W_1}$. Since $\delta_3 =
\delta_{W_3}\oplus \delta_C$, therefore
 $$ \delta_{2,1} =
\delta_C \oplus \delta_{2,1}' \ .$$ Recall $\varepsilon_2*\delta_{-C}\equiv \varepsilon_1\oplus
\delta_{2g-3,1}$ up to translates of constant sheaves. Further notice ${}^p
\delta_{2g-3,1}\oplus \delta_{\{e\}} \cong \delta_C*\delta_{e-C} = \delta_{W_2}\oplus A\oplus
\delta_{\{e\}}$. Hence ${}^p m(\varepsilon_2*\delta_{-C})=3$. For $\varepsilon_2*\delta_C
=\delta_{2,1}\oplus \varepsilon_3$ the multiplicity formula ${}^p m(\varepsilon_2*\delta_C) =
{}^p m(\varepsilon_2*\delta_{-C})$ of lemma \ref{swip} therefore implies ${}^p
m(\varepsilon_2*\delta_{C})=3$. Hence ${}^p m(\varepsilon_3)=1$ and ${}^p m(\delta_{2,1})=2$,
and $\delta_{2,1}'$ is irreducible. This implies the assertion.

\section{The perverse sheaf $\varepsilon_{g-1}$}
In this section we suppose, that  $C$ is not hyperelliptic. Then $\varepsilon_{g-1}$ is an
irreducible perverse sheaf on $X$ by corollary \ref{MA}. We compare the properties of
$\delta_{g-1}$ with those of $\delta_{g-1}$. The proofs are geometric. The results of this
section are obtained independently and more generally  from theorem \ref{main}.

\bigskip\noindent

\bigskip\noindent
\begin{Lemma} {\it If $C$ is not a hyperelliptic curve, $\varepsilon_{2g-2}$ does not contain
the skyscraper sheaf $\delta_{\{\kappa\}}$.} \end{Lemma}

\bigskip\noindent
\underbar{Proof}: The fiber $p_{2g-2}^{-1}(\kappa)=\Proj^{g-1}$,
is canonically isomorphic to the linear system of the canonical
class of $C$
$$ p_{2g-2}^{-1}(\kappa) = \vert K \vert \cong \Proj^{g-1} \ .$$
 We
have to show ${\cal H}^0(\Proj^{g-1},\varepsilon_{2g-2} )= 0$, or
equivalently that the  sheaf $F_\varepsilon$ (genic of rank one on
$C^{(2g-2)}$) has nonconstant restriction to this fiber.

\bigskip\noindent
We fix a canonical divisor $K$ on $C$. By the theorem of Riemann-Roch one can choose an
arbitrary positive divisor $D\in C^{(g-1)}$ of degree $g-1$ on $C$, such that $h(C,K-D)= h^0(D)
\geq 1$. For a divisor $D$ of general position in $C^{(g-1)}$ evidently $h^0(D)=1$. In other
words, there exist $U\subseteq C^{(g-1)}$ dense open, such that for $D\in U$ there is a unique
positive divisor $\theta(D)\in C^{(g-1)}$, such that $D+\theta(D)=K$. This defines a regular
morphism $U\to C^{(g-1)}$. Since $C^{(g-1)}$ is normal, it extends to a morphism $\theta:
C^{(g-1)} \to C^{(g-1)}$ making the diagram
$$ \xymatrix@+0.5cm{ C^{(g-1)}\ar[d]_{p_{g-1}} \ar[r]^\theta &
C^{(g-1)}\ar[d]^{p_{g-1}} \cr W_{g-1} \ar[r]^{x\mapsto \kappa-x} & W_{g-1} \cr} $$ commutative.
Obviously then $\theta^2=id$. Hence $\theta$ is an automorphism of $C^{(g-1)}$ of order two.

\bigskip\noindent
Hence there exists the morphism
$$ \xymatrix{ C^{(g-1)} \ar[r]^-f &  \vert K\vert\cong\Proj^{g-1} \subseteq C^{(2g-2)} \cr} $$
defined by $f(D)= D + \theta(D) \in C^{(2g-2)}$. By dimension
reasons this is a surjective morphism.

\bigskip\noindent
Now suppose $F_\varepsilon$ restricted to $\vert K\vert$ is
trivial. Then, the restriction via $f$ is trivial. Since $f$
factorizes over $C^{(g-1)} \to C^{(g-1)}\to C^{(g-1)}\times
C^{(g-1)} \to C^{(2g-2)}$, where the first map is $id\times
\theta$, and since $F_{\varepsilon}$ on $C^{(2g-2)}$ restricts to
$F_\varepsilon \boxtimes F_\varepsilon$ on $C^{(g-1)} \times
C^{(g-1)}$, this would imply  $F_\varepsilon \cong
\theta^*(F_\varepsilon)$ on $C^{(g-1)}$.

\bigskip\noindent
Suppose this were the case. Then, the ramification locus $V$ of
$F_\varepsilon$ and the ramification locus  $V'$ of
$\theta^*(F_\varepsilon)$ in $C^{(g-1)}$ must coincide. The
divisor $F$ is the image of the finite surjective map
$$ C\times C^{(g-3)} \to C^{(g-1)}$$
defined by $(P,D)\mapsto 2P+D$. But then $\theta(2P+D) = 2P' + D'$
must hold for some $P'\in C$ and some $D'\in C^{(g-3)}$. For $P$
and $D$ in general position (contained in $U\cap \theta(U)$)
however the points $P'\in C$ and $D'\in C^{(g-3)}$ defined by this
are uniquely determined. This defines a birational map $C\times
C^{(g-3)} \to C\times C^{(g-3)}$, which is regular on $U\cap
\theta(U)$. Again by normality this birational map extends to a
regular morphism, making the diagram
$$ \xymatrix { C\times C^{(g-3)} \ar[d]\ar[r]^F & C\times C^{(g-3)} \ar[d]\cr
V \ar[d]\ar[r] & \theta(V)=V \ar[d]\cr C^{(g-1)} \ar[r]^\theta &
C^{(g-1)} \cr}
$$ commutative, such that the compositions
$C\ni P\mapsto \varphi(P,D)= pr_C(F(P,D))$ defined by the
composition $C\to C\times C^{(g-3)} \to C\times C^{(g-3)} \to C$
extend to automorphisms of $C$. Notice, then $F^2=id$. Since
$g\geq 2$, also $\varphi(P,D)=\varphi(P)$ does not depend on $D\in
C^{(g-3)}$, and defines an involution of $C$. Hence
$F(P,D)=(\theta(P),\psi_P(D))$, such that $\psi_{\theta(P)}\circ
\psi_P =id_{C^{(g-3)}}$. In particular $\psi_P: C^{(g-3)}\to
C^{(g-3)}$ is an automorphism of $C^{(g-3)}$.  But then, if we put
$P=0\in C$ and $Q=\theta(P)$, we get from the commutative diagram
above $$ \kappa - W_{g-3} = 2Q + W_{g-3} \ .$$ Now $$(W_r +u) \cap
(W_r +v) = W_{r+1}^2 \cup (W_{r-1} + u+v) $$ or
$$W_r  \cap
(W_r +v-u) = (W_{r+1}^2 -u) \cup (W_{r-1} +v) $$ holds for all $r\leq g-1$ and all $u\neq v\in
C$. See \cite{Gunning}, page 141. Furthermore $dim(W_{r+1}^2) \leq r-1$  holds for all $r \leq
g-2$ with equality reached only for  hyperelliptic curves (\cite{Gunning}, p.142). Hence for
non-hyperelliptic curves $C$ we conclude from some equality $-W_r = x_r + W_r$ where $x_r\in X$
and $1\leq r\leq g-3$, that
$$ - W_{r-1} - u - v \subseteq (-W_r - u)\cap (-W_r - v)$$ $$ =
(W_r + x_r - u)\cap (W_r+x_r - v) = x_r -u + W_r \cap (W_r+u - v)
$$ and
$$  x_r -u + (W_{r-1} +u) = x_r + W_{r-1} $$
are both the same irreducible component of maximal dimension of
$x_r -u + W_r \cap (W_r+u - v)$. Hence
$$ - W_{r-1} - u - v  = x_r -u + (W_{r-1} +u) = x_r + W_{r-1} \ .$$
This gives immediately a contradiction. This proves the lemma.

\bigskip\noindent
Since $C^{(g-1)} \to W_{g-1}$ is a birational map, the proof of $F_\varepsilon \not\cong
\theta(F_\varepsilon)$ on $C^{(g-1)}$ implies

\bigskip\noindent
\begin{Corollary} \label{cor18}{\it If $C$ is not hyperelliptic $\varepsilon_{g-1}$ and
$\varepsilon_{g-1}^-$ are not isomorphic.} \end{Corollary}

\bigskip\noindent
As a consequence

\bigskip\noindent
\begin{Corollary} {\it If $C$ is not hyperelliptic, then ${\cal H}^0(\delta_\alpha)=0$ holds for all $\alpha$ of degree $deg(\alpha)=2g-2$ with
$\alpha_1\leq 2$.} \end{Corollary}

\bigskip\noindent
\underbar{Proof}:    The convolution $\varepsilon_{g-1} * \varepsilon_{g-1}$ is  a direct sum
by the general multiplication formula $ \varepsilon_{g-1}
* \varepsilon_{g-1} \ = \ \sum_\alpha m(\alpha) \delta_\alpha \
$ for certain $\alpha$ of degree $2g-2$, for which $\alpha_1\leq 2\leq g-1$. Every
$\delta_\alpha$ with $deg(\alpha)=2g-2$ and $\alpha_1\leq 2$ appears with $m(\alpha)\neq 0$ in
this sum according to the Littlewood-Richardson rules. The claim follows, since ${\cal
H}^0(\varepsilon_{g-1}
* \varepsilon_{g-1}) = 0$ holds by corollary \ref{cor18} and \ref{cor1}.

\section{Other examples}

\bigskip\noindent
\begin{Lemma} {\it For a BN-sheaf $K$ ($K$ is irreducible say with
irreducible support $Y$ of dimension $d$) suppose $n_X(K*\delta_{C})\neq 0$. Then either
$d=g-1$ and $K=\delta_Y$ or $Y=X$. In both cases $ n_X(K*\delta_{C})= dim(H^{g-1}(X,E))$.}
\end{Lemma}

\bigskip\noindent
\underbar{Proof}: Lemma \ref{tri}  forces $d\geq g-1$, since $H^{g-1}(X,K)$ vanishes unless
$d\geq g-1$. $H^{g-1}(X,K)= H^{2g-2}(Y,E)=0$ also for $d=g-1$ unless the irreducible
coefficient system $E$ is trivial. Notice $E=j_*(E\vert_U)$, where $j:U=Y\setminus
V\hookrightarrow Y$ is the complement of a divisor $V$ in $Y$, such that $E\vert_U$ is a smooth
etale $\overline\Q_l$-sheaf on $U$. Then $H^{2g-2}(Y,E)= H^{2g-2}(Y,j_*(E\vert_U)) =
H^{2g-2}(Y,Rj_*(E\vert_U)) = H^{2g-2}(U,E\vert_U) = H^0_c(U,(E\vert_U)^\vee)$. Since $E\vert U$
is irreducible also on $U$, the dual etale sheaf $D(E\vert_U)$ is also irreducible. Hence
$H^0_c(U,D(E\vert_U))$ vanishes unless $E=\overline\Q_l$ is the constant coefficient system.

\bigskip\noindent
\begin{Lemma} {\it Suppose $C$ is not hyperelliptic and suppose $K=\delta_E$ is a BN-sheaf.
Suppose $Y$ is the irreducible support of $K$ of dimension $d$. Then
$$ {}^p m(K*\delta_{\pm C})\ \leqq\ {}^p m(K)\ +\ \frac{1}{2g}\cdot dim\
H^{2d-1}(Y,E\otimes \tilde E)\ $$ for $E=j_*(E\vert_U)$ and $\tilde E=j_*(D(E)\vert_U)$.}
\end{Lemma}

\bigskip\noindent
\underbar{Proof}: Let $l= dim(Hom_{Perv(X)}(\Omega, {}^p H^0(K*K^\vee)))$. Then $l\cdot \Omega
\hookrightarrow {}^p H^0(K*K^\vee)$ by Gabber's theorem. Comparing stalk dimensions using
$dim({\cal H}^{-1}(\Omega)_0)= 2g$ this implies $2g\cdot l \leq dim({\cal H}^{-1}({}^p
H^0(K*K^\vee)_0)$. On the other hand in \ref{help} we have shown,that the stalk ${\cal
H}^{-1}\bigl(K*K^\vee\bigr)_0$ is a quotient of the cohomology group $H^{2d-1}(Y,E\otimes
\tilde E)$. The assertion now follows from lemma \ref{L}.

\bigskip\noindent
\underbar{If $E$ has generic rank one}: For the rest of this section, we assume that $C$ is not
hyperelliptic. By assumption $E$ corresponds to a one dimensional $\overline\Q_l$-adic
representation $\rho_E$ of the etale fundamental group $\pi_1(Y\setminus V,y_0)$, where $V$ is
the ramification divisor of the corresponding covering of $Y$. Let $U=Y\setminus V$ be the open
complement, so that $E\vert_U$ is a smooth etale sheaf on $U$. Suppose $Y$ is a normal variety.
We can assume that $E\otimes \tilde E\big\vert_{V}$ is an etale sheaf in the ordinary sense
with support of dimension $< dim(V)$ and $dim(V)= dim(Y)-1$ or $V$ is empty (by the purity of
branch points). Hence $ H^{\nu}(V,E\otimes \tilde E\vert_V)=0$ vanishes for $\nu\geq 2d-2$.
From the long exact sequence of cohomology with compact support attached to $(U,X,Y)$ we
therefore obtain
$$ H^{2d-1}(Y,E\otimes \tilde E) \cong H^{2d-1}_c(U, E\otimes \tilde E\vert_U) \
.$$ But $E\otimes \tilde E\vert_U \cong \overline\Q_{l,U}$ holds, hence we  have an exact
sequence
$$ \to H^{2d-2}(Y)\overset{res}{\longrightarrow} H^{2d-2}(V)\to H^{2d-1}_c(U, E\otimes \tilde E\vert_U)
\to H^{2d-1}(Y) \to 0 .$$ If in addition the ramification divisor $V$ is irreducible, and its
proper transform $\tilde V \hookrightarrow \tilde Y$ in a desingularization $\tilde Y$ of $Y$
has nontrivial Chern class $c_1(\tilde V)\in H^2(\tilde Y)$ (or if $V$ is empty), then the
restriction map $res$ in this exact sequence is surjective. Together this implies $
H^{2d-1}(Y,E\otimes \tilde E) \cong H^{2d-1}(Y)$.

\bigskip\noindent
\underbar{Example}: $K=\delta_{r}=\delta_{W_r}$ for $C$ not hyperelliptic. Then $V$ is empty
and $d=r\leq g-1$ and $H^{2d-1}(Y)\cong H^1(X)$. Hence by the corollary above ${}^p m(\delta_r
*\delta_{\pm C}) \leq 2$. Therefore the multiplication formula implies, that for $1\leq r\leq
g-1$ the sheaves perverse sheaves $ {}^p \delta_{r,1}$ and ${}^p\delta_{2g-3,r} $ are
nontrivial irreducible perverse sheaves on $X$.

\bigskip\noindent
\underbar{Example}: For $K=\varepsilon_{g-1}$ the situation above also applies. Since $\tilde
V$ is the image of $2C\times C^{(g-3)} \to \tilde Y= C^{(g-1)}$, the conditions above are
satisfied. This shows by the corollary above, that the perverse sheaves $\varepsilon_g$ and
${}^p \delta_{2,1,..,1}$ (of degree $g$) are irreducible perverse sheaves on $X$, obtained as
the two BN-constituents of $\varepsilon_{g-1}*\delta_C$. Similar, convolution with
$\delta_{-C}$ proves, that the perverse sheaves $\varepsilon_{g-1}$ and ${}^p
\delta_{2g-3,1,...,1}$ (of degree $3g-4$) are irreducible perverse sheaves. Repeating the
argument above, inductively with $\varepsilon_r$ and $1\leq r\leq g-1$ instead of
$\varepsilon_{g-1}$, one obtains

  %To compute $L*(-1)^L$
%assume $\kappa=0$ for simplicity of notation. Then for $\beta =
%(\chi-1,...,\chi-1)$ we have $\delta_\beta = {}^p \delta_\beta
%\oplus T_\beta$
%$$ L*(-1)^L \oplus L*T_\beta = \delta_\alpha*\delta_\beta \ .$$
%Notice
%$$ \delta_\alpha*\delta_\beta =\bigoplus_{i=0}^r
%\delta_(\chi,...,\chi,\chi-1,..,\chi-1,1,..,1) \ ,$$ where the
%cardinality of the entries $\chi-1$  is equal to the cardinality
%of the entry $1$, which is denoted $i$ and varies from $0$ to $r$.
%It has no
%
%Now $\gamma_1\leq \alpha_1+\beta_1 \leq \chi +\alpha_1 -\alpha_r$.

\bigskip\noindent
\begin{Lemma} \label{spec}{\it If $C$ is not hyperelliptic, then for $r=0,..,g$ the perverse sheaves $\varepsilon_r$
and $\delta_{2,1,..,1}$ (both of degree $r$) and $\delta_{2g-3,1,..,1}$ (of degree $2g-3+r$),
and ${}^p\delta_{2g-3,r}$ and ${}^p\delta_{r,1}$ for $r\leq g-1$ are irreducible perverse
sheaves on $X$.} \end{Lemma}

\bigskip\noindent
\underbar{Remark}: Another point of view is the following. Recall
that $p_g: C^{(g)}\to X$ is an isomorphism outside
$\kappa-W_{g-2}$. The automorphism
$$ \theta: C^{(g-1)} \to C^{(g-1)} $$
maps $C^{(g-2)}$ to $\theta(C^{(g-2)})\subseteq C^{(g-1)}$, hence
$$ p_g^{-1}(\kappa - W_{g-2}) = \theta(C^{(g-2)}) \subseteq
C^{(g-1)} \subseteq C^{(g)} \ .$$ The Abel-Jacobi map is defined by the choice of an auxiliary
point $P_0\in C$. This point also defines $C^{(r-1)}\hookrightarrow C^{(r)}$ via $D\mapsto
D+P_0$. For $r=g-1$ notice, that $\theta(D+P_0)+P_0+D=K$. Hence $P_0+ \theta(C^{(g-2)}+P_0)\in
C^{(g)}$ maps to $\kappa-W_{g-2}$. By dimension reasons, this is the fiber
$p_{g}^{-1}(\kappa-W_{g-2})$, at least over the generic point of $\kappa-W_{g-2}$. Now, it is
enough to show, that over the algebraic closure $\overline\eta$ of the generic point $\eta$ of
$\kappa-W_{g-2}$ the fiber cohomology (*)
$$ H^2(p_{g}^{-1}(\eta), F_\varepsilon \vert_{p_{g}^{-1}(\eta)}) =
0 \ $$ vanishes. $\varepsilon_g$  is an irreducible perverse sheaf
and $$ \fbox{$ \varepsilon_g = j_{!*} (F_\varepsilon) $} \
$$ for the open inclusion $j: X\setminus (\kappa - W_{g-2})
\hookrightarrow X$, if and and if this holds. Argue as in section
\ref{RR} to conclude this:
 Since the fiber $p_{g}^{-1}(\overline\eta)$ is
a one-dimensional projective space over $\overline\eta$, the
vanishing statement (*) is equivalent to the fact, that the
restriction of the rank one coefficient system $F_\varepsilon$ on
$C^{(g)}$ to $\theta(C^{(g-2)}) \subseteq C^{(g-1)} \subseteq
C^{(g)}$ is nontrivial. Since $F_\varepsilon$ on $C^{(g)}$
restricts to the corresponding coefficient system on $C^{(g-1)}$,
this means that $\theta^*(F_\varepsilon)$ (for $F_\varepsilon$ on
$C^{(g-1)}$) remains nontrivial after restriction to $C^{(g-2)}$.

\chapter{Tensor categories}\label{Tann}

\section{Rigid tensor categories}
The category ${\cal BN}'$ of $BN$-sheaves of a Jacobian variety $X$ is a semisimple abelian
full subcategory ${\cal BN}'$ of the category of pure perverse sheaves on $X$, and it is
equivalent to the quotient  category ${\cal BN}$ contained in $\overline{Perv}(X)$ (see section
\ref{APP}).

\bigskip\noindent
\underbar{Tensor functor}: As explained in section \ref{APP} the convolution product of two
BN-sheaves is isomorphic to a direct sum of $BN$-sheaves and a constant perverse sheaf, which
become zero in the quotient category ${\cal BN}$. Hence the convolution
$$ K*L \ =\ Ra_*(K\boxtimes L)
$$ of objects $K$ and $L$ in ${\cal BN}$
defines an object in ${\cal BN}$. The category ${\cal BN}$ thus becomes a tensor category
$({\cal BN},*)$ with commutativity and associativity constraints inherited from the category
${\cal BN}'$, defined by direct images of constraints for the outer tensor product $\boxtimes$.
See section \ref{Begin}. The induced tensor functor will be written
$$ *: {\cal BN} \times {\cal BN} \to {\cal BN} \ .$$
According to section \ref{Begin} the associativity and commutativity constraints make ${\cal
BN}$ into a strict $\overline\Q_l$-linear tensor category in the sense of \cite{ML},
\cite{ML2}. The tensor functor $*$ has canonical extensions
$$*_{i\in I}\ \ {\cal BN}^I \to {\cal BN}\ .$$
The symmetric group $\Sigma_r$ acts on the $r$-fold tensor product of objects $A$ of a tensor
category. This action $$
 \sigma: \ A*\cdots *A \longrightarrow
 A*\cdots *A  \quad , \quad \sigma\in \Sigma_r $$ is induced from the action of the transpositions $(i,i+1)\in\Sigma_r$
 defined by
$$ \sigma_{i,i+1} = id_{A^{*i-1}} * \Psi_{A,A} * id_{A^{*r-i-1}}\ $$
for the commutativity constraint $\Psi_{A,A}:A*A\cong A*A$. Since $\sigma_{i,i+1}\circ
\sigma_{i,i+1}=id$ and since the $\sigma_{i,i+1}$ satisfy the Coxeter relations (\cite{CP}, p.
153), this extends to a well defined  action of the symmetric group $\Sigma_r$ (\cite{D}
section 7).

\bigskip\noindent
%$$ \xymatrix{ A\times \cdots \times A \ar[r]\ar[d]_\sigma & A*\cdots *A \ar@{..>}[d]^\sigma\cr
%A\times \cdots \times A \ar[r] & A*\cdots *A \cr} $$
\underbar{Schur Projectors}: The Schur projector attached to an irreducible representation
$\sigma_\alpha$ for some partition $\alpha$ of degree $r$ projects onto some
$\alpha$-idempotent $A^\alpha$ contained in $A^{*r}=A*\cdots *A$ ($r$ copies). This allows to
define objects $A^\alpha$ for  arbitrary $A$ in a $K$-linear tensor category, provided
$char(K)=0$ (\cite{D} 2.1). As a special case one obtains $\Lambda^r(A)=A^\alpha$, in case of
the partition $\alpha=([1]^r)$.

\bigskip\noindent
\begin{Lemma} {\it For an arbitrary partition $\alpha$ let $\alpha^*$ be the dual partition. Then
$\bigl(\delta_C\bigr)^\alpha$ and ${}^p\delta_{\alpha^*}$ are isomorphic in ${\cal BN}$.}
\end{Lemma}

\bigskip\noindent
\underbar{Proof}:  There are two different actions of the symmetric group $\Sigma_r$ on the
complex
$$ \delta_C * \delta_C * \cdots * \delta_C * \delta_C = Ra_{r,*}\Bigl(\boxtimes_{i=1}^r \delta_C\Bigr) \ $$
($r$ copies), where $a_r:X^r\to X$ is the addition map $a_r(x_1,..,x_r)=x_1+\cdots +x_r$. This
is a special case of $F_1*\cdots *F_r = Ra_{r,*}(\boxtimes_{i=1}^r pr_i^*(F_i))$.

\bigskip\noindent
The first action is induced from the geometric action of the symmetric group $\Sigma_r$ on
$C^r$. One has a commutative diagram
$$ \xymatrix@+0,5cm{ C^r \ar[d]_{p\times\cdots\times p} \ar[r]^{f} & C^{(r)\ar[d]^{p_r}} \cr
X^r \ar[r]^{a_r} &\ X \ . \cr} $$  This diagram is equivariant with respect to the geometric
permutation action $\sigma: C^r \to C^r$, and the perverse sheaf $\delta_{C^r}$ is a
$\Sigma_r$-equivariant sheaf on $C^r$ with respect to this action. Hence on $R(p_r\circ
f)_*\delta_{C^r}= Ra_{r,*}(R(p\times \cdots \times p)_*\delta_{C^r})$ there is an action,
induced from the action of $\Sigma_r$ on the equivariant complex $\boxtimes_{i=1}^r\ \delta_C$
on $X^r$
$$ \sigma^*(\boxtimes_{i=1}^r\ \delta_C) \ \overset{\sim}{\longrightarrow}\  \boxtimes_{i=1}^r\ \delta_C
\ $$ in the sense of \cite{KW}, definition 15.1 and remark 15.3.a) and d). This induced action
of $\Sigma_r$ on the complex $Ra_{r,*}\Bigl(\boxtimes_{i=1}^r \delta_C\Bigr)$ defined by
equivariance does not involve commutativity constraints! The objects $\delta_\alpha$ and
${}^p\delta_\alpha$ were defined in section \ref{ML} via this geometric action of $\Sigma_r$ on
$Ra_{r,*}\Bigl(\boxtimes_{i=1}^r \delta_C\Bigr) = \delta_C^{*r}$.

%\bigskip\noindent
%For $A=\delta_C\in {\cal BN}$ compare this perverse sheaf
%$A^\alpha$ with our previously defined complex
%$$ {}^p\delta_\alpha \hookrightarrow R(p_{r}\circ f)_* \delta_{C^{r}} \ .$$
%As complexes on $X$ they coincide. However, if we consider the action of the group $S_r$
%obtained by permutation, they differ.

\bigskip\noindent
Remark: For different complexes $F_i$ we could not imitate this, since we would only get maps
induced by $\sigma\in \Sigma_r$ of the form $ \boxtimes_{i=1}^r\ F_i \ \longrightarrow\
\boxtimes_{i=1}^r\ F_{\sigma^{-1}(i)}$.  $F_1=\cdots ... \cdots F_r=\delta_C$ was used to
obtain equivariance in this first construction.

\bigskip\noindent
Admitting the use of \lq{additional\rq}\ commutativity constraints the symmetric group acts on
$A^{*r}$ for an arbitrary object $A$ of a $K$-linear tensor category. In our case $A=\delta_C$
this gives almost the same action of $\Sigma_r$ on $\delta_C^{*r}$. Since $\delta_C=\overline
\Q_{l,C}[1]$  is a complex concentrated in the odd degree $-1$, this implies that the
commutativity constraints $\Psi_{K,L}$ for complexes give an additional sign factor
$sign(\sigma)$ compared to the first action. For this it is enough to look at the special case
of the involutions $\sigma_{i,i+1}$. The commutativity constraints
$\Psi_{A,A}=\Psi_{\delta_C,\delta_C}: Ra_*(\delta_C\boxtimes \delta_C) \to
Ra_*(\delta_C\boxtimes \delta_C) $ involved in the definition of $\sigma_{i,i+1}$ for
$A=\delta_C$ gives the additional factor $sign(\sigma_{i,i+1})=(-1)^{pq}=-1$ for $p=q=-1$. This
proves the claim.

\bigskip\noindent
\begin{Corollary}\label{aler} {\it The objects $\Lambda^i(\delta_C\bigr)$  of the category ${\cal BN}$ are isomorphic to
${}^p\delta_i$ for all $i\geq 0$. Hence $\Lambda^i(\delta_C\bigr)=0$ in ${\cal BN}$ if and only
if $i\geq 2g-1$.} \end{Corollary}

\bigskip\noindent
\begin{Corollary} {\it Suppose $C$ is not hyperelliptic. Then the dimension of $\delta_C\in {\cal BN}$ in the
sense of \cite{D}, section 7 is $$ dim(\delta_C) = 2g -2  \ .$$} \end{Corollary}

\bigskip\noindent
\underbar{Proof}: See \cite{D} page 165 and 167, preuve de 7.1 (iii)$\Longrightarrow$ (ii). In
fact this argument shows, that $\Lambda^i(A)=0$ implies $dim(A)\in \{0,..,i-1\}$. Conversely,
we use \cite{TW}, lemma 5.1

\bigskip\noindent
\begin{Lemma} {\it In rigid semisimple ribbon tensor categories simple objects have nonzero dimension.}
\end{Lemma}

\bigskip\noindent
Since the symmetric braiding of a $\overline\Q_l$-linear tensor
category is a special case of a rigid ribbon category, and since
our category is semisimple we can use \cite{D} (7.1.2)
$$ dim(\Lambda^i(A)) = dim(A)\cdot (dim(A)-1) \cdots
(dim(A)-i+1)/i! \ $$ to show $dim(\delta_C)$ is the smallest
positive integer $i$, such that $\Lambda^i(\delta_C)=0$. Here we
use, that all $\Lambda^i(\delta_C)$ are simple, if $C$ is not a
hyperelliptic curve. Also  recall $\Lambda^{2g-2}(\delta_C)=0$ and
$\Lambda^{2g-3}(\delta_C) \neq 0$.

\bigskip\noindent
We now summarize the properties of the category ${\cal BN}$, established so far, which imply
that ${\cal BN}$ is a $\overline \Q_l$-linear tensor category. Recall

\bigskip\noindent
\underbar{Unit element}: The skyscraper sheaf $\delta_0=1$ with support in zero  admits an
isomorphism $u_0: \delta_0 \cong \delta_0*\delta_0$, defined by $u=Ra_*(v)$, where $v:\delta_0
\cong \delta_0 \boxtimes \delta_0$ induces by the isomorphism $\overline \Q_l
\otimes_{\overline\Q_l} \overline\Q_l \to \overline\Q_l$ in the stalk over the origin $0\in
X\times X$. Obviously, the functor $u_K: K \mapsto \delta_0
* K$ is an auto-equivalence of ${\cal BN}$.
Hence $(\delta_0,u_0)$ is an identity object for the tensor category $({\cal BN},*)$. Notice
$End_{\cal BN}(1)=\overline\Q_l$.

\bigskip\noindent
\underbar{Internal Hom}:  Assume $\kappa$ normalized  to be zero. Furthermore ignore Tate
twists. Then the category ${\cal BN}'$ is stable under the duality $K\mapsto K^\vee$ by lemma
\ref{stable}. The contravariant functors $M \mapsto Hom_{{\cal BN}}(M*K,L) = Hom_{{\cal
BN}'}(M*K,L)$ are representable by the objects
$$ \underline{{\cal H}om}(K,L) = K^\vee * L $$
in ${\cal BN}$ defining internal Hom-objects. See section \ref{VT}, in particular corollary
\ref{zuvor},\ref{zuvor2} and \ref{Canonical}.

%To show this we have to find functorial isomorphisms $\psi_M$ commuting with the maps induced
%by morphisms $f:M'\to M$ in ${\cal BN}$
%$$ \xymatrix{ Hom_{\cal BN}(M*K,L) \ar[d]\ar[r]^-{\psi_M}
%& Hom_{\cal BN}(M,(DK)^-*L) \ar[d]\cr Hom_{\cal BN}(M'*K,L)
%\ar[r]^-{\psi_{M'}} & Hom_{\cal BN}(M',(DK)^-*L)\cr}
%$$
%%Since the abelian category ${\cal BN}$ is $\overline\Q_l$-linear
%%and semisimple, the functoriality will an immediate consequence of
%%the existence of $\psi_M$ by looking at the stalks at the origin
%%$0\in X$.
%Recall from \ref{APP} and \ref{VT}, that there exists a functorial
%isomorphism between the bi-covariant functors ${\cal BN}\times
%{\cal BN} \to Vec_{\overline \Q_l}$
%$$  Hom_{\cal BN}(D(M),N)
% \ \cong \ ker\Bigl(cosp: {\cal H}^0\bigl(M^-*N\bigr)_0
%\to {\cal H}^0\bigr(M^-*N \bigl)_\eta\Bigr) \ .$$  Hence there
%exist functorial isomorphisms
%$$ Hom(M*K,L) \ \cong \ ker\Bigl(cosp: {\cal H}^0\bigl(D(M*K)^-*L\bigr)_0 \to
%{\cal H}^0\bigr(D(M*K)^-*L \bigl)_\eta\Bigr) \ $$ and
%$$ Hom(M,DK^-*L) \ \cong \ ker\Bigl(cosp: {\cal H}^0\bigl(DM^-*DK^-*L\bigr)_0 \to
%{\cal H}^0\bigr(DM^-*DK^-*L \bigl)_\eta\Bigr) \ .$$ Since there
%are functorial isomorphisms $(K*L)^-\cong K^-*L^-$ and $(DK)^-
%\cong D(K^-)$, the existence of the isomorphisms $\psi_M$ now
%follows.

\bigskip\noindent
\underbar{Rigid tensor categories}: Then we define the dual object
$$K^\vee = (DK)^- = \underline{{\cal H}om}(K,1)$$ for $K\in {\cal BN}$.
Obviously there exists a functorial isomorphism $i_K:(K^{\vee})^\vee \cong K$, which shows that
every object of ${\cal BN}$ is reflexive. Now, since the category ${\cal BN}$ is semisimple,
the tensor category $({\cal BN},*)$ is a rigid tensor category by corollary \ref{sppo}, lemma
\ref{thereafter} and the remarks thereafter. Since the tensor functor is
$\overline\Q_l$-bilinear, it is a $\overline\Q_l$-linear tensor category (categorie tensorielle
sur $\overline\Q_l$ in the sense of \cite{D} 2.1). In addition it is semisimple and satisfies
the finiteness conditions of \cite{D} 2.12.1.

\bigskip\noindent
\underbar{Fiber functors}: A fiber functor from a
$\overline\Q_l$-linear abelian tensor category ${\cal C}$ with
values in a $\overline\Q_l$-algebra $R$ is an
$\overline\Q_l$-linear exact faithful tensor functor
 $$ \omega: {\cal C} \to Mod_R \ .$$
To be faithful, it is enough that $\omega(K)=0$ implies $K=0$ on objects, since ${\cal C}$ is
abelian. If ${\cal C}$ is semisimple, it is enough that $\omega$ is an additive functor, in
order to be an exact functor. There is a more general concept of fiber functor and Tannakian
duality in \cite{D}, theorem 8.17 replacing functors to $Mod_R$ by more general exact (then
faithful) functors to more general  $\overline\Q_l$-linear abelian tensor categories.

\bigskip\noindent
\section{The symmetric powers $\varepsilon_r$}\label{SP}

\bigskip\noindent
Recall $\varepsilon_1=\delta_C$. The direct sum of perverse sheaves $\varepsilon_r\in Perv(X)$
$$ \bigoplus_{r=0}^\infty \varepsilon_r $$
is equipped with  \lq{multiplication}\rq\ maps $m$
$$ \xymatrix@+0,5cm{  \varepsilon_{r-1} * \delta_C \ar[dr]_m\ar[r]^-\simeq &
\varepsilon_r \oplus \delta_{2,1^{r-1}} \ar[d]^{pr_1}\cr   &
 \varepsilon_r\cr}  $$
and \lq{differentiation\rq}\ maps $d$, using $\delta_C^\vee = \delta_{2g-3}$,
$$ \xymatrix@+0,5cm{ \varepsilon_r*\delta_C^\vee \ar@{^{(}->}[r] \ar[dr]_d & \varepsilon_r *\delta_{2g-3} \ \cong \ \delta_{2g-2,1^{r-1}} \oplus \delta_{2g-3,1^r} \ar[d]^D\cr  & \varepsilon_{r-1}\cr}
 $$
$m$ is defined by symmetrization (the Schur projector $pr_1$ to $\varepsilon_r$). $d$ is
defined by the multiplication law (section \ref{ML}) using $\delta_C^\vee \equiv \delta_{2g-3}$
in $K_*(X)$. The vertical map $D$ is the composite of  the projection $pr_1$ to
$\delta_{2g-2,1^{r-1}}$ followed by the projection map onto the perverse constituent
${}^p\delta_{2g-2,1^{r-1}} \cong \varepsilon_{r-1}$
$$\delta_{2g-2,1^{r-1}} = T \oplus {}^p\delta_{2g-2,1^{r-1}} \to
{}^p\delta_{2g-2,1^{r-1}} \cong \varepsilon_{r-1}\  $$ using $\varepsilon_{r-1} =
\delta_{0}*\varepsilon_{r-1} \cong \delta_{2g-2}*\varepsilon_r = {}^p\delta_{2g-2,1^{r-1}}
\oplus T'$, where as usual $T$ and $T'$ denote certain translates of constant perverse sheaves
on $X$.

\bigskip\noindent
\begin{Lemma} \label{superp}{\it For $m=(m_r)$ and $d=(d_r)$ the following holds
\begin{enumerate} \item $S^r(V)=H^\bullet(X,\varepsilon_r)$
for $V=H^\bullet(C)$ is the  super-symmetric  polynomial algebra $S^\bullet(V)$ on $V^\vee$
with multiplication induced by $m$ on $S^1(V)\times S^\bullet(V)$.
\item The map $d$ is super-differentiation
$$ d: S^r(V) \otimes V^\vee \to S^{r-1}(V) \ .$$
\end{enumerate}} \end{Lemma}

\bigskip\noindent
\underbar{Proof}: All projection maps are defined via idempotents of the symmetric group (Schur
projectors).  $m$ is a special case of the multiplication maps
$$ m_{a,b} : \varepsilon_a * \varepsilon_b \to \varepsilon_{a+b}
\ $$ defined by $ \varepsilon_a * \varepsilon_b \hookrightarrow \delta_C^{*(a+b)}\to
\varepsilon_{a+b}$ via the symmetrizing Schur projector, which projects to $\varepsilon_{a+b}$.
Now the first claim follows from the K\"unneth theorem
$$ \xymatrix{ H^\bullet(X,\varepsilon_a*\varepsilon_b)\ar[d] \ar@{=}[r] & S^a(V)\otimes S^b(V)\ar[d] \cr
H^\bullet(X,\varepsilon_{a+b}) \ar@{=}[r] & S^{a+b}(V)\ . \cr} $$ By abuse of notation $m$ will
now also denote the  maps induced on the cohomology groups
$$ m_{a,b} : S^a(V) * S^b(V) \to S^{a+b}(V)
\ .$$

\bigskip\noindent
\underbar{Concerning $d$}: First show, that $d$ defines a super-derivation of the
super-polynomial algebra $S^\bullet(V)$. This boils down to an identity in the group ring of a
symmetric group. Recall $V^{(2g-3)}\otimes V^{(1^a)} = \Lambda^{2g-3}S^{a} (V^{\otimes (2g-3 +
a)})$ in the tensor category of super vectorspaces, where $\Lambda^{2g-3}$ denotes alternation
over the first $2g-3$ tensor factors with the symmetric group $\Sigma_{2g-3}$, and $S^{a}$
denotes symmetrization over the last $a$ tensor factors for the corresponding group $\Sigma_a$.
Up to some constant the Schur projector
$$ \Pi_a: V^{(1^{2g-3})}\otimes V^{(a)} \to V^{(a-1,1^{2g-2})}   $$
is defined  by $ S^{a-1} \Lambda^{2g-2} S^{a}\Lambda^{2g-3}$ in the group ring
$\overline\Q_l[\Sigma_{2g-3+a}]$. Here $S^{a-1}$ is the symmetrizer over the last $a-1$
elements, and $\Lambda^{2g-2}$ is the alternator over the first $2g-2$ elements. $d$ arises
from the composite map
$$ V^{(1^{2g-3})}\otimes V^{(a)} \to V^{(a-1,1^{2g-2})} \to
V^{(a-1)} \ .$$ Ignoring a constant, the element $ S^{a-1} \Lambda^{2g-2} S^{a}\Lambda^{2g-3}$
in the group ring becomes
$$ S^{a-1} \Lambda^{2g-2} \Bigl(\sum_{i=1}^{a} \Sigma_{a-1}\cdot \sigma_{(2g-2)(2g-3+i)}\Bigr)
\Lambda^{2g-3}
$$
$$ = (a-1)!\cdot S^{a-1} \Lambda^{a+1}
\sum_{i=1}^{2g-2} \sigma_{(2g-2)(2g-3+i)} \Lambda^a \ .$$ Here
$\Sigma_{a-1}$ is is the group of permutation of the last $a-1$
tensor factors, and $\sigma_{ij}$ are transpositions of the places
$i$ and $j$. Notice $\Lambda^{2g-2}
  \Sigma_{a-1} = \Sigma_{a-1} \Lambda^{2g-2}$.
From this it easily follows, that
$$ \Pi_{a+b}\circ (id \times m_{a,b}) = m_{a-1,b}\circ \Pi_a
\ +\ m_{a,b-1}\circ \Pi_b \circ (id\times c_{a,b})  $$ as maps $ V^{(1^{2g-3})}\otimes V^{(a)}
\otimes V^{(b)} \to V^{(a+b-1)}$, where
$$c_{a,b}:
 V^{\otimes a} \otimes
V^{\otimes b}
 \to
 V^{\otimes b} \otimes V^{\otimes a}
  $$ is the permutation of the last blocks.
This finishes the sketch of the proof for the super derivation property, which now easily
follows.

\bigskip\noindent
\underbar{Second step}: It remains to compute $d$ on the generators $V=S^1(V) \subseteq
S^\bullet(V)$, i.e. we have to show that the map
$$  d: S^1(V)*V^\vee \to S^0(V) \to k \ .$$
defines a nondegenerate pairing $V\times V^\vee \to \overline\Q_l$. $d$ is the map induced on
cohomology by the underlying projections
$$ \delta_C* \delta_C^\vee \to \delta_{2g-2} \to \delta_0 \ .$$
Since $\delta_C\in Perv(X)$ is irreducible,
$$ Hom_{Perv(X)}(\delta_C*\delta_C^\vee,\delta_0) \cong Hom_{{\cal BN}}(\delta_C*\delta_C^\vee,\delta_0) $$
has dimension one by section \ref{rig}. This defines $d$  uniquely up to a constant. Recall $
(\delta_C)^\vee = D(\delta_{-C}) \cong \delta_{-C}$. The induced map on cohomology defines a
nondegenerate pairing on the cohomology groups by theorem \ref{th1}: For the closed smooth
subvariety $C\hookrightarrow X$ the map $f:C\times C\to X$ defined by $f(x,y)=x-y$ blows down
the diagonal $\Delta_C\hookrightarrow C\times C$ to the point zero: $f^{-1}(0)=\Delta_C$. Since
$\delta_C*\delta_{-C}=Rf_*(\delta_C\boxtimes \delta_C)$, we can compute $\delta_C*\delta_{-C}
\to \delta_0$ using the proper basechange diagram
$$ \xymatrix{ C\times C \ar[r]^-f & X \cr \Delta_C\ar[u]^i \ar[r] & \{0\} \ar[u]_{i_0}\cr} $$
On the cohomology groups this induces the nondegenerate cup-product pairing. See \cite{SGA5},
p.101, \cite{SGA45}, cycle, section 2.3 and also \cite{FK}, p.155 or \cite{Mi2}, theorem 12.3
$$ \xymatrix{ H^\bullet(C)\times H^\bullet(C)\ar[d]_{i^*}\ar[r]  & H^\bullet(X,\delta_0)\ar@{=}[d] \cr H^\bullet (\Delta_C) \ar[r]^{S} & \overline\Q_l \cr} $$
for the trace map $S: H^2(\Delta_C) \cong \overline\Q_l$ as explained in section \ref{rig}.

\goodbreak
\section{$0$-Types in the symmetric algebra}\label{Nulltypes}

%\bigskip\noindent
%Recall, that a perverse subsheaf is a direct summand by purity. In general, if
%$$ f: L \to L' $$
%is a homomorphisms between pure perverse sheaves, then $f$ factorizes over the projection to a
%summand of $L$, an isomorphism and the inclusion of a summand of $L'$. Since for a BN-sheaf its
%cohomology is zero if and only if the perverse sheaf is zero itself, we get

By lemma \ref{Ffaith} the cohomology functor $H^\bullet(X,-)$ is exact and faithful on the
category ${\cal BN}'$. For $r>0$ assume, that $ L\hookrightarrow \varepsilon_r$ is a nontrivial
irreducible perverse constituent. Suppose $L$ is of $0$-type. Then by definition
$$ H^\bullet(X,L) = H^0(X,L) \ .$$
%We already know $L=0$, if and only if $H^0(L)=0$. Notice $L^\vee$
%and $L*L^\vee$ are again of $0$-type. Hence
%$$ Hom_{Perv{X}}(\Omega, L*L^\vee) = 0 \ ,$$
%again since the cohomology functor $H^\bullet(X,-)$ is fully
%faithful by Gabber's theorem. By \ref{APP} the sheaf perverse
%$K=L*\delta_C^\vee$ is nonzero. By the multiplicity formula given
%in \ref{MF} therefore
For the convolution $K=L*\delta_C^\vee$ corollary \ref{Mul} and $n_X(K)=0$ (lemma \ref{tri})
imply
$$ {}^p m(K)={}^p m(L)=1 \quad , \quad K=L*\delta_C^\vee \ .$$
Hence $K$ is irreducible. One has the following  morphisms
$$ K=L*\delta_C^\vee \ \overset{i}{\hookrightarrow}\ \varepsilon_r *\delta_C^\vee
\ \cong\ \delta_{(\chi-1,1^r)} \oplus \varepsilon_{r-1} \ \overset{\pi=pr_2}{\longrightarrow}\
\varepsilon_{r-1} \ .$$ The composite $ \pi\circ i: K=L*\delta_C^\vee \to \varepsilon_{r-1} $
is nonzero. It is enough to check this on cohomology. Since $H^\bullet(X,L)\neq 0$, it is
enough to show that the super-differentiation $d: H^\bullet(X,L)\otimes V^\vee \to S^{r-1}(V)$
is nontrivial. However this is obvious, since a super polynomial $P\in H^\bullet(X,L) \subseteq
S^r(V)$ with vanishing partial super derivatives  is zero. Since $K$ is simple, we obtain

\begin{Lemma}
$ \pi\circ i: \ K=L*\delta_C^\vee\ \hookrightarrow\ \varepsilon_{r-1} $ is an inclusion.
\end{Lemma}

\bigskip\noindent
Now we revert things. Consider
$$K*\delta_C = L*\delta_C^\vee*\delta_C = L\ \oplus\ (L*\Omega) \ $$
and the morphisms
$$ K*\delta_C\ \overset{i'}{\hookrightarrow}\ \varepsilon_{r-1} *\delta_C \ =\
\varepsilon_r \oplus \delta_{(2,1^{r-2})}\ \overset{\pi'}{\longrightarrow}\ \varepsilon_r \ .$$
Lemma \ref{stable} implies $L*\delta_C\neq 0$. Looking at cohomology  the composed map
$\pi'\circ i'$ is nonzero, since $\pi': \varepsilon_{r-1}*\delta_C \to \varepsilon_r$ induces
the multiplication $S^{r-1}(V)\otimes V\to S^r(V)$ (lemma \ref{superp}, part 1). In fact, the
multiplication map $W\otimes V \to S^r(V)$ is nonzero for any nonzero linear subspace
$W\subseteq S^{r-1}(V)$; that $S^\bullet(V)$ has zero divisors fortunately is irrelevant for
this. Therefore $\pi'\circ i'$
$$ \pi'\circ i' :\  K*\delta_C\ =\ L\ \oplus\ (L*\Omega) \
\longrightarrow \ \varepsilon_r \ $$ is nonzero. Concerning $\Omega$ we refer to section
\ref{CundA} and \ref{Hyperelliptic}.

\bigskip\noindent
\begin{Lemma}\label{lastl} $\pi'\circ i'(L*\Omega)=0\ $ for $\ L*\Omega \subseteq K*\delta_C$ (nonhyperelliptic case)
and $\pi'\circ i'(L*\delta_{C-C})=0\ $ (hyperelliptic case). \end{Lemma}

\bigskip\noindent
\underbar{Proof}: Suppose $f: L*\Omega \to \varepsilon_r $ were nontrivial. Recall $Hom_{\cal
BN}(L*\Omega,\varepsilon_r) = Hom_{\cal BN}( \Omega, L^\vee * \varepsilon_r) \ $ via
$\varphi\mapsto \Phi$, using
$$ \xymatrix{ \Omega \ar[r]_-{\delta*id}\ar@/^7mm/[rr]^-\Phi & L^\vee*L*\Omega \ar[r]_-{id*\varphi} & L^\vee *\varepsilon_r\cr}  \ .$$
Therefore Gabber's theorem gives a contradiction in the non-hyperelliptic case, since
$H^{-2}(X,\Omega)\neq 0$ and $\Omega$ is irreducible, whereas $L^\vee*\varepsilon_r$ is of
$1$-type: i.e. $H^\nu(X,L^\vee*\varepsilon_r)=0$ holds for $\vert \nu\vert
>1$ (use the K\"{u}nneth formula, and that $L^\vee$ is of $0$-type
and $\varepsilon_r$ is of $1$-type). In the hyperelliptic case $\Omega=A\oplus \delta_{C-C}$
for $A\cong T_e^*(\varepsilon_2)$, and the same argument as above can be applied for the
summand $\delta_{C-C}$. This proves the lemma.

\bigskip\noindent
The morphism $ \pi'\circ i': \ K*\delta_C \longrightarrow \varepsilon_r $
 is nonzero and has a nontrivial image
$$ \fbox{$ L'=\pi'\circ i'\Bigl(K*\delta_C\Bigr)  \hookrightarrow \varepsilon_r $} \ .$$
Since under the cohomology functor $H^\bullet(X,-)$ the maps $\pi\circ i$ and $\pi'\circ i'$
induce differentiation respectively the multiplication on the super-symmetric algebra, we get
commutative diagrams
$$ \xymatrix@+0,5cm{ S^{r-1}(V)  &  S^r(V)\otimes V^\vee\ar@{->>}[l]_d  \cr
H^\bullet(X,K)\ar@{^{(}->}[u] & H^0(X,L)\otimes V^\vee = H^{\bullet}(X,L)\otimes V^\vee
\ar@{^{(}->}[u]\ar[l]_-{\varphi} \cr}
$$
and
$$ \xymatrix@+0,5cm{ S^{r-1}(V)\otimes V \ar[r]^m  & S^r(V) \cr
H^\bullet(X,K)\otimes V\ar@{^{(}->}[u]\ar[r]^-\psi  & H^{\bullet}(X,L')\ . \ar@{^{(}->}[u]\cr}
$$
\begin{Lemma}\label{added}  $\psi$ is a surjection and $\varphi$ an isomorphism of graded vector spaces.
\end{Lemma}

\bigskip\noindent
\underbar{Proof}: This follows from $\pi\circ i(L)=K$ and $\pi'\circ i'(K)=L' $, since
$H^\bullet(X,\pi\circ i)=d$ and $H^\bullet(X,\pi'\circ i')=m$  (adding the extra observation,
that for all sheaves considered in this section $\delta_{\alpha}={}^p\delta_\alpha$ holds by
$\alpha_1 \leq 2\leq g-1$ due to our general assumption $g\geq 3$). This completes the proof of
the claim.

\bigskip\noindent
Lemma \ref{lastl} implies
$$  L \ \cong \ L' $$
in the non-hyperelliptic case. In the hyperelliptic case $L'$ is a perverse quotient of
$L\oplus (L*A)$ . Since $L$ is of $0$-type and $A$ is of $1$-type, in any case $L'$ is of
$1$-type. Recall $L\neq 0$ implies $H^\bullet(X,L)\neq 0$.  Since $L$ is of $0$-type we always
have
\begin{enumerate}
\item $H^\bullet(X,L)=H^0(X,L)\cong \overline\Q_l^d $ for some $d\geq 1$.
\item $ H^\bullet(X,K)= H^\bullet(X,L)\otimes H^\bullet(C)$ has
dimensions $d, 2dg,d$ in the degrees $1,0,-1$ and and vanishes in all other degrees. \item $L'$
is of 1-type ($0$-type in the non-hyperelliptic case) and contains $L$.
\end{enumerate}

\bigskip\noindent
Notice $H^\bullet(X,L)\subseteq H^\bullet(X,L')$, which is a consequence of lemma \ref{added}
and the Euler formula for superpolynomials $Q$, suffices to show the last assertion $L\subseteq
L'$.

\bigskip\noindent  Now suppose
$Q\in H^\bullet(X,L)=H^0(X,L)$ such that
$$Q=P + \eta_+\cdot \eta_{-}\cdot \tilde P \quad , \quad P\neq 0 $$ with polynomials $P\in S^r(H_-)$ and
$\tilde P\in S^{r-2}(H_-)$ in the usual sense. Here $\eta_{\pm }$ denote the generators of two
odd parts $H^{\pm 1}(C)\subseteq S^1(H_+)$.

\bigskip\noindent
{\it The non-hyperelliptic case}. Since $r>0$, there exists a partial derivative $\partial_i$
for some $i=1,..2g$ with respect to the even variables in $H_-$, such that $\partial_i P \neq
0$. Since $\partial _i Q\neq 0$ is in the image of $H^0(X,L)\otimes V^\vee$ (under
differentiation), it is contained in $H^0(X,K)$ by lemma \ref{added}. Therefore
$$ \eta_1 \partial_i(Q) = \partial_i (\eta_1 Q) = \partial_i(\eta_1 P) = \eta_1 \partial_i P
\neq 0 $$ is contained in $ H^1(X,\varepsilon_r) \cap H^\bullet(X,L')$, once more by lemma
\ref{added}. In the non-hyperelliptic case this intersection is zero, since $H^1(X,L')\cong
H^1(X,L)=0$ . This contradicts our assumption $P\neq 0$.

\bigskip\noindent
{\it The hyperelliptic case}. For simplicity assume $d= dim(H^0(X,L))=1$. Then $Q$ necessarily
generates $H^0(X,L)$.  The odd derivatives of $Q$ generate $H^{\pm 1}(X,K)$ by lemma
\ref{added}. Hence  $H^{\pm 1}(X,K)= \eta_{\pm}\tilde P \overline\Q_l$. The multiples
$x_j\eta_{\pm}\tilde P$ obtained from the even generators $x_j$ of $H_-\subseteq V$ define $2g$
linear independent elements in $H^{\pm 1}(X,L')$. Hence
$$ 2g \leq dim(H^{\pm 1}(X,L')) \ .$$ Since $L'$ is a quotient of $L\oplus (L*A)$, the cohomology
$H^{\pm 1}(X,L')$ is a quotient of  the $2g$-dimensional space $H^{\pm 1}(X,L*A)\cong  H^{\pm
1}(X,A)$. Thus
$$ dim(H^{\pm 1}(X,L')) \ \leq\  2g \ .$$ Therefore the kernel $B$ of
$\pi'\circ i': L*A \to \varepsilon_r$ is of $0$-type. But $B\hookrightarrow L*A$ implies
$Hom_{\cal BN}(B,L*A) \cong Hom_{\cal BN}(A^\vee,B^\vee*L) \neq 0$. Since $B^\vee*L$ is of
$0$-type and $A$ is irreducible of $1$-type,  Gabber's theorem gives a contradiction  unless
$B=0$. Hence $L*A \hookrightarrow L'$. By the same argument $L*A\subseteq L'$ and $L\subseteq
L'$  do not intersect
$$  L \oplus (L*A) \overset{\sim}{\longrightarrow} L' \ .$$
Hence $$ H^{\pm 1}(X,L') = V\cdot  \eta_{\pm} \tilde P$$ and  $ dim(H^0(X,L')) = 1 +
dim(H^0(X,A)) = 2 + dim(S^2(\overline\Q_l^{2g}))$. This dimension formula will give a
contradiction: $H^0(X,K)$ is of dimension $2g$, generated by the $2g$ even derivatives of $Q$
$$ < {\partial}_i Q > = H^0(X,K) \ .$$  Since $\eta_\pm H^0(X,K) \subseteq  H^{\pm 1}(X,L') = V\cdot \eta_\pm \tilde P$,
$ \eta_\pm {\partial}_i Q =  \eta_\pm {\partial}_i P = \eta_{\pm} \lambda_i(x) \tilde P(x) $
holds for linear functions $\lambda_i(x)\in V$. Thus $
\partial_i P(x) = \lambda_i(x) \tilde P(x)$, and the Euler formula implies $ P(x) = q(x) \cdot
\tilde P(x)$ for the quadratic polynomial $q(x)= r^{-1}\sum_{i=1}^{2g} x_i \lambda_i(x) $. Thus
$ q(x) {\partial}_i \tilde P(x) = \tilde\lambda_i(x) \tilde P(x)$ for $\tilde\lambda_i(x) =
\lambda_i(x) - {\partial}_i q(x)$, hence
$$ \frac{ d\tilde P(x)}{\tilde P(x)} = \frac{ \sum_{i=1}^{2g} \tilde\lambda_i(x) dx_i}{q(x)}
\ .$$ By decomposing $\tilde P(x)$ into a product of irreducible polynomials, this forces
$$ \tilde P(x) = c \cdot q(x)^n $$
for some constant $c\in \overline\Q_l$ and an integer $n$, if $q(x)$ is irreducible. Or $\tilde
P = c\cdot \lambda(x)^n\tilde\lambda(x)^m $ for the linear factors $\lambda(x)\tilde
\lambda(x)=q(x)$ in case $q(x)$ is reducible. Hence
$$ Q(x) = q(x)^{n+1} + c\cdot q(x)^n\eta_+\eta_- $$
or $Q(x)= \lambda^{n+1}\tilde\lambda^{m+1} + c\cdot \lambda^n\tilde\lambda^m \eta_+\eta_- $. In
this second case the partial derivatives $\partial_i Q$ only span a space of dimension $\leq
2$. This contradicts $dim(H^0(X,K))=2g>2$. Thus $q(x)$ is irreducible. By lemma \ref{added} the
terms $x_j {\partial}_i Q(x)$ and $\eta_+\eta_-\tilde P$ generate $H^0(X,L')$. Since for
varying $i$ and $j$ the span of the
$$x_j {\partial}_i Q(x)  = \Bigl( x_j\cdot {\partial}_i q(x)\Bigr)  \cdot R(x) \quad , \quad R(x)= (n+1)q(x)^n + cn q(x)^{n-1} \eta_+\eta_- $$
has dimension at most $dim(S^2(\overline \Q_l^{2g}))$, this implies
$$ dim\bigl(H^0(X,L')\bigr) \leq 1 + dim\bigl(S^2(\overline \Q_l^{2g})\bigr) \ $$
contradicting the dimension formula above. At this stage we can drop the assumption, that $L$
is irreducible, and get

\bigskip\noindent
\begin{Theorem} \label{0}{\it Suppose $L$ is a perverse sheaf  of $0$-type contained in $\varepsilon_r$ for $r>0$.
Suppose $dim(H^0(X,L))=1$ in case $C$ is hyperelliptic. Then the cohomology of $L$ satisfies
$$ H^\bullet(X,L)\ \ \subseteq\ \ \eta_1\eta_{-1}\cdot S^{r-2}(H_-)\ \subseteq\ S^r(V)_{even} \
.$$} \end{Theorem}

\bigskip\noindent
\begin{Corollary} {\it The natural map induced by the inclusion $i$ composed with the multiplication $m_{r,r}$
$$  L * L\ \overset{i*i}{\hookrightarrow}\ \varepsilon_r * \varepsilon_r
\overset{m_{r,r}}{\longrightarrow}\ \varepsilon_{2r}
$$
is the zero map (in the situation of theorem \ref{0}).} \end{Corollary}

\bigskip\noindent
\underbar{Proof}: Since the cohomology functor is faithful on ${\cal BN}'$, it suffices that
the induced map on the cohomology groups is zero. Since $m_{r,r}$ induces the super
multiplication
$$ m_{r,r}: S^r(V)\otimes S^r(V) \to S^{2r}(V) \ $$
(see section \ref{SP}), this map is zero on $H^\bullet(X,L)\otimes H^\bullet(X,L) \subseteq
S^r(V)\otimes S^r(V) $ by the nilpotency
$$(\eta_1\eta_{-1})^2=0\ .$$ This proves the corollary.

\bigskip\noindent
\section{Structure theorem}
Let $K$ be an algebraically closed field of characteristic zero (in our case
$K=\overline\Q_l$). Let ${\cal T}$ be a $K$-linear tensor category and finitely
$\otimes$-generated (in our case ${\cal BN}$ is generated by $\delta_C$). Then

\bigskip\noindent
\begin{Theorem} \label{Deligne}{(\cite{De}) For ${\cal T}$ the following assertions are equivalent \begin{enumerate} \item $\cal T$ is
equivalent as a tensor category to the category of representation $Rep(G,\varepsilon)$ of a
supergroup $G$. \item Every object $A$ of ${\cal T}$ is annihilated by a Schur projector
(depending on $A$).
\end{enumerate}} \end{Theorem}

\bigskip\noindent

\bigskip\noindent
This result also implies, that the following conditions for an object $A$ of ${\cal T}$ are
equivalent (see \cite{De}):
\begin{enumerate}
\item $A$ is annihilated by some Schur functor $A^\alpha=0$.
\item There exists an integer $N$ such that $$lenght(A^{\otimes n}) \leq N^n$$ holds
for all $n\geq 0$.
\end{enumerate}

\bigskip\noindent
Hence $A^{\otimes m}$ is Schur finite provided $A$ is Schur finite, since $ lenght((A^{\otimes
m})^{\otimes n}) = lenght(A^{\otimes mn})\leq N^{mn} = (N_1)^n$ for $N_1=N^m$.

\bigskip\noindent
For the category ${\cal BN}$ this has the following consequence: Since $\delta_C$ is a
generator for ${\cal BN}$ and since $\Lambda^{2g-1}(\delta_C)=0$, every object of ${\cal BN}$
is Schur finite. Hence by \cite{De}

\bigskip\noindent
\begin{Corollary} {\it The category ${\cal BN}$ is a semisimple super-Tannakian category equivalent to the
category of representations $Rep(G,\varepsilon)$ of some supergroup $G$ and some datum $\mu_2
\to G$.} \end{Corollary}

%\bigskip\noindent
%By the same reason the category ${\cal T}(\chi)$ is semisimple
%super-Tannakian, and equivalent to the category of representations
%$Rep(G(\chi),\varepsilon)$ of a supergroup $G(\chi)$. The
%Grothendieck ring of this category $Rep(G(\chi),\varepsilon)$ is
%the Grothendieck group of the category $Rep(Gl(2g-2))$. Being
%semisimple is therefore one of the braided categories described in
%\cite{WF}.

\section{The Tannaka category ${\cal BN}$}

\bigskip\noindent
In this section we suppose, that the Riemann constant $\kappa$ is normalized to be zero. For
the tensor category $ {\cal BN} $ we know \begin{enumerate} \item  ${\cal BN}$ is generated by
a simple object ${\cal X}=\delta_C$. \item $\Lambda^{\chi+1}({\cal X})=0$ for some $\chi>0$.
\end{enumerate} This implies ${\cal BN}$ to be super-Tannakian. Let $$\omega: {\cal BN} \to
Vec_{\overline\Q_l}^{\pm} $$ be a super fiber functor. Since $\omega$ is an exact tensor
functor $\omega({\cal X})^{\otimes a} \cong \omega({\cal X}^{\otimes a})$, and it commutes with
the projectors $\Lambda^{a}$. Hence
$$ \Lambda^{a}\Bigl(\omega({\cal X})\Bigr) \cong \omega\Bigl(\Lambda^{a}({\cal X})\Bigr) \ .$$
For $\omega({\cal X})= V^+ \oplus V^-$ this implies  $\Lambda^{a}(V^+\oplus V^-)=0$ for the
case $a=\chi +1$. For $S^a(V^-) \subseteq \Lambda^a(V^+\oplus V^-)$ hence $S^a(V^-)=0$,
therefore $V^-=0$. Since $\omega$ is a tensor functor this implies $\omega({\cal X}^{\otimes
n})= (V^+)^{\otimes n}$ is even for all $n$. The same holds for all subquotients. Since ${\cal
X}$ is a tensor generator of the tensor category ${\cal BN}$, this implies that
$$ \omega : {\cal BN} \to Vec_{\overline\Q_l} $$
is a fiber functor with values in the category of finite dimension $\overline\Q_l$-vector
spaces.

\bigskip\noindent
\begin{Corollary} {\it The category ${\cal BN}$ is Tannakian, and isomorphic to the category $Rep(H)$ of finite
dimensional $\overline\Q_l$ representations of an algebraic group $H$ over $\overline\Q_l$. The
Zariski connected component $H^0$ of $H$ is reductive.} \end{Corollary}

\bigskip\noindent
\underbar{Proof}: This follows from \cite{DM}, theorem 2.11, prop.
2.20 (b), prop 2.21 (b) and remark 2.28.

\bigskip\noindent
\begin{Corollary}\label{lastc} {\it Suppose $L=\delta_{0}\hookrightarrow \varepsilon_r$, then $r=0$.} \end{Corollary}

\bigskip\noindent
\underbar{Proof}: By the main result of section \ref{Nulltypes} the natural projection
$\varepsilon_r * \varepsilon_r \to \varepsilon_{2r} $ is the zero map on the $0$-type
$L=\delta_0$ for $r>0$. By the existence of a fiber functor $\omega$ this projection induces
the multiplication on the symmetric polynomial algebra
$$ \overline\Q_l[W] = \omega\Bigl(\bigoplus_{r=0}^\infty
\varepsilon_r\Bigr) $$ where $W=\omega(\delta_C)$ is the canonical representation space of $G$
of dimension $2g-2$ defined by the fiber functor $\omega$. Since now $\overline\Q_l[W]$ is a
polynomial ring in the usual sense, hence an integral domain, this implies that
$\omega(L)\subseteq \overline\Q_l[W]$ is zero. Since $\omega$ is faithful, $\omega(L)=0$
implies $L=0$. A contradiction for $r>0$.

\bigskip\noindent
\begin{Corollary}\label{inv} {\it $\overline\Q_l[W]^G=0$.} \end{Corollary}

\bigskip\noindent
\underbar{Proof}: Any invariant subspace arises from some embedding
$\delta_0\hookrightarrow \varepsilon_r$. Since $\delta_0$ is of
$0$-type, the last corollary \ref{lastc} proves the claim.

\bigskip\noindent

\section{Representations}\label{Rep}
Let $G$ be an algebraic group over $\overline\Q_l$ or more generally over an algebraically
closed base field of characteristic zero. In fact, one can choose an isomorphism $\tau:
\overline\Q_l \cong \C$. Hence it will suffice to consider the field $\C$ of complex numbers in
the following. Let $W$ be an irreducible faithful representation of $G$ over the field $\C$ of
dimension $dim(W)>1$, such that
\begin{enumerate} \item $G^0$ is reductive.
\item
$W\otimes W$ has at most three irreducible constituents. If it has three irreducible
constituents, then $\Lambda^2(W)$ contains the trivial representation, but $\Lambda^2(W)$ is
not trivial.
\item $\Lambda^3(W)$ does not
contain the trivial representation. If $W\otimes W$ has three irreducible constituents, then
$W\otimes W\otimes W$ has (at most) 7 irreducible constituents.
\item
$\Lambda^i(W)$ is the trivial representation for $i=dim(W)$.
\item $\C[W]^G = \C$.
\end{enumerate}
These conditions come from 1) the semisimplicity of ${\cal BN}$, 2) corollary \ref{MA} and
section \ref{Hyperelliptic} concerning the decomposition of $\delta_C*\delta_C$, 3)  corollary
\ref{MA} and corollary \ref{triple} and for $g=3$ also lemma \ref{spec} describing the
decomposition of $\delta_C*\delta_C*\delta_C$. The condition on $\Lambda^3(W)$ follows from
lemma \ref{hy} and corollary \ref{MA}. Furthermore 4) follows from corollary \ref{aler} and
finally 5) from corollary \ref{inv}.

\bigskip\noindent

\goodbreak
\bigskip\noindent
\underbar{Remark}: Under these assumptions on $(G,W)$ we now prove step by step the following
assertions
\begin{enumerate}
\item Suppose a semidirect product $(G_1\times\cdots \times G_r)\rtimes \Delta$,
with a finite group $\Delta$ normalizing each factor $G_i$, surjects to the group $G$. Then $W$
defines an irreducible representation of this semidirect product, which is induced from a
tensor product of irreducible representations $(W_i,\lambda_i)$ of $G_i$ for $I=1,..,r$ such
that at most one of the representations $\lambda_i$ has dimension $>1$.
\item $G^0$ is semisimple.
\item $G^0$ is a product of simple groups of the same
Lie-type. \item $G^0$ is simple.
\item $G^0$ is either $Sl(2g-2)$ or $Sp(2g-2)$.
\item $G= G^0\times \pi_0(G)$ with a cyclic group $\pi_0(G)$
\end{enumerate}

\bigskip\noindent
We start with

\begin{Lemma}\label{43}
For projective algebraic representations $(W,\rho)$ of an algebraic group $G$ the following
holds
\begin{enumerate}
\item If $\rho=\rho_1\otimes \rho_2$ for projective representations $\rho_1,\rho_2$ the multiplicities
satisfy
$$ m_G(\rho\otimes \rho) \geq m_{G}(\rho_1\otimes\rho_1) \cdot m_{G}(\rho_2\otimes\rho_2) \
.$$ \item $m_G(\rho\otimes\rho)\geq 2$ for $dim(\rho)>1$. \item $m_{G'}(Ind_G^{G'}(\rho)
\otimes Ind_G^{G'}(\rho)) \geq m_G(\rho\otimes \rho) +1$ if $G \subsetneqq G'$ is a subgroup of
finite index.
\end{enumerate}
\end{Lemma}

\bigskip\noindent
\underbar{Proof}: The proof of 1  is obvious. For two notice $W\otimes
W=S^2(W)\oplus\Lambda^2(W)$. Concerning 3. observe that for
$$ W = Ind_{G}^{G'}( \rho)  = \oplus_{\delta\in
G'/G} \delta(V) $$ the tensor square representation is $$ W\otimes W = \bigoplus_{\delta\in
G'/G} \delta(V)\otimes \delta(V) \ \oplus \underset{\delta_1\neq \delta_2\in
G'/G}{\bigoplus\bigoplus} \delta_1(V)\otimes \delta_2(V) \ .$$ Since $\delta(V)\otimes
\delta(V) = \delta(V\otimes V)$, the first summand is $Ind_G^{G'}(\rho\otimes \rho)$. Obviously
$m_{G'}(Ind_G^{G'}(\rho\otimes \rho)) \geq m_{G}(\rho\otimes \rho)$ and
$m_G(\underset{\delta_1\neq \delta_2\in G'/G}{\bigoplus\bigoplus} \delta_1(V)\otimes
\delta_2(V))\geq 1$ for $G'\neq G$.

%Now $Y\otimes Y = (\lambda_1\otimes \lambda_1)\otimes \cdots \otimes (\lambda_r\otimes
%\lambda_r) $ and each
%$$\lambda_r\otimes \lambda_r = S^2(\lambda_r)\oplus \Lambda^2(\lambda_r)\ .$$ Unless
%$\lambda_i$ is one dimensional, this gives a splitting in two nontrivial factors. Hence, if
%more than one $\lambda_i$ have dimension $>1$, the multiplicity of $W\otimes W$ is $\geq 4$
%contradicting our assumptions. (Furthermore the multiplicity of $W\otimes W$ is $\geq 3$, if
%$\Delta\neq \Delta_Y$).

\begin{Corollary}\label{30}
For $W=Ind_{G}^{G'}(\rho)$, $dim(\rho)>1$ and $1< [G':G]<\infty$ suppose $W\otimes W$ has two
irreducible constituents or three, and in the latter case it case it contains the trivial
representation. Then $G'=G$ and $dim(\rho)=2$.
\end{Corollary}

\bigskip\noindent
\underbar{Proof}: Unless $G'=G$ and $dim(\rho)=2$ the subrepresentations
$Ind_G^{G'}(S^2(\rho\otimes \rho))$ and $Ind_G^{G'}(\Lambda^2(\rho\otimes \rho))$ and
${\bigoplus\bigoplus}_{\delta_1\neq \delta_2\in G'/G} \delta_1(V)\otimes \delta_2(V))$ of
$W\otimes W$ have dimension $>1$.

\bigskip\noindent
\begin{Lemma}\label{44}(Mackey) Suppose $(W,\rho)$ is an irreducible representation of $G$, and suppose
$A$ is a normal subgroup of $G$. Then either a) there exists a subgroup $H\subsetneqq G$
containing $A$ and an irreducible representation $\rho'$ of $H$ such that $\rho
=Ind_{H}^G(\rho')$, or b) the restriction of $\rho$ to $A$ is isotypic. In the latter case
there exists an irreducible representation $\rho'$ of $A$, which can be extended to a
projective representation of the group $G$ on the same vectorspace, and there exists an
irreducible projective representation $\rho''$ of $G/A$ such that $\rho\cong \rho'\otimes
\rho''$.
\end{Lemma}

\bigskip\noindent
\underbar{Proof}: Well known. % Similar to Serre Lineare Darstellung endlicher Gruppen, \S 9.1, theorem 16.

\bigskip\noindent
We say a representation $(W,\rho)$ of an algebraic group $G$ is almost faithful, if the kernel
of $\rho$ is a finite group.

\begin{Lemma}\label{remains}
For an almost faithful irreducible representation  $(W,\rho)$ of $G$ with the property
$m_G(W\otimes W)\leq 3$ the restriction of $(W,\rho)$ to $G^0$ remains irreducible unless $G^0$
is a torus.
\end{Lemma}

\bigskip\noindent
\underbar{Proof}: For $G^0\triangleleft G$ lemma \ref{44} implies, that either the restriction
of $\rho$ to $G^0$ is isotypic and $\rho=\rho'\otimes \rho''$, or $\rho=Ind_H^G(\rho')$ for
$G^0\subseteq H \subsetneqq G$. In the second case $G=G^0$ and $dim(W,\rho)=2$ by corollary
\ref{30}. In the first case lemma \ref{43} part 1 and 2 give $m_G(W\otimes W)\geq 4$ unless
either $dim(\rho)=1$ or $dim(\rho'')=1$.  If $dim(\rho)=1$, then $G^0$ is a torus, since $\rho$
is faithful. Hence $dim(\rho'')=1$, which proves the assertion.

\begin{Lemma}\label{46}
Under the assumptions made on the representation $(W,\rho)$ of $G$, the connected component
$G^0$ is a semisimple algebraic group of isotypic Lie-type (i.e. is isogenious to a product
$\prod_{i=1}^n H$ for a simple algebraic group)
\end{Lemma}

\bigskip\noindent
\underbar{Proof}: We replace the group $G$ by a group $G'$, which admits a surjective
homomorphism $G'\to G$ with finite kernel, and then we consider the almost faithful pullback
representation $(W,\rho)$ of $G'$.

\bigskip\noindent
\underbar{Construction of G'}: Notice $G^0=(G^0)_{der}\cdot Z(G^0)$ (with finite intersection).
The finite group $\pi_0(G)=G/G^0$ acts on $(G^0)_{der}$ and $Z(G^0)$. Using a splitting (see
\cite{Bo}) we get
$$ 1 \to Int({(G^0)}_{der}) \to Aut({(G^0)}_{der}) \to
Aut({(G^0)}_{der},B_{der},T_{der},\{x_\alpha\}_{\alpha\in \Delta}) \to 1
$$
$$ 1 \to Int({G^0}) \to Aut({G^0}) \to
Aut({G^0},B,T,\{x_\alpha\}_{\alpha\in \Delta}) \to 1 \ .$$ Notice $Int({(G^0)}_{der}) =
Int({G^0})$, and the natural map
$$ Aut({G^0}) \to Aut({(G^0)}_{der}) $$
induces a map
$$ Aut({G^0},B,T,\{x_\alpha\}_{\alpha\in \Delta}) \to Aut({(G^0)}_{der},B_{der},T_{der},\{x_\alpha\}_{\alpha\in
\Delta}) \ .$$ The group $Aut({(G^0)}_{der},B_{der},T_{der},\{x_\alpha\}_{\alpha\in \Delta})$
can be identified with a subgroup of the group of graph automorphisms of the Dynkin diagram of
${(G^0)}_{der}$. The group $G/G^0$ acts on $G$ by conjugation, which defines a map from $G/G^0$
to the automorphism group of $G^0$. The image of $G/G^0$ under
$$ G/G^0 \to Aut({(G^0)}_{der},B_{der},T_{der},\{x_\alpha\}_{\alpha\in \Delta}) \ $$
defines a group $\Delta \subseteq Aut({(G^0)}_{der},B_{der},T_{der},\{x_\alpha\}_{\alpha\in
\Delta})$ of automorphisms preserving a splitting. $\Delta$ acts by graph automorphisms on the
Dynkin diagram of ${(G^0)}_{der}$, and defines a finite subgroup of $G$. The kernel
$\Delta^0\subseteq G/G^0$ of this map is the subgroup, which acts trivially on the Dynkin
diagram, and it is contained in the centralizer of ${(G^0)}_{der}$ by the definition of
splitting \cite{Bo}. Hence $Z(G^0)\rtimes \Delta_0$ commutes with ${(G^0)}_{der}$. Hence we get
a surjective map from
$$ \Bigl({(G^0)}_{der}  \times
(Z(G^0)\rtimes \Delta_0)\Bigr) \rtimes \Delta \
$$
to $G$. So we end up in the situation of step 1, where $G_1={(G^0)}_{der}$ and $G_2 =
Z(G^0)\rtimes \Delta_0$. We can further break up $G_1$ into isotypic factors, for which each
term is isogenious to a product of simple groups of the same Lie-type. Or, with a more refined
decomposition of $G_1$  the factors are in 1-1 correspondence with the orbits of the action of
the group $\Delta$ on the set of connected components of the Dynkin graph of $G$. We may also
replace each simple factor by its simply connected covering group. This defines
$$  G' = (\prod_{i=0}^r G_i)\rtimes \Delta \ ,$$
where $(G_0)^0$ is a torus, where $(G_i)^0$ is isotypic semisimple for all $i=1,..,r$ such that
$\Delta$ normalizes each factor $G_i$ for $i=0,..,r$.

\bigskip\noindent
\underbar{Continuation of proof}: We can now apply lemma \ref{43}, which implies that the
restriction of $\rho$ to the subgroup $\prod_{i=0}^r G_i$ remains irreducible. Hence
$$ \rho\vert_{\prod_{i=0}^r G_i} \ \cong\ \bigotimes_{i=0}^r \ \rho_i $$
for irreducible representations $\rho_i$ of $G_i$. Since $\Delta$ normalizes each factor $G_i$
each representation $\rho_i$ can be extended to a projective representations $\rho_i^{ext}$ of
the semidirect product $G_i\rtimes \Delta$ respectively. Hence
$$ \rho \ \cong\ \bigotimes_{i=0}^r \rho_i^{ext} \ .$$
Since $m_{G'}(W\otimes W)\leq 3$ part 1 of lemma \ref{43} implies $dim(\rho_i^{ext})=1$ except
for one $i$. Since $\rho$ is almost faithful, therefore either $(G')^0$ is a torus, or $(G')^0$
is isotypic semisimple. Hence for the proof of the lemma it now suffices to exclude the case
that $G^0 =Z(G^0)=T$ is a torus. For this notice
$$ dim(Spec(\C[W]^G)) =  dim(Spec((\C[W]^{G^0})^{G/G^0}))
= dim(Spec(\C[W]^{G^0})) = 0 \ .$$ Hence $\C[W]^T=\C$. Since $\Lambda^i(W)$ for $i=dim(W)$ is
the trivial representation, there exist roots $\alpha_1,...,\alpha_i$ of $T$ on $W$, such that
$\sum_{\nu=1}^i \alpha_\nu =0$. Then $\C[W]^{T}$ has dimension $\geq 1$. A contradiction. This
proves lemma \ref{46}.

\bigskip\noindent
% Hence $G$ can be replaced, without restriction of generality, by $(G_1\times
%G_2)\rtimes \Delta$, where $G_1$ is isogenious to $H^n$ for a simple group $H$, and
%where $G_2=Z(G^0)\rtimes \Delta_0$. Since $\lambda_2$ is one dimensional and faithful,
%the group $G_2$ is cyclic abelian. $G_2$ is a cyclic group of order dividing $i$, since
%$\Lambda^i(W)$ is the trivial representation for $i=dim(W)$. This implies $(\lambda_2)^i=1$,
%since $H$ is simple and $G_2$ is cyclic abelian. Since $G_2$ commutes with $G_1$, we are now
%free to replace $H$ by its simply connected covering group, although the irreducible
%representation $W$ will then not be faithful any longer. Then
%$$  G^0 = H\times \cdots \times H $$
%for $H$ connected, simply, and simply connected. This proves the second and third claim.
%Furthermore $G_2=\Delta_0$ is a finite cyclic group of order dividing $i$ and
%$$ G = (H^n\times \Delta_0)  \rtimes \Delta $$
%where $\Delta$ acts on $\Delta_0$ and $H^n$ by conjugation. The action of $\Delta$ on $H^n$
%preserves a splitting, and acts transitively on the connected components of the Dynkin diagram.
%We now change notation. We may replace $\Delta$ by the semidirect product $\Delta_0
%\rtimes \Delta $,  such that for the new $\Delta$ we have
Hence we can assume
$$ G' = H^n \rtimes \Delta \ ,$$
where $\Delta\cong \pi_0(G')$ acts by conjugation on $H^n$ and acts transitively on the set of
the $n$ connected components of the Dynkin diagram of $G'$. Since the restriction of $\rho$ to
$(G')^0$ is irreducible by lemma \ref{remains}
$$ \rho\vert_{H^n} \cong \bigotimes_{i=1}^n \nu \quad , \quad dim(\nu)>1 $$
holds for an irreducible representation $\nu$ of the connected simply connected group simple
group $H$.
%The first component $H$ defines a point in the set of connected components of the Dynkin
%diagram of $G'$. Let $\Delta_1$ be the stabilizer in $\Delta$ of this point.
Via inclusions $i_\nu : H\to H^n$ and projections $\pi_\mu:H^n \to H$ each $\pi_\mu\circ
Int(g)\circ i_\nu$ for $g\in G'(k)$ defines an automorphism of $H$. Let ${\cal H}$ denote the
group generated by these automorphisms. It contains the subgroup $Int(H)$ as a normal subgroup
of finite index. By Schur's lemma for each $\sigma\in {\cal H}$ there exists  $A_\sigma\in
Gl(W_\nu)$ with $\nu(\sigma(h)) = A_\sigma \nu(h) A_{\sigma}^{-1}$, which is uniquely
determined by $\sigma$ up to a constant in $\overline\Q_l^*$. Thus the representation $\nu$ of
$H$ extends to a projective representation $\nu^{ext}$ of the group ${\cal H}$ on $W_\nu$. It
is clear that $g\mapsto Ad_g\vert_{ H^n}$ induces a map to the semidirect product
$$  \phi: G' \to {\cal H}^n \rtimes \Sigma_n  \ .$$
The representation $\nu^{ext}\otimes \cdots \otimes \nu^{ext}$ on $W_\nu \otimes \cdots \otimes
W_\nu$ extends to a projective representation $\rho'$ of ${\cal H}^n \rtimes \Sigma_n$ on which
the permutation groups acts on $(W_\nu)^{\otimes n}$ via the Weyl action, the group ${\cal H}$
acts via the matrices $A_\sigma$. From the Schur lemma it follows, that for $g\in G'(k)$ the
matrices $\rho(g)$ and $\rho'(\phi(g))$ conincide up to a scalar. The same holds for tensor
powers of $\rho$ respectively $\rho'\circ\phi$. Hence to obtain lower bounds for the
multiplicitities of the tensor powers $W^{\otimes r}$ of the representation $(W,\rho)$ we may
replace $(\rho,W)$ by $(\rho'\circ\phi,W)$ and $G'$ by the \lq{larger\rq}\ group ${\cal
H}^n\rtimes \Sigma_n$ in the following. Since $m_{{\cal H}^n\rtimes \Sigma_n}(W^{\otimes r})
\leq m_{G'}(W^{\otimes r})$ our assumption on the multiplicities for $ r=2$ and $r=3$ carry
over. For simplicity of notation we now write  $\Delta$ instead of $\Sigma_n$ and $\nu$ instead
of $\nu^{ext}$, and consider the action of $\Delta$ on
$$ W\otimes W = \bigotimes_{i=1}^n \ W_\nu\otimes W_\nu \ .$$

\bigskip\noindent
Our first aim is to show $n\leq 2$: Since the multiplicity of a representation is the same as
its multiplicity as a projective representation, in the following no distinctions will be made
between representations and projective representations. The projective representation
$W_\nu\otimes W_\nu$ is a projective representation of ${\cal H}$, which decomposes into at
least two invariant subspaces $A\oplus B$ by part two of lemma \ref{43}, since $dim(W_\nu)>1$.
This induces invariant subspaces in the tensorproduct $\rho\otimes \rho$. For $n=2$ these are
the subspaces $A\otimes A$ and $B\otimes B$ and $(A\otimes B)\oplus (B\otimes A)$. In general
there are at least $n+1$ different such nontrivial subspaces invariant under the group ${\cal
H}^n\rtimes \Delta$. This forces $n\leq 2$ by the bound $m_{{\cal H}^n\rtimes
\Sigma_n}(W\otimes W) \leq m_{G'}(W\otimes W)\leq 3$. Hence
\begin{enumerate}
\item either $n=2$ and $G={\cal H}^2\rtimes \Delta$ and $\Delta$
permutes the two factors of ${\cal H}^2$. In this case the multiplicity of $W\otimes W$ is 3.
\item or $n=1$ and $G^0$ is simple.
\end{enumerate}
In the first case we repeat this argument,  now with the third tensor power of $W$ instead of
the second tensor power. For $Y=W_\nu$ notice $Y\otimes Y\otimes Y = S^3(Y) \oplus 2 T^{2,1}
\oplus\Lambda^3(Y)$ decomposes into a direct sum of at least four different  subspaces
invariant under $H\rtimes\Delta_1$, and all four summands are nontrivial for $dim(Y)\geq 3$.
Since $dim(Y)\geq 2$ the last of the four subspaces vanishes only for $dim(Y)=2$. In this case
${\cal H}=Sl(2)$. So assume first $dim(Y)\geq 3$. Then the third tensor power $W\otimes
W\otimes W$ of the representation of $G={\cal H}^2 \rtimes \Delta$ admits a decomposition
$W^{\otimes 3}=(Y\otimes Y)^{\otimes 3}=(Y^{\otimes 3})\otimes (Y^{\otimes 3})$ into at least
${4\choose 2}=10$ different nontrivial ${\cal H}^2\rtimes \Sigma_2$-invariant constituents. By
assumption there are at most 7 irreducible constituents in $W\otimes W\otimes W$. This excludes
this case. It remains to consider the case $G'=G=(Sl(2)^2)\rtimes \Delta$ and $dim(Y)=2$. In
this case $\Delta$ is the cyclic group of order two, which permutes the two $Sl(2)$-factors.
The representation of $H=Sl(2)$ on $Y$ must be the standard representation. Hence $W\otimes W$
is induced from $st^{\otimes 2}\otimes st^{\otimes 2} = (S^2\oplus 1) \otimes (S^2\oplus 1)$.
But then, since $m_{G'}(W\otimes W)=3$, by our assumptions the trivial representation must be
contained in $\Lambda^2(W)$. This gives a contradiction, since now, evidently, the trivial
representation can only be contained in $1\otimes 1 \subseteq S^2(W)$. This proves step 4.

\bigskip\noindent
Concerning step 5. Hence
$$  G^0 = H $$
is a nontrivial connected simple algebraic group $H$ (which we may replace by its simply
connected covering). Our assumptions imply
$$dim(Spec(\C[W]^H))=dim(Spec(\C[W]^G))=0\ .$$ Since $H$ is
simple and $W$ is irreducible and since
$$ \C[W]^H = \C \ ,$$
the classification of Kac-Popov-Vinberg (see \cite{KPV}, \cite{S1}, \cite{S2}) shows, that
$(W,H)$ appears in the following list \begin{enumerate} \item $W$ is the standard
representation $st$ of the simple group $H$ or its dual (excluding the case $E_8$), or
\item the symmetric and alternating square of $A_n$. \item The second fundamental representation in $\Lambda^2(st)$
for $C_n$
\item the third fundamental representation of $C_3$ on $\Lambda^3(\C^6)$
or $\Lambda^3(st)$ for $A_5,A_6,A_7$. \item the spin
representation of $B_3, B_4, B_5, B_6$ \item the spin
representation of $D_4,D_5,D_6,D_7$.
\end{enumerate}
Notice, we can restrict ourselves to representations from \cite{S1}, 1.2 part (1), since
$\C[W]^H =\C$ implies that $0$ is the only closed orbit, hence $H$ is the principal isotropy
group in the sense of loc. cit. Also we omitted cases of loc. cit., where the dimension of
$\C[W]^H$ obviously is $>0$ as for example the adjoint representations, where this dimension of
is equal to the rank of $G$. In the following we freely use notation from \cite{S1}. For
representations $W$ of type $2)$ the third tensor power $W\otimes W\otimes W$ has more than 7
irreducible constituents by the Littlewood-Richardson rule. Since the representations of type 3
are obtained from the alternating square of the linear group as the highest weight constituents
of the restriction, one easily eliminates also this possibility (just compare highest weight
restrictions of the decomposition obtained above in the $A_n$-case). Concerning the
representations $W$ of type 4) it follows once more from the Littlewood-Richardson rules, that
$W=\Lambda^3(st)$ has a tensor square $W\otimes W$ with $4$ or more irreducible constituents
contrary to the assumption, that there should be at most 3 irreducible constituents. This
eliminates the representations of type 4). Similarly one can exclude the spin representations
$W$ of the types 5) and 6): The tensor square $W\otimes W$ of a spin representation is a sum of
the alternating powers $1+\Lambda +\cdots \Lambda^n$ in the $B_n$ case. Since $n\geq 3$ this
gives $>3$ irreducible constituents in $W\otimes W$ and therefore excludes the representations
in group 5). $W\otimes W= \Lambda^n_+ \oplus \Lambda^{n-2} \oplus \cdots $ holds for the tensor
square $W\otimes W$ of a spin representation $W$ of a group of $D_n$-type. This immediately
excludes $D_6,D_7$. Although for $D_4,D_5$ there are not more than three constituents, a closer
look at $W\otimes W$ exhibits the trivial representation to be contained in $S^2(W)$ and not in
$\Lambda^2(W)$. This also excludes type 6). Hence $W$ is the standard representation
$st=\varphi_1$ (or its dual) of a simple group (except $E_8$). Obviously the orthogonal cases
can be excluded by looking at the trivial constituent in $W\otimes W$. We also claim, that
neither of the standard representations $G_2,F_4,E_6,E_7$ are possible. For this we refer to
the decomposition rules given in the tables of \cite{S1}. For $G_2$ the third alternating power
and the second symmetric power of the standard representation $\varphi_1$ both contain the
trivial representation. For $F_4$ the second symmetric power of $\varphi_1$ contains the
trivial representation. For $E_6$ the trivial representation is contained in the third
symmetric power of $\varphi_1$. For $E_7$ the forth symmetric power of $\varphi_1$ contains the
trivial representation. By assumption $\Lambda^3(W)$ and $S^r(W)$ do not contain the trivial
representation ($\C[W]^G=0$ implies $\C[W]^{G^0}=\C$). This excludes all exceptional groups.

\bigskip\noindent
Since $dim(W)=2g-2$, therefore
\begin{enumerate} \item either $G^0=Sp(2g-2,\C)$ and $W$ is the selfdual standard
representation (if $W\otimes W$ has three irreducible constituents)
\item or $G^0=Sl(2g-2,\C)$ and $W$ is the standard representation or its dual (if $W\otimes W$ has two irreducible constituents) \end{enumerate}

\bigskip\noindent
In the second case, since $2g-2\geq 4$ the dual of the standard representation is not
isomorphic to the standard representation. Therefore by lemma \ref{remains} the homomorphism
$\pi_0(G) \to Out(G^0)$ must be trivial in both cases. Hence each element in $\pi_0(G)$ admits
a representative in $G$, which respects a splitting of $G^0$. Hence $G\cong G^0 \times
\pi_0(G)$. Since $\rho$ is faithful the group $\pi_0(G)$ must be a cyclic group of order
dividing $dim(W)=2g-2$. In case $G^0=Sp(2g-2)$ even stronger $\pi_0(G)$ has to be of order
$\leq 2$. In fact $m_G(W\otimes W)=m_{G^0}(W\otimes W)$ has three irreducible constituents in
this case. Hence by assumption the trivial representation of $G=G^0\times \pi_0(G)$ is
contained in $\Lambda^2(W)$. This forces $\pi_0(G)$ to act on $W$ faithfully by a quadratic
character. This proves the claim.

\bigskip\noindent
Since the center of $G^0$ is a cyclic group of order $2$ in the first case $G^0=Sp(2g-2)$ and a
cyclic group of order $2g-2$ in the second case $G^0=Sl(2g-2)$, the trivial representation of
$G^0$ only occurs in a tensor power $W^{\otimes r}$ for $2\vert r$ respectively $2g-2\vert r$
in the first respectively second case. But for these $r$ the finite cyclic group $\pi_0(G)$ has
acts trivially on $W^{\otimes r}$, since $\pi_0(G)$ is cyclic of order dividing $2$
respectively $2g-2$. Therefore a one dimensional representation of $G$, which is a constituent
of some tensor power $W^{\otimes r}$ of the generator $W$ of the tensor category $Rep(G)$, is
the trivial representation. This implies, that $\pi_0(G)$ is the trivial group.

\bigskip\noindent
Thus, tacitly assuming that the $\kappa\in X(k)$ has been normalized to become zero, we obtain

\bigskip\noindent
\begin{Theorem}\label{main} {\it For smooth projective curves $C$ of genus $g\geq 3$ the $\overline\Q_l$-linear tensor category ${\cal BN}$ of
BN-sheaves is equivalent to the category $Rep(G)$ of finite dimensional
$\overline\Q_l$-representations of the linear group $G=Sl(2g-2,\overline\Q_l)$, if $C$ is not
hyperelliptic, respectively the symplectic group $G=Sp(2g-2,\overline\Q_l)$, if $C$ is
hyperelliptic.} \end{Theorem}

\bigskip\noindent
\begin{Corollary} \label{IRR}{\it If $C$ is not hyperelliptic, the perverse sheaves ${}^p\delta_{\alpha}$ for $\alpha_1 <
2g-2$ are a full set of inequivalent representatives of the set of isomorphism classes of
irreducible perverse $BN$-sheaves.} \end{Corollary}

\bigskip\noindent
Immediate consequences are the theorem of Martens, the theorem of Torelli, and formulas for the
intersection cohomology groups of $W_r-W_r$ (see section \ref{ANEX}).

\bigskip\noindent
\underbar{Remark}: For $\kappa\neq 0$ the Tannaka group has the form $(G\times D)/Z$, where
$D=D(\kappa)$ is a diagonalizable group only depending on the order of $\kappa\in X(k)$, which
contains the center $Z=\mu_{2g-2}$ respectively $Z=\mu_2$ of the respective group $G$, that is
obtained (as above) after  $\kappa$ has been normalized to become zero for a suitable translate
$C$ of the given curve.

\bigskip\noindent
\underbar{The theta divisor}: Let  ${\cal T}(X,\Theta)$ for $X=J(C)$ be the Tannakian
subcategory of ${\cal BN}$ generated by $\delta_\Theta=\delta_{W_{g-1}}$. The generator
$\delta_\Theta$ of this tensor category corresponds to the representation $\Lambda^{g-1}(st)$
of $Sl(2g-2)$ respectively the highest weight constituent of $\Lambda^{g-1}(st)$ for
$Sp(2g-2)$. The subcategory ${\cal T}(\Theta)$ is therefore equivalent as a tensor category to
$$ {\cal T}(X,\Theta) \ \approx \ Rep(G') \ ,$$ where in the hyperelliptic case $G'=G$ or
$G'=G_{ad}$ depending on whether $g$ is even or odd. In the nonhyperelliptic case $
G'=G/\mu_{g-1}$ holds for $G=Sl(2g-2)$. In this latter case a perverse sheaf
${}^p\delta_\alpha$ is isomorphic to an object in ${\cal T}(X,\Theta)$ if and only if
$deg(\alpha)$ is divisible by $g-1$. For example $ \varepsilon_{g-1}$ is in ${\cal
T}(X,\Theta)$. In fact
$$ \varepsilon_{g-1} \hookrightarrow \delta_\Theta^{*(2g-1)} \ .$$
Conversely, in the non-hyperelliptic case, $ \delta_\Theta \hookrightarrow
\varepsilon_{g-1}^{*(2g-1)}$ holds. The adjoint representation $A$ always is in $Rep(G')$.
Since $A$ can be described intrinsically,  $A$ is  determined by the category ${\cal
T}(X,\Theta)$, hence determined by $(X,\Theta)$. Since $A$ is an irreducible perverse sheaf
with support $C-C$ (under our assumption $\kappa=0$) this determines $C$ from the data
$(X,\Theta)$ (Torelli). See \cite{W} for a more elementary version of this argument.

\goodbreak
\section{The categories ${\cal T}(X)$}\label{categor}

\bigskip\noindent
We say a projective (not necessarily smooth) curve $\tilde C\hookrightarrow X$ is generating,
if $\tilde C$ generates $X$ as an abelian variety. Theorem \ref{Geomorigin} allows to define a
tensor category ${\cal T}(X)$, generated by all admissible perverse sheaves supported on
generating curves of $X$ -- as a subcategory of the quotient category of the category of
semisimple perverse sheaves of geometric origin modulo the Serre subcategory of semisimple
translation-invariant perverse sheaves of geometric origin -- in a similar way as we did this
for a single smooth projective curve $C$ in its Jacobian.

\bigskip\noindent
For these more general tensor categories ${\cal T}(X)$ there are tensor functors $$f_*: {\cal
T}(Y)\to {\cal T}(X)$$ attached to any surjective homomorphism
$$ f: Y\to X $$
between abelian varieties defined over $k$, induced from the map $Rf_*: \overline
D_c^b(Y,\overline\Q_l)\to \overline D_c^b(X,\overline\Q_l)$. Furthermore there exist tensor
functors
$$ f_*: {\cal BN}(C) \to {\cal T}(X) $$
for any map $f: C\to X$ over $k$ from a smooth projective curve $C$ to the abelian variety $X$,
such that the image $\tilde C= f(C)$ generates $X$. Any of the generators $\delta_E$ of the
tensor category ${\cal T}(X)$ is a subobject of some $f_*(\delta_C')$. Since
$\Lambda^r(\delta_{C'})=0$ holds in ${\cal BN}(C')$ for some $r$ depending on $C'$, and since
the property $\Lambda^i(X)=0$ is inherited, when the tensor functor $f_*$ is applied to $X$, we
see again that the category ${\cal T}(X)$ admits a super fiber functor. Since the generators
$X_i,i\in I$ satisfy $\Lambda^{r_i}(X_i)=0$ for certain $r_i\in I$ the category is a Tannakian
category over the algebraically closed field $\overline\Q_l$, hence admits a fiber functor over
$\overline\Q_l$ by theorem \ref{Deligne} and its consequences.

\bigskip\noindent
\begin{Corollary} {\it Let $J(C)$ be the Jacobian of a smooth projective curve over $k$, and let
$f: J(C) \to X$ be a surjective homomorphism of abelian varieties over $k$. Then the direct
image $Rf_*({}^p\delta_\alpha) \in D_c^b(X,\overline\Q_l)$ of any of the perverse BN-sheaves
${}^p\delta_\alpha \in Perv(J(C))$ is of the form
$$ Rf_*({}^p\delta_\alpha) \ =\ P \oplus T \ $$
for some $P\in Perv(X)$ and a direct sum $T = \oplus_\nu T_\nu[\nu]$ of translates of
translation-invariant sheaves $T_\nu$ on $X$ both depending on $\alpha$.} \end{Corollary}

\bigskip\noindent
\underbar{Proof}: This follows from theorem \ref{Geomorigin} and the functoriality of
convolution, i.e the fact that
$$Rf_*: D_c^b(J(C),\overline\Q_l)\to D_c^b(X,\overline\Q_l)$$ is a $\overline\Q_l$-linear tensor functor.

\bigskip

\section{Coverings of curves}

\bigskip\noindent
Let $\varphi: C'\to C$ be a (possibly ramified) covering of projective curves over $k$, such
that $g(C)\geq 2$. We assume $f(P')=P$ for the distiguished points defining the Abel-Jacobi
map. Then the covering $\varphi$ induces a homomorphisms $f$ between abelian varieties
$$ f: X'=J(C') \longrightarrow X=J(C) \ .$$
Consider the associated $\overline\Q_l$-linear tensor categories with $\overline\Q_l$-linear
fiber functors

\begin{itemize}
\item $({\cal C}',\omega_{G'})$ for ${\cal C}'=Rep(G')$, where ${\cal C}'$ is the full subcategory generated by
$\delta_{C'}$ in ${\cal T}(X')$, and
\item $({\cal C},\omega_{G})$ for ${\cal C}=Rep(G)$, where ${\cal C}$ is the full subcategory generated by
$\delta_C$ in ${\cal T}(X)$. \end{itemize}

\bigskip\noindent
Let ${\cal C}''$ be the full $\overline\Q_l$-linear Tannakian tensor subcategory of ${\cal
T}(X)$ generated by $Rf_*({\delta_{C'}})$. Notice, $\delta_{C}$ is a direct summand of
$Rf_*(\delta_{C'})=R\varphi_*(\delta_{C'})$. So we have
 the fully faithful inclusion functor $v:{\cal C} \to {\cal C}''$. Let $\omega_{G}=
\omega_{G''}\circ v $ be the fiber functor induced by $\omega_{G''}$. $({\cal
C''},\omega_{G''})$ defines an algebraic group $G''=G''_\varphi$ over $\overline\Q_l$, such
that ${\cal C}'' = Rep(G'')$.

\bigskip\noindent
We get a commutative diagram of $\overline\Q_l$-linear tensor functors
$$ \xymatrix{ & & {\cal C}'' \ar@/^3mm/[rrd]^{\omega_{G''}} \ar@/_3mm/[lld]_{\omega_{G''}} & & \cr
Vec_{\overline\Q_l} & {\cal C}'\ar[ru]_u \ar[l]_-{\tilde\omega_{G'}} &  & {\cal C} \ar[lu]^{v}
\ar[r]^-{\omega_G} & Vec_{\overline\Q_l}\ , \cr }$$ where $\tilde\omega_{G'}$ is defined by
$\tilde\omega_{G'}= \omega_{G''}\circ u$, for $u=Rf_*\vert_{{\cal C}'}$.

\bigskip\noindent
Notice, that $\tilde\omega_{G'}$ again is a $\overline\Q_l$-linear fiber functor of ${\cal
C}'$. This follows from the fact, that any exact tensor functor $F: {\cal C}'\to {\cal C}''$
between rigid, abelian $\overline\Q_l$-linear tensor categories ${\cal C}', {\cal C}''$ with
$End(1')=\overline\Q_l$ and $1''\neq 0$ is automatically faithful (\cite{DM}, prop. 1.19).
Since $u$ is an exact $\overline\Q_l$-linear functor, the functor $\tilde\omega_{G'}$ again is
a fiber functor of ${\cal C}'$. By \cite{Sav}, prop. 3.1.1.1 and \cite{Mi}, prop.III.4.2 and
cor.III.4.7 the equivalence of the category of fiber functors (of the neutral Tannaka category
${\cal C}'$) with the category of $G'$-torsors over $\overline\Q_l$ (\cite{DM}, Thm 3.2)
implies, that $\omega_{G'}$ and $\tilde\omega_{G'}$ are isomorphic fiber functors
$$ \omega_{G'} \approx \tilde\omega_{G'}\ .$$
We therefore assume $\omega_{G'}=\tilde\omega_G$ in the following.

\bigskip\noindent
By \cite{DM}, cor. 2.9 for any tensor functor $u: {\cal C}'=Rep(G')\to {\cal C}''=Rep(G'')$,
such that $\omega_{G''}\circ u = \omega_{G'}$, there exists a unique group homomorphism $\mu:
G''\to G'$ which induces $u$. This applies to $u$ in the situation above, and similarly applies
to the tensor functor $v: {\cal C}=Reg(G)\to {\cal C}''=Rep(G'')$:

\bigskip\noindent
This gives rise to a correspondence $(G'',\nu,\mu)$ between $G'$ and $G$ attached to the
covering $\varphi: C'\to C$ via homomorphisms $\mu$ and $\nu$ of algebraic groups over
$\overline\Q_l$
$$  \xymatrix{ & G''_\varphi\ar@{_{(}->}[dl]_\mu \ar@{->>}[dr]^\nu & \cr G' & & G \ .\cr } $$
By \cite{Sav} 4.3.2 and \cite{DM}, prop. 2.2.1 we have

\begin{Lemma} The homomorphism $\mu$ is a closed immersion, and the homomorphism $\nu$
is a faithfully flat epimorphism. \end{Lemma}

\bigskip\noindent
\underbar{Proof}: For the first assertion it is enough by loc. cit., that every object of
${\cal C}''$ is isomorphic to a subquotient of an object $u(X')$ for $X'\in {\cal C}'$. But
this holds by the definition of ${\cal C}''$. For the second assertion it is enough by loc.
cit., that $v$ is fully faithful and that every subobject of $v(X)$, $X\in {\cal C}$ is
isomorphic to the image of a subobject of $X$. Again this holds by definition.

\bigskip\noindent
If $C'\to C$ is a covering of smooth projective curves, then the groups $G,G'$ depend only on
the genus $g'=g(C'),g=g(C)$, the Riemann constants and whether $C',C$ are hyperelliptic or not.
This gives a diagram defining representations $\rho_\mu$ and $\rho_\nu$
$$ \xymatrix{ & & G''_\varphi(\overline\Q_l)\ar@{_{(}->}[dl]_\mu \ar@{->>}[dr]^\nu \ar@/^3mm/[rrd]^{\rho_\nu} \ar@/_3mm/[lld]_{\rho_\mu} & & \cr Gl(2g'-2,\overline\Q_l) & G'(\overline\Q_l) \ar@{_{(}->}[l] & & G(\overline\Q_l) \ar@{^{(}->}[r] & Gl(2g-2,\overline\Q_l)\ . \cr } $$
We remark, that in the case of a Galois covering $C'\to C$ this commutes with the Galois
action. Indeed this gives an underlying refined structure to the Artin representation \cite{Se}
of the Galois group.

\bigskip\noindent
\goodbreak
\section{Another Application}\label{ANEX}

\bigskip\noindent
Assume $C$ is a projective smooth non-hyperelliptic curve, and assume $\kappa=0$ (using a
translate of $C$). For $1\leq r\leq s\leq g-1$ consider the difference map $f(x,y)=x-y$ on
$X=J(C)$
$$ f: W_r \times W_s \to X \ .$$
By the Littlewood-Richardson rules  the direct image $ Rf_*\bigl(\delta_{W_r\times W_s}\bigr)$
is  $ \delta_s^- * \delta_r\ =\ \sum_{\nu=0}^s \delta_{\chi-s+\nu,r-\nu}\  =\ \delta_{\chi-s,r}
+
 \delta_{s-1}^- * \delta_{r-1} $. Hence $$Rf_* (\delta_{W_r\times W_s})\ = \ \delta_{\chi-s,r}\ \oplus \ Rf_*
(\delta_{W_{r-1}\times W_{s-1}})\ .$$  Recall the sheaf complex $\delta_s^-
* \delta_r= Rf_*(\delta_{W_r\times W_s})$ is perverse for $r+s<g$ by corollary \ref{pin}.
%Consider $\sigma = g-1 -s$ and
%$\rho = g-1-r$. Then the pairs $(r,s)$ with $r+s\geq g-1$
%correspond to the pairs $(\rho,\sigma)$ with $\rho +\sigma \leq
%g-1$. Notice, that $\delta_\sigma$ is a translate of
Now $\delta_{\chi-s,r}={}^p \delta_{\chi-s,r}\oplus T(r,s)$, where  $T(r,s)$ are translates of
constant sheaves on $X$. By corollary \ref{IRR} the perverse sheaf  ${}^p \delta_{\chi-s,r}$ is
irreducible. This implies

\begin{Lemma} If $C$ is not hyperelliptic, then ${}^p\delta_{2g-2-r,r} = \delta_{W_r - W_r}$ for $ r < g/2$.
\end{Lemma}

\bigskip\noindent
Thus the expression for  the cohomology given by lemma \ref{dimension} combined with the
results of section \ref{Compu} gives an  explicit description of the intersection cohomology of
the singular subvariety $W_r-W_r \subseteq X$.

\bigskip\noindent
\underbar{Proof}: In fact by induction, since $\delta_{W_r-W_r}$ is a nonconstant irreducible
direct summand of $Rf_*(\delta_{W_r\times W_r})$ and not a summand of $Rf_*(\delta_{W_s\times
W_s})$, it must be a summand of ${}^p\delta_{\chi-r,r}$. Notice  its support has dimension $2r$
and therefore is not contained in $W_s-W_s$ for $s<r$. This follows from the well known
dimension formula $dim(W_s-W_s)=min(2s,g)$ \cite{ACGH}, p.223. Since ${}^p\delta_{\chi-r,r}$ is
irreducible by corollary \ref{IRR} the equality ${}^p\delta_{\chi-r,r} \cong\delta_{W_r - W_r}$
follows.

\end{document}